\documentclass[11pt]{article}
\usepackage{}

\usepackage{amsmath}
\usepackage{amsfonts}
\usepackage{graphicx}
\usepackage{amssymb}
\usepackage{bm}
\usepackage{epstopdf}
\usepackage{color}
\usepackage{dsfont}

\newtheorem{condition**}{A*}
\newtheorem{condition***}{C*}
\newtheorem{condition*}{C}

\newtheorem{definition}{Definition}[section]
\newtheorem{theorem}{Theorem}[section]
\newtheorem{lemma}{Lemma}[section]
\newtheorem{remark}{Remark}[section]
\newcommand{\ignore}[1]{}

\newenvironment{keywords}{{\bf Key words: }}{}

\begin{document}

\title{Linear-Quadratic Delayed Mean-Field Social Optimization}
\author{Tianyang Nie\thanks{Tianyang Nie is with the School of Mathematics, Shandong University, Jinan, Shandong 250100, China (e-mail: nietianyang@sdu.edu.cn). }\qquad
Shujun Wang \thanks{Corresponding author. Shujun Wang is with the School of Management, Shandong University, Jinan, Shandong 250100, China (e-mail: wangshujun@sdu.edu.cn).}
\qquad Zhen Wu
\thanks{Zhen Wu is with the School of Mathematics, Shandong University, Jinan, Shandong 250100, China (e-mail: wuzhen@sdu.edu.cn).}}

\date{}

\maketitle

\begin{abstract}
A linear quadratic (LQ) stochastic optimization problem with delay involving weakly-coupled large population is investigated in this paper. Different to classic mean field (MF) game, here agents cooperate with each other to minimize the so-called \emph{social} objective. With the aid of \emph{delayed person-by-person optimality} principle, one arrives at an auxiliary LQ delayed control problem by decentralized information. A decentralized strategy is obtained by feat of an MF type anticipated forward-backward stochastic differential delay equation (AFBSDDE) consistency condition. The discounting method with delay feature is employed to solve the consistency condition system. Finally, by some estimates of AFBSDDEs we derive the asymptotic social optimality.
\end{abstract}

\section{Introduction}

In this section, firstly we give some notations and introduce main motivations of this work. Then a \emph{stochastic linear quadratic delayed MF game problem} is posed, which will be investigated in this paper.

\subsection{Notation}

Consider a finite time horizon $[0,T]$ for fixed $T>0$. Assume that $(\Omega,\mathcal F,$ $\{\mathcal F_t\}_{0\leq t\leq T},\mathbb P)$ is a complete filtered probability space satisfying the usual conditions and $\{W_i(t),1\leq i\leq N\}_{0\leq t\leq T}$ is a $N\times 1$-dimensional Brownian motion on this space.  Let $\mathcal F_t$ be the filtration generated by $\{W_i(s),1\leq i\leq N\}_{0\leq s\leq t}$ and augmented by $\mathcal N_{\mathbb P}$ (the class of all $\mathbb P$-null sets of $\mathcal F$). Let $\mathcal F_t^i$ be the augmentation of $\sigma\{W_i(s),0\leq s\leq t\}$ by $\mathcal N_{\mathbb P}$. Let $\delta>0$ and $\theta>0$ be two given fixed time horizon. $\mathbb E$ denotes the expectation under $\mathbb P$ and $\mathbb E^{\mathcal F_t}[\cdot]:=\mathbb E[\cdot|\mathcal F_t]$ denotes the conditional expectation. Define $\mathcal F_t:=\mathcal F_0$ for $\forall\ t\in[-(\delta\vee \theta),0)$.

Let $\langle\cdot,\cdot\rangle$ denote standard Euclidean inner product. $x^\top$ denotes the transpose of a vector (or matrix) $x$. For a matrix $M\in \mathbb R^{n\times m}$, we define the norm $|M|:=\sqrt{\texttt{tr}(M^\top M)}$. $M\in \mathbb{S}^n$ denotes the set of symmetric $n\times n$ matrices with real elements. $M> (\geq) 0$ denotes that $M\in \mathbb{S}^n$ which is positive (semi)definite, while $M\gg 0$ denotes that, for some $\varepsilon>0$, $M - \varepsilon I \geq 0$. For any times $t,\ t_1,\ t_2$ valued in $[-(\delta\vee \theta),T]$, denote by
\begin{itemize}
  \item $L_{\mathcal F}^2(\Omega;\mathbb R^n):=\Big\{\zeta:\Omega\rightarrow\mathbb R^n|\zeta\text{ is }\mathcal F\text{-measurable such that }\mathbb E|\zeta|^2<\infty\Big\}$
   \item $L^2_{\mathcal {F}}(\Omega;C([t_1,t_2];\mathbb R^n)):=\Big\{\zeta(\cdot):[t_1,t_2]\times\Omega\rightarrow\mathbb R^n|\zeta(\cdot)$ \text{is} $\mathcal F_t$\text{-adapted, continuous, such that } \\ $\mathbb E \Big[\sup\limits_{s\in[t_1,t_2]}|\zeta(s)|^2\Big]<\infty\Big\}$
  \item $L_{\mathcal F}^2(t_1,t_2;\mathbb R^n):=\Big\{\zeta(\cdot):[t_1,t_2]\times\Omega\rightarrow\mathbb R^n|\zeta(\cdot)\text{ is }\mathcal F_t\text{-progressively measurable}$ $\text{process such that}\\
      \mathbb E\int_{t_1}^{t_2}|\zeta(\cdot)|^2dt<\infty\Big\}$
  \item $L^2(t_1,t_2;\mathbb R^n):=\Big\{\zeta(\cdot):[t_1,t_2]\rightarrow\mathbb R^n|\zeta(\cdot)\text{ is }\text{deterministic  function s.t.}$ $ \int_{t_1}^{t_2}|\zeta(\cdot)|^2dt<\infty\Big\}$
  \item $L^\infty(t_1,t_2;\mathbb R^{n\times m}):=\Big\{\zeta(\cdot):[t_1,t_2]\rightarrow\mathbb R^{n\times m}|\zeta(\cdot)\text{ is uniformly bounded}\Big\}$
 \end{itemize}

\subsection{Motivation}

Consider a controlled large population system in which the dynamic of agent $\mathcal A_i$ is modelled by a stochastic differential equation (SDE)
\begin{equation}\left\{\label{eq00}\begin{aligned}
dX_i(t)=\ &b\Big(t,X_i(t),X^{(N)}(t),u_i(t),u^{(N)}(t)\Big)dt+g\Big(t,X_i(t),X^{(N)}(t),u_i(t),u^{(N)}(t)\Big)dW_i(t),\\
X_i(0)=\ &x_{i0}
\end{aligned}\right.\end{equation}
with cost functional
\begin{equation}\label{eq01}\begin{aligned}
\mathcal J_i(x_{i0},u(\cdot))=\ &\mathbb E\Bigg\{\int_0^TL\Big(t,X_i(t),X^{(N)}(t),u_i(t),u^{(N)}(t)\Big)dt+\Phi\big(x_{i0},X_i(T),X^{(N)}(T)\big)\Bigg\},
\end{aligned}\end{equation}
where $X^{(N)}(\cdot)=\frac{1}{N}\sum_{i=1}^NX_i(\cdot)$ and $u^{(N)}(\cdot)=\frac{1}{N}\sum_{i=1}^Nu_i(\cdot)$ denote the average state and average strategy of agents, respectively; $u(\cdot)=(u_1(\cdot)$, $\cdots,u_N(\cdot))$, $u_i(\cdot)\in \mathcal U_i:=\Big\{u_i(\cdot)|u_i(\cdot)\in L^2_{\mathcal F}(0,T;\mathbb R^m)\Big\}$, $i=1,\cdots,N$; $x_{i0}$ is deterministic.
Then we can pose a classical MF optimal control problem.\\

\textbf{Problem 00.} Find a strategy $\bar u=(\bar u_1,\cdots,\bar u_N)$ where $\bar u_i(\cdot)\in \mathcal U_i$, $1\leq i\leq N$ such that
\begin{equation}\nonumber
\mathcal J_i(x_{i0},\bar u(\cdot))=\inf_{u_{i}\in\mathcal U_{i}}\mathcal J_i(x_{i0},\bar u_1(\cdot),\cdots,u_i(\cdot),\cdots,\bar u_{N}(\cdot)).
\end{equation}

The controlled large population systems are widely involved in many fields, such as economics, engineering, biology, physics, etc. It is remarkable that because of the highly complicated coupling structure, it is infeasible and ineffective to obtain the classical \emph{centralized} strategies by combining all agents' exact dynamics in the large population problem. As a substitute, it is more effective and tractable to investigate the corresponding \emph{decentralized} strategies by taking into account his/her own individual dynamic (information) and some off-line quantities only. In this research direction, many researchers have taken a lot of efforts on MF game study.  The interested readers may refer to \cite{BSYS2016,BLPR2017,CD2013,HHN2018,LL2007,NNY2020,NC2013} and the references therein for MF game theory.

We further define
$$\mathcal J_{soc}^{(N)}(u(\cdot))=\sum_{i=1}^N\mathcal J_i(x_{i0},u(\cdot))$$
as the aggregate functional of $N$ agents corresponding to \eqref{eq00}-\eqref{eq01}. Then one poses a classical cooperative-type MF optimal control problem.\\

\textbf{Problem 01.} Find a strategy $\bar u=(\bar u_1,\cdots,\bar u_N)$ where $\bar u_i(\cdot)\in \mathcal U_i$, $1\leq i\leq N$ such that
\begin{equation}\nonumber\mathcal J_{soc}^{(N)}(\bar u(\cdot))=\inf_{u_{i}\in\mathcal U_{i},1\leq i\leq N}\mathcal J_{soc}^{(N)}(u_1(\cdot),\cdots,u_i(\cdot),\cdots,u_{N}(\cdot)).
\end{equation}
Contrast to aforementioned works for classical MF game where the agents are competitive, cooperative team optimization problem has attracted a lot of attentions in last ten years, which is so-called \emph{social optima problem}. One can refer to \cite{AM2015,DW22,FHW2021,HWY2021,HCM2012,HN2019,WZZ2020,WZ2017} for corresponding research in different frameworks. For more researches and applications of MF social optima problems, readers can refer to \cite{BDT2019,CBM2016,NM2018} and corresponding references.

To proceed, both competitive behaviors or cooperative trades are influenced by time factor. For example, investors or producers' decisions are affected by macro economic policies. There might occur some time delay due to the policy transmission mechanism. In many phenomena, the evolution of a controlled system depends not only on the present state or decision policy at time $t$, but also on their past trajectories on some history interval $[t-\theta,t]$, where $\theta$ denotes some lagged index. Correspondingly, a controlled system may contain \emph{output delay} (state delay) or \emph{input delay} (control delay).
In particular, there are various practical problems, in which the systems evolved by some stochastic differential delay equations (SDDEs).
To be specific, we formulate the dynamic $X(\cdot)$ with \emph{output} delay $\delta$ and \emph{input} delay $\theta$ as
\begin{equation}\left\{\label{eq02}\begin{aligned}
dX_i(t)=\ &b\Big(t,X_i(t),X_i(t-\delta),X^{(N)}(t-\delta),u_i(t),u_i(t-\theta),u^{(N)}(t-\theta)\Big)dt\\
&\qquad+g\Big(t,X_i(t),X_i(t-\delta),X^{(N)}(t-\delta),u_i(t),u_i(t-\theta),u^{(N)}(t-\theta)\Big)dW_i(t),\\
X_i(0)=\ &x_{i0},\ X_i(t)=X^0(t),\ t\in[-\delta,0),\ u_i(t)=u^0(t),\ t\in[-\theta,0)
\end{aligned}\right.\end{equation}
with cost functional
\begin{equation}\label{eq03}\begin{aligned}
\mathcal J_i(x_{i0},u(\cdot))=\ &\mathbb E\Bigg\{\int_0^TL\Big(t,X_i(t),X_i(t-\delta),X^{(N)}(t-\delta),u_i(t),u_i(t-\theta),u^{(N)}(t-\theta)\Big)dt\\
&\qquad+\Phi\big(x_{i0},X_i(T)\big)\Bigg\}.
\end{aligned}\end{equation}

Nevertheless, the delayed responses are challenging and difficult to be processed because of absence of It\^o's formula to deal with the delayed part of the trajectory. In recent years, the study of SDDEs has attracted the attention of many researchers, see e.g. \cite{Mao97}, \cite{OS00}, \cite{OS08}, etc.
It is well known that there is a perfect duality between SDEs and backward SDEs (BSDE). Peng and Yang \cite{PY09} introduced a new type of BSDEs called anticipated BSDEs (ABSDEs) and established a duality between SDDEs and ABSDEs. By virtue of the theory of ABSDEs and dual method, rich literature studies emerge. More details of SDDE theory and its wide applications may be found in the monograph of \cite{CW10}, \cite{HL2018}, \cite{MP15},\cite{Yu12}, \cite{ZXS21} and the references therein.

\subsection{Problem Formulation}\label{formulation}

Based on the above discussion, with consideration to obtain some explicit results, in this paper we turn to study an LQ large population system with $K$-type \emph{discrete} heterogeneous agents $\{\mathcal A_i:1\leq i\leq N\}$. The dynamics of the agents are given by a system of linear SDDEs with mean-field coupling:
that is, for $1\leq i\leq N,$
\begin{equation}\left\{\label{eq1}\begin{aligned}
dx_i(t)=\ &\Big[A_{\theta_i}(t)x_i(t)+\hat{A}_{\theta_i}(t)x_i(t-\delta)+\tilde A(t)x^{(N)}(t-\delta)+B(t)u_i(t)+\hat{B}(t)u_i(t-\theta)\\
&+\tilde B(t)u^{(N)}(t-\theta)\Big]dt+\Big[D(t)u_i(t)+\hat D(t)u_i(t-\theta)\Big]dW_i(t),\ t\in[0,T],\\
x_i(0)=\ &\xi_i,\ x_i(t)=x^0(t),\ t\in[-\delta,0),\ u_i(t)=u^0(t),\ t\in[-\theta,0),
\end{aligned}\right.\end{equation}
where $x^{(N)}(\cdot)=\frac{1}{N}\sum_{i=1}^Nx_i(\cdot)$ denotes the forward state-average of the agents. Note that while the coefficients $(A_{\theta_i}(\cdot),\hat{A}_{\theta_i}(\cdot),\tilde A(\cdot),B(\cdot),\hat B(\cdot),\tilde B(\cdot),D(\cdot),\hat D(\cdot))$ depend on the time variable $t$, in what follows the variable $t$ will usually be suppressed if no confusion would occur. The number $\theta_i$ is a parameter of the agent $\mathcal A_i$ to model a heterogeneous population. For simplicity, we only assume that the coefficients $A$ to be dependent on $\theta_i$. Similar analysis can be proceeded in case that all other coefficients are also dependent on $\theta_i$. Moreover, we assume that $\theta_i$ takes values in a finite set $\Theta:=\{1,2,\cdots,K\}$. We call $\mathcal A_i$ a $k$-type agent if $\theta_i=k\in\Theta$. In this paper, we are interested in the asymptotic behavior as $N$ tends to infinity. For $1 \leq k \leq K$,  introduce$$\mathcal{I}_k=\{i|\theta_i=k, 1 \leq i \leq  N\}, \quad \quad N_k=|\mathcal{I}_k|,$$where $N_k$ is the cardinality of index set $\mathcal{I}_k$. For $1\leq k\leq K$, let $\pi_k^{(N)}=\frac{N_k}{N}$, then $\pi^{(N)}=(\pi_1^{(N)}, \cdots, \pi_K^{(N)})$ is a probability vector representing the empirical distribution of $\theta_1, \cdots, \theta_N.$ We introduce the following assumption:
\begin{description}
  \item[(A1)] There exists a probability mass vector $\pi=(\pi_1, \cdots, \pi_K)$ such that $\displaystyle{\lim_{N\rightarrow+\infty}}\pi^{(N)}=\pi$, and $\displaystyle{\min_{1 \leq k \leq K}}\pi_{k}>0.$
  \item[(A2)] For $i=1,\cdots,N$, $\xi_i\in L_{\mathcal F}^2(\Omega;\mathbb R^n)$, $x^0\in L^2(-\delta,0;\mathbb R^n)$, $u^0\in L^2(-\theta,0;\mathbb R^d)$. If $\theta_i=\theta_j=k$,
   $\xi_i$ and $\xi_j$ are independent identically distributed (i.i.d) and they have the same distribution as the random variable $\xi^{(k)}$.
   \item[(A3)] $A_{\theta_i}\in L^\infty(0,T;\mathbb R^{n\times n})$, $\hat{A}_{\theta_i},\tilde A\in L^\infty(0,T+\delta;\mathbb R^{n\times n})$, $B,D\in L^\infty(-\theta,T;\mathbb R^{n\times d})$, $\hat B,\tilde{B},\hat{D}\in L^\infty(0,T+\theta;\mathbb R^{n\times d})$, $i=1,\cdots,N$.
  \end{description}
It follows that under (A1)-(A3), the state equation \eqref{eq1} admits a unique solution for all $u_i$.
In fact, if we denote by
\begin{equation}\left\{\nonumber\begin{aligned}
&\mathbb{X}=\left(
    \begin{smallmatrix}
      x_1 \\
      \vdots \\
      x_N \\
    \end{smallmatrix}
  \right),
\mathbb{U}=\left(
    \begin{smallmatrix}
      u_1 \\
      \vdots \\
      u_N \\
    \end{smallmatrix}
  \right), \mathbb{W}=\left(
    \begin{smallmatrix}
      W_1 \\
      \vdots \\
      W_N \\
    \end{smallmatrix}
  \right),\hat{\xi}=\left(
    \begin{smallmatrix}
      \xi_1 \\
      \vdots \\
      \xi_N \\
    \end{smallmatrix}
  \right),
\mathds{1}=\left(
    \begin{smallmatrix}
      1 \\
      \vdots \\
      1 \\
    \end{smallmatrix}
  \right),\mathbb A=\left(
              \begin{smallmatrix}
                A_{\theta_1} &  &  \\
                 & \ddots &  \\
                 &  & A_{\theta_N} \\
              \end{smallmatrix}
            \right),
\\
&\hat{\mathbb A}=\left(
              \begin{smallmatrix}
                \hat A_{\theta_1} &  &  \\
                 & \ddots &  \\
                 &  & \hat A_{\theta_N} \\
              \end{smallmatrix}
            \right),\tilde{\mathbb A}=\frac{1}{N}\left(
              \begin{smallmatrix}
                \tilde A & \cdots & \tilde A  \\
                \vdots &  & \vdots \\
                \tilde A  & \cdots & \tilde A  \\
              \end{smallmatrix}
            \right),
\mathbb B=\left(
              \begin{smallmatrix}
                B &  &  \\
                 & \ddots &  \\
                 &  & B \\
              \end{smallmatrix}
            \right),
\hat{\mathbb B}=\left(
              \begin{smallmatrix}
                \hat B &  &  \\
                 & \ddots &  \\
                 &  & \hat B \\
              \end{smallmatrix}
            \right),
\tilde{\mathbb B}=\frac{1}{N}\left(
              \begin{smallmatrix}
                \tilde B & \cdots & \tilde B  \\
                \vdots &  & \vdots \\
                \tilde B  & \cdots & \tilde B  \\
              \end{smallmatrix}
            \right),\\
&\mathbb D(\mathbb{U})=\left(
              \begin{smallmatrix}
                Du_1 &  &  \\
                 & \ddots &  \\
                 &  & Du_N \\
              \end{smallmatrix}
            \right),
\hat{\mathbb D}(\mathbb{U}_\theta)(\cdot)=\left(
              \begin{smallmatrix}
               \hat D(\cdot)u_1(\cdot-\theta) &  &  \\
                 & \ddots &  \\
                 &  & \hat D(\cdot)u_N(\cdot-\theta) \\
              \end{smallmatrix}
            \right),
\end{aligned}\right.\end{equation}
then \eqref{eq1} can be rewritten as
\begin{equation}\label{eq2}\left\{\begin{aligned}
d\mathbb X(t)=\ &\Big[\mathbb A(t)\mathbb X(t)+\big(\hat{\mathbb A}(t)+\tilde{\mathbb A}(t)\big)\mathbb X(t-\delta)+\mathbb {B}(t)\mathbb {U}(t)+\big(\hat{\mathbb B}(t)+\tilde{\mathbb {B}}(t)\big)\mathbb U(t-\theta)\Big]dt\\
&+\Big(\mathbb {D}(\mathbb U)(t)+\hat{\mathbb D}(\mathbb U_\theta)(t)\Big) d\mathbb {W}(t),\ t\in[0,T],\\
\mathbb {X}(0)=\ &\hat{\xi},\quad \mathbb X(t)=x^0(t)\mathds{1},\ t\in[-\delta,0),\ \mathbb U(t)=u^0(t)\mathds{1},\ t\in[-\theta,0),
\end{aligned}\right.\end{equation}
which is a linear SDDE of vector value.
%Here, ``$\circ$" denotes the Hadamard product. It is well known that Hadamard product (also called Schur product or entry-wise product) is a binary operation between two matrices of the same dimensions, and it may produce another matrix in which each element $(i,j)$ is the product of the elements $(i,j)$ in the original matrices.
By the classic theory of SDDE (see e.g. \cite{Mao97}), \eqref{eq2} admits a unique solution $\mathbb X(\cdot)\in L^{2}_{\mathcal{F}}(\Omega;C([-\delta,T];\mathbb R^{Nn}))$ for $\mathbb U(\cdot)\in L^{2}_{\mathcal{F}}(-\theta, T; \mathbb{R}^{Nd})$.
Thus, for any $1\leq i\leq N$, the state equation \eqref{eq1} admits a unique solution $x_i(\cdot)\in L^2_{\mathcal {F}}(\Omega;C([-\delta,T];\mathbb R^n))$.

Let $u=(u_1,\cdots,u_N)$ be the set of strategies of all $N$ agents and  $u_{-i}=(u_1,\cdots,u_{i-1},u_{i+1}$, $\cdots,u_N)$, $1\leq i\leq N$. The cost functional for $\mathcal A_i$, $1\leq i\leq N$, is given by
\begin{equation}\label{eq3}\begin{aligned}
\mathcal J_i(u_i(\cdot),u_{-i}(\cdot))=\ &\frac{1}{2}\mathbb E\Big\{\int_0^T\big[\big\langle Q(t)(x_i(t)-S(t) x^{(N)}(t)),x_i(t)-S (t) x^{(N)}(t)\big\rangle\\
&\quad+\big\langle \tilde{Q}(t)(x_i(t-\delta)-\tilde{S}(t) x^{(N)}(t-\delta)),x_i(t-\delta)-\tilde{S} (t) x^{(N)}(t-\delta)\big\rangle\\
&\quad+\big\langle R_{\theta_i}(t)u_i(t),u_i(t)\big\rangle+\big\langle \tilde{R}_{\theta_i}(t)u_i(t-\theta),u_i(t-\theta)\big\rangle\big]dt\\
&\quad+\big\langle G (x_i(T)-\Gamma x^{(N)}(T)),x_i(T)-\Gamma x^{(N)}(T)\big\rangle\Big\}.
\end{aligned}\end{equation}
The aggregate team functional of $N$ agents is
\begin{equation}\label{eq4}
\mathcal J_{soc}^{(N)}(u(\cdot))=\sum_{i=1}^N\mathcal J_i(u_i(\cdot),u_{-i}(\cdot)).
\end{equation}

We impose the following assumptions on the coefficients of the cost functionals:
\begin{description}
\item[(A4)] $Q(\cdot)\in L^\infty(0,T;\mathbb S^n)$, $\tilde{Q}(\cdot)\in L^\infty(0,T+\delta;\mathbb S^n)$, $Q(\cdot)\geq0,\ \tilde{Q}(\cdot)\geq0$, $\tilde{Q}(t)=0,\ t\in(T,T+\delta]$, $S(\cdot)\in L^\infty(0,T;\mathbb R^{n\times n})$, $\tilde{S}(\cdot)\in L^\infty(0,T+\delta;\mathbb R^{n\times n})$, $\tilde{S}(t)=0,\ t\in(T,T+\delta]$, $R_{\theta_i}(\cdot)\in L^\infty(-\theta,T;\mathbb S^d)$, $\tilde{R}_{\theta_i}(\cdot)\in L^\infty(0,T+\theta;\mathbb S^d)$, $R_{\theta_i}(\cdot)\gg 0,\ \tilde{R}_{\theta_i}(\cdot)\gg 0$, $\tilde{R}_{\theta_i}(t)=0,\ t\in(T,T+\theta]$ ($i=1,\cdots,N$), $G \in \mathbb S^{n\times n}$, $G\geq0$, $\Gamma\in \mathbb R^{n\times n}$.
\end{description}
For $i=1,\cdots,N$, the centralized admissible strategy set for the $i^{th}$ agent is given by
$$\mathcal U_i^c=\Big\{u_i(\cdot)|u_i(\cdot)\in L^2_{\mathcal F}(0,T;\mathbb R^d)\Big\}.$$
Correspondingly, the decentralized admissible strategy set for the $i^{th}$ agent is given by
$$\mathcal U_i^d=\Big\{u_i(\cdot)|u_i(\cdot)\in L^2_{\mathcal F^i}(0,T;\mathbb R^d)\Big\}.$$
We propose the following optimal problem:\\

\textbf{Problem 1.} Find a strategy set $\bar u=(\bar u_1,\cdots,\bar u_N)$ where $\bar u_i(\cdot)\in \mathcal U_i^c$, $1\leq i\leq N$, such that
\begin{equation}\label{eq5}\mathcal J_{soc}^{(N)}(\bar u(\cdot))=\inf_{u_{i}\in\mathcal U_{i}^c,1\leq i\leq N}\mathcal J_{soc}^{(N)}(u_1(\cdot),\cdots,u_i(\cdot),\cdots,u_{N}(\cdot)).
\end{equation}

\begin{definition}
A strategy $\widetilde u_i(\cdot)\in\mathcal U_i^d$, $i=1,\cdots,N$ is an $\varepsilon$-social decentralized optimal strategy if there exists $\varepsilon=\varepsilon(N)>0$,
$\lim_{N\rightarrow\infty}\varepsilon(N)=0$ such that
$$\frac{1}{N}\Big(\mathcal J_{soc}^{(N)}
(\widetilde u(\cdot))-\inf_{u_{i}(\cdot)\in\mathcal U_{i}^c,1\leq i\leq N}\mathcal J_{soc}^{(N)}(u(\cdot))\Big)\leq\varepsilon.$$
\end{definition}
\begin{remark}
In reality, the linear SDDE \eqref{eq1} stands for the dynamics of some investment behaviors such as in stocks and bonds in a self-financed market with some delay effects. In the problems such as finance, optimal control, etc., the SDDE dynamics have been deeply studied in the existing literature, such as \cite{CW10}, \cite{Mao97}, \cite{OS00} and so on. In addition, in \eqref{eq1} $x^{(N)}(\cdot)$ and $u^{(N)}(\cdot)$ appear as the delay items, the mainly consideration of which is due to the practical significance. We think the collective behaviours and integrated environment should produce delayed effects to individual agent.
\end{remark}

\begin{remark}
In \eqref{eq1}, $x_i,$ $x^{(N)}(\cdot)$ and $u^{(N)}(\cdot)$ do not enter in the diffusion term of the dynamics. The reason is that there exists essential difficulties while doing the error estimations, because if $x_i,$ $x^{(N)}(\cdot)$ or $u^{(N)}(\cdot)$ enters in the diffusion term, we have to deal with some error estimations like $\sum_{i=1}^N z_i(\cdot)-\mathbb E z_i(\cdot)$, where $z_i(\cdot)$ stands for a part of solution of BSDE derived by the martingale representation theorem. It seems an impossible task based on the existing BSDE theory. In the future, one may focus on this problem and try to derive some new technique to overcome this difficulty.
\end{remark}

\begin{figure}[!h]
\centering
\setlength{\abovecaptionskip}{-0.1cm}
\hspace*{0.4cm}\includegraphics[width=6in,height=2.2in]{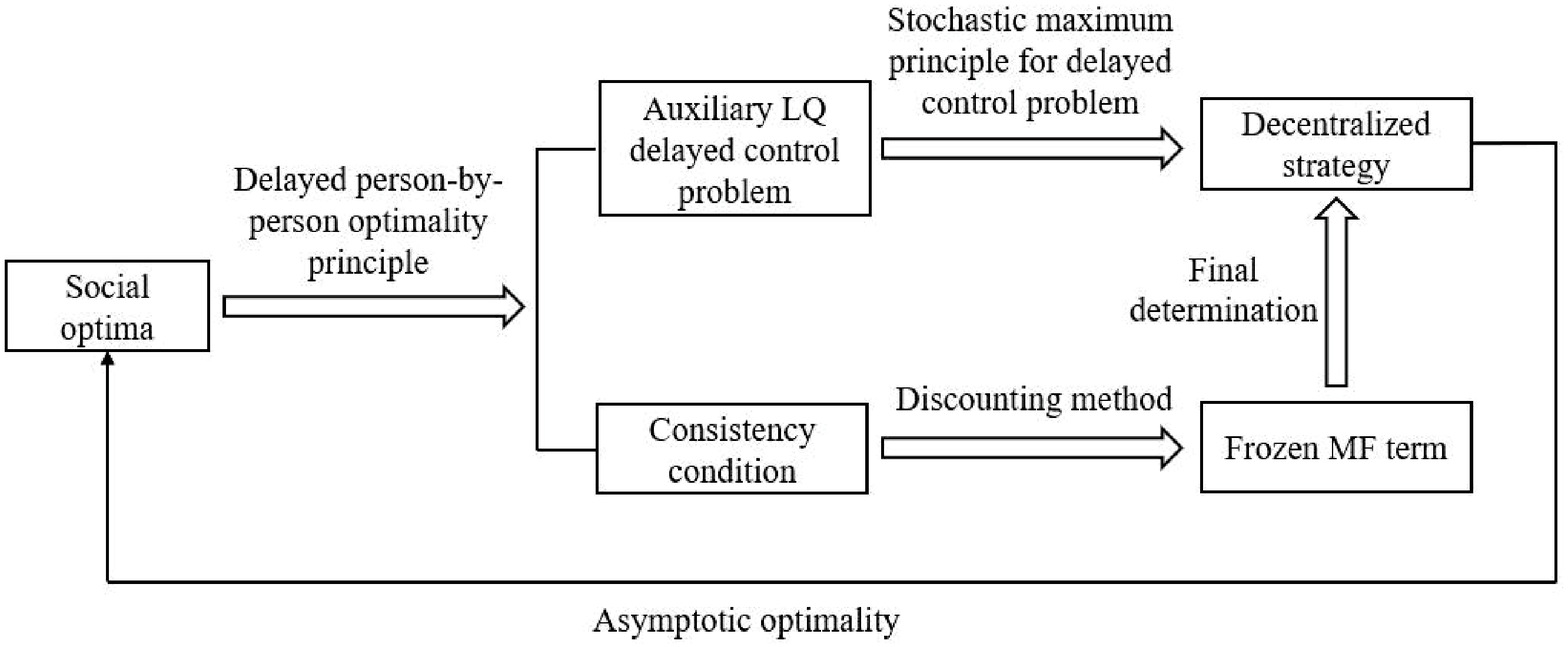}
\caption{The research roadmap.}\label{f1}
\end{figure}

Now we briefly show the process of studying \textbf{Problem 1} (see Figure \ref{f1}):
\begin{itemize}
 \item Firstly, with the help of \emph{delayed person-by-person optimality} principle and variational technique, we obtain an auxiliary LQ control problem with delay. Then stochastic maximum principle (\cite{CW10}, etc) are applied to derive the optimal control of auxiliary problem.
 \item Secondly, we establish and investigate the consistency condition (CC) system by discounting method (e.g. \cite{HHN2018}, \cite{PT1999}) to determine such frozen MF term and some off-line variables.
 \item Thirdly, by virtue of standard estimations of AFBSDDE (\cite{CW10}, \cite{Yu12}), we verify the decentralized strategy is the asymptotic optimality of centralized strategy.
\end{itemize}

The contributions of this work are illustrated as follows.
\begin{enumerate}
  \item In this work, we use indicative functions to overcome the gaps in \cite{CW10, ZXS21}, etc. In fact, for delay problems the value of adjoint equations (BSDEs indeed) are always supposed to be zero after terminal time $T$ (e.g. \cite{CW10, ZXS21}), i.e.  $Y(t)=0$ for $t\in[T,T+\delta]$, which will make solution $Y$ be noncontinuous at time $T$. Here we adopt indicative functions $I_{[0,T-\delta]}$ to guarantee the continuity of $Y$ since we do not give any assumptions for $Y$ after terminal time $T$.
  \item The monotonic condition and Riccati equation method are invalid due to the delay structure, the discounting method is applied to derive the wellposedness of CC system which is recognized as a MF-AFBSDDE. Besides, there exists the circular dependency of parameter $\rho$ and we overcome the difficulty.
  \item We give a new description about directional derivative of social cost, which is more clearer and accurate than that of existing works.
\end{enumerate}

The rest of this paper is organized as follows. In Section 2, we apply person-by-person optimality to derive the auxiliary control problem with delay of the individual agent. The CC system and its well-posedness are established by discounting method in Section 3. In Section 4, the asymptotic optimality of decentralized strategy is obtained. Section 5 concludes this paper.

\section{Stochastic delayed optimal control problem for $\mathcal A_i$}\label{auxiliary problem}

In this section, we try to solve the optimal delayed control problem and derive the decentralized control.

\subsection{Person-by-person optimality}\label{p-b-p optimality}

In MF social optima scheme (MF team), person-by-person optimality is a critical technique, which has been used in the recent social optima literature, e.g. \cite{WZ2017}, etc. There is significant difference between MF team scheme and MF game scheme, where the auxiliary control problem is usually derived directly by fixing the state-average. This would lead to some ineffective control in social optima scheme. Thus, in this section under the delayed person-by-person optimality principle, variation method will be applied to analyze the MF approximation.

Let $\{\bar u_i,\bar u_{-i}\in \mathcal U_i^c\}_{i=1}^N$ be
 centralized optimal strategy of all the agents. Now consider the perturbation that the agent $\mathcal A_i$ use the strategy $u_i\in \mathcal U_i^c$ and all the other agents still apply the strategy $\bar u_{-i}=(\bar u_1,\cdots,\bar u_{i-1},\bar u_{i+1},\cdots,\bar u_N)$. The realized states \eqref{eq1} corresponding to $u:=(u_i,\bar u_{-i})$ and $\bar u:=(\bar u_i,\bar u_{-i})$ are denoted by $(x_1,\cdots,x_N)$ and $(\bar x_1,\cdots,\bar x_N)$, respectively.
 For $j=1,\cdots,N$, denote the perturbation
$$\delta u_j=u_j-\bar u_j,\quad \delta x_j=x_j-\bar x_j,\quad\delta x^{(N)}=x^{(N)}-\bar{x}^{(N)},\quad \delta u^{(N)}=u^{(N)}-\bar{u}^{(N)}.$$
Therefore, the variation of the state for $\mathcal A_i$ is given by
\begin{equation}\left\{\label{eq6}\begin{aligned}
&d\delta x_i=\ \Big[A_{\theta_i}\delta x_i+\hat A_{\theta_i}\delta x_i(t-\delta)+\tilde A\delta x^{(N)}(t-\delta)+B\delta u_i+
\hat B\delta u_i(t-\theta)\\
&\qquad\qquad\quad+\frac{\tilde B}{N}\delta u_i(t-\theta)\Big]dt+\Big[D\delta u_i(t)+\hat D\delta u_i(t-\theta)\Big]dW_i(t),\ t\in[0,T],\\
&\delta x_i(0)=\ 0,\ \delta x_i(t)=0,\ t\in[-\delta,0),\ \delta u_i(t)=0,\ t\in[-\theta,0),
\end{aligned}\right.\end{equation}
and for $\mathcal A_j$, $j\neq i$,
\begin{equation}\label{eq7}\left\{\begin{aligned}
d\delta x_j=\ &\Big[A_{\theta_j}\delta x_j+\hat A_{\theta_j}\delta x_j(t-\delta)+\tilde A\delta x^{(N)}(t-\delta)+\frac{\tilde B}{N}\delta u_i(t-\theta)\Big]dt,\ t\in[0,T],\\
\delta x_j(0)=\ &0,\ \delta x_j(t)=0,\ t\in[-\delta,0),\ \delta u_i(t)=0,\ t\in[-\theta,0].
\end{aligned}\right.\end{equation}
Here, we have used the fact $\delta u^{(N)}(t-\theta)=\frac{1}{N}\delta u_i(t-\theta)$.
For $k=1,\cdots,K$, define $\delta x_{(k)}=\sum_{j\in\mathcal I_k,j\neq i}\delta x_j$, thus we have
\begin{equation}\label{eq8}\left\{\begin{aligned}
&d\delta x_{(k)}=\ \Big[A_{k}\delta x_{(k)}+\hat A_{k}\delta x_{(k)}(t-\delta)+(N_k-I_{\mathcal I_k}(i))\tilde{A}\delta x^{(N)}(t-\delta)\\
&\qquad\qquad+(N_k-I_{\mathcal I_k}(i))\frac{\tilde{B}}{N}\delta u_i(t-\theta)\Big]dt,\ t\in[0,T],\\
&\delta x_{(k)}(0)=\ 0,\ \delta x_{(k)}(t)=0,\ t\in[-\delta,0), \ \delta u_i(t)=0,\ t\in[-\theta,0],
\end{aligned}\right.\end{equation}
where $I_{\mathcal I_k}(\cdot)$ denotes the indicative function, that is $I_{\mathcal I_k}(i)=\left\{\begin{array}{cc}
                                                                                                1, & i\in \mathcal I_k\\
                                                                                                0, & i\notin \mathcal I_k
                                                                                              \end{array}\right.
$.
By some elementary calculations, we can further obtain the variation of the cost functional of $\mathcal A_i$ as follows
\begin{equation*}\begin{aligned}
\Delta \mathcal J_{i}
=\ &\mathbb E\Big\{\int_0^T\big[\big\langle Q(\bar x_i-S \bar x^{(N)}),\delta x_i-S \delta x^{(N)}\big\rangle+\big\langle \tilde{Q}(\bar x_i(t-\delta)-\tilde{S} \bar x^{(N)}(t-\delta)),\\
\quad&\ \ \delta x_i(t-\delta)-\tilde{S} \delta x^{(N)}(t-\delta)\big\rangle+\big\langle R_{\theta_i}\bar u_i,\delta u_i\big\rangle+\big\langle \tilde{R}_{\theta_i}\bar u_i(t-\theta),\delta u_i(t-\theta)\big\rangle\big]dt\\
\quad&+\big\langle G(\bar x_i(T)-\Gamma \bar x^{(N)}(T)),\delta x_i(T)-\Gamma \delta x^{(N)}(T)\big\rangle\Big\}.
\end{aligned}\end{equation*}
For $j\neq i$, the variation of the cost functional of $\mathcal A_j$ is given by
\begin{equation*}\begin{aligned}
\Delta \mathcal J_{j}=\ &\mathbb E\Big\{\int_0^T\big[\big\langle Q(\bar x_j-S \bar x^{(N)}),\delta x_j-S \delta x^{(N)}\big\rangle+\big\langle \tilde{Q}(\bar x_j(t-\delta)-\tilde{S} \bar x^{(N)}(t-\delta)),\delta x_j(t-\delta)\\
&\quad-\tilde{S} \delta x^{(N)}(t-\delta)\big\rangle\big]dt+\big\langle G(\bar x_j(T)-\Gamma \bar x^{(N)}(T)),\delta x_j(T)-\Gamma \delta x^{(N)}(T)\big\rangle\Big\}.
\end{aligned}\end{equation*}
Therefore, by combining above equalities, the variation of the social cost satisfies
\begin{equation}\label{eq10}\begin{aligned}
\Delta \mathcal J_{soc}^{(N)}=\ &\mathbb E\Big\{\int_0^T\Big[\big\langle Q(\bar x_i-S \bar x^{(N)}),\delta x_i-S \delta x^{(N)}\big\rangle+\sum_{j\neq i}\big\langle Q(\bar x_j-S\bar x^{(N)}), \delta x_j-S\delta x^{(N)}\big\rangle\\
&+\big\langle \tilde{Q}(\bar x_i(t-\delta)-\tilde{S} \bar x^{(N)}(t-\delta)),\delta x_i(t-\delta)-\tilde{S} \delta x^{(N)}(t-\delta)\big\rangle\\
&+\sum_{j\neq i}\big\langle \tilde{Q}(\bar x_j(t-\delta)-\tilde{S} \bar x^{(N)}(t-\delta)),\delta x_j(t-\delta)-\tilde{S} \delta x^{(N)}(t-\delta)\big\rangle\\
&+\big\langle R_{\theta_i}\bar u_i,\delta u_i\big\rangle+\big\langle \tilde{R}_{\theta_i}\bar u_i(t-\theta),\delta u_i(t-\theta)\big\rangle \Big]dt+\big\langle G(\bar x_i(T)-\Gamma \bar x^{(N)}(T)),\\
&\ \ \delta x_i(T)-\Gamma \delta x^{(N)}(T)\big\rangle+\sum_{j\neq i}\big\langle G(\bar x_j(T)-\Gamma \bar x^{(N)}(T)),\delta x_j(T)-\Gamma \delta x^{(N)}(T)\big\rangle\Big\}.
\end{aligned}\end{equation}

\textbf{Step 1:}
First, replacing $\bar x^{(N)}$ and $\bar x^{(N)}(t-\delta)$ in \eqref{eq10} by some MF terms $\hat x$ and $\hat x(t-\delta)$, respectively, which will be determined later,
 \begin{equation*}\begin{aligned}
\Delta \mathcal J_{soc}^{(N)}
=\ &\mathbb E\Big\{\int_0^T\Big[\big\langle Q\bar x_i,\delta x_i\big\rangle-\big\langle (QS+ S^\top Q-S^\top QS)\hat x,\delta x_i\big\rangle\\
&-\sum_{k=1}^K\big\langle (QS+ S^\top Q -S^\top QS)\hat x,\delta x_{(k)}\big\rangle+ \sum_{k=1}^K\frac{1}{N_k}\sum_{j\in\mathcal I_k,j\neq i}\big\langle Q \bar x_j,N_k\delta x_{j}\big\rangle\\
&+\big\langle \tilde{Q}\bar x_i(t-\delta),\delta x_i(t-\delta)\big\rangle-\big\langle (\tilde{Q}\tilde{S}+ \tilde{S}^\top \tilde{Q}-\tilde{S}^\top \tilde{Q}\tilde{S})\hat x(t-\delta),\delta x_i(t-\delta)\big\rangle\\
&-\sum_{k=1}^K\big\langle (\tilde{Q}\tilde{S}+ \tilde{S}^\top \tilde{Q} -\tilde{S}^\top \tilde{Q}\tilde{S})\hat x(t-\delta),\delta x_{(k)}(t-\delta)\big\rangle\\
&+ \sum_{k=1}^K\frac{1}{N_k}\sum_{j\in\mathcal I_k,j\neq i}\big\langle \tilde{Q} \bar x_j(t-\delta),N_k\delta x_{j}(t-\delta)\big\rangle\\
&+\big\langle R_{\theta_i}\bar u_i,\delta u_i\big\rangle+\big\langle \tilde{R}_{\theta_i}\bar u_i(t-\theta),\delta u_i(t-\theta)\big\rangle \Big]dt+\big\langle G\bar x_i(T),\delta x_i(T)\big\rangle\\
&-\big\langle (G\Gamma+\Gamma ^\top G-\Gamma ^\top G\Gamma)\hat x(T),\delta x_i(T)\big\rangle-\sum_{k=1}^K\big\langle (G\Gamma+\Gamma ^\top G-\Gamma ^\top G\Gamma)\hat x(T),\delta x_{(k)}(T)\big\rangle\\
&+ \sum_{k=1}^K\frac{1}{N_k}\sum_{j\in\mathcal I_k,j\neq i}\big\langle G \bar x_j(T),N_k\delta x_{j}(T)\big\rangle\Big\}+\sum_{l=1}^3\varepsilon_l,
\end{aligned}\end{equation*}
where
\begin{equation*}\left\{\begin{aligned}
&\varepsilon_1=\mathbb E\int_0^T\big\langle (QS+ S^\top Q- S^\top Q S)(\hat x-\bar x^{(N)}),N\delta x^{(N)}\big\rangle dt,\\
&\varepsilon_2=\mathbb E\int_0^T\big\langle (\tilde{Q}\tilde{S}+ \tilde{S}^\top \tilde{Q} -\tilde{S}^\top \tilde{Q}\tilde{S})(\hat x(t-\delta)-\bar x^{(N)}(t-\delta)),N\delta x^{(N)}(t-\delta)\big\rangle dt,\\
&\varepsilon_3=\mathbb E\big\langle (G\Gamma+\Gamma ^\top G-\Gamma ^\top G\Gamma)(\hat x(T)-\bar x^{(N)}(T)),N\delta x^{(N)}(T)\big\rangle .
\end{aligned}\right.\end{equation*}

\textbf{Step 2:}
Next, for $k=1,\cdots,K$, introduce the limit
 $x^{**}_k$ to replace $\delta x_{(k)}$, and for $j\in\mathcal I_k$, $j\neq i$, introduce the limit $x_j^*$ to replace $N_k\delta x_j$, where
\begin{equation}\label{eq11}\left\{\begin{aligned}
&dx_k^{**}=\Big[A_kx_k^{**}+\hat A_kx_k^{**}(t-\delta)+\tilde{A}\pi_k\delta x_i(t-\delta)+\tilde{A}\pi_k\sum_{l=1}^Kx_l^{**}(t-\delta)+\tilde{B}\pi_k\delta u_i(t-\theta)\Big]dt,\\
&dx_j^*=\Big[A_kx_j^*+\hat A_kx_j^*(t-\delta)+\tilde{A}\pi_k\delta x_i(t-\delta)+\tilde{A}\pi_k\sum_{l=1}^Kx_l^{**}(t-\delta)+\tilde{B}\pi_k\delta u_i(t-\theta)\Big]dt,\ t\in[0,T],\\
&x_k^{**}(0)=0,\ x_j^*(0)=0,\ x_k^{**}(t)=0,\ x_j^*(t)=0,\ t\in[-\delta,0).
\end{aligned}\right.\end{equation}
Therefore,
\begin{equation}\label{eq12}\begin{aligned}
\Delta \mathcal J_{soc}^{(N)}
  =&\mathbb E\Bigg\{\int_0^T\Big[\langle Q\bar x_i,\delta x_i\rangle-\langle (QS+ S^\top Q -S^\top QS)\hat x,\delta x_i\rangle
-\sum_{k=1}^K\langle (QS+ S^\top Q\\
&\qquad-S^\top QS)\hat x,x^{**}_k\rangle+ \sum_{k=1}^K\frac{1}{N_k}\sum_{j\in\mathcal I_k,j\neq i}\langle Q \bar x_j, x_j^* \rangle+\langle \tilde{Q}\bar x_i(t-\delta),\delta x_i(t-\delta)\rangle\\
&-\langle (\tilde{Q}\tilde{S}+ \tilde{S}^\top \tilde{Q} -\tilde{S}^\top \tilde{Q}\tilde{S})\hat x(t-\delta),\delta x_i(t-\delta)\rangle\\
  &
-\sum_{k=1}^K\langle (\tilde{Q}\tilde{S}+ \tilde{S}^\top \tilde{Q} -\tilde{S}^\top \tilde{Q}\tilde{S})\hat x(t-\delta),x^{**}_k(t-\delta)\rangle\\
&+ \sum_{k=1}^K\frac{1}{N_k}\sum_{j\in\mathcal I_k,j\neq i}\langle \tilde{Q} \bar x_j(t-\delta), x_j^*(t-\delta) \rangle+\langle R_{\theta_i}\bar u_i,\delta u_i\rangle+\langle \tilde{R}_{\theta_i}\bar u_i(t-\theta),\\
&\qquad\delta u_i(t-\theta)\rangle \Big]dt+\langle G\bar x_i(T),\delta x_i(T)\rangle-\big\langle (G\Gamma+\Gamma ^\top G-\Gamma ^\top G\Gamma)\hat x(T),\delta x_i(T)\big\rangle\\
&-\sum_{k=1}^K\big\langle (G\Gamma+\Gamma ^\top G-\Gamma ^\top G\Gamma)\hat x(T),x_k^{**}(T)\big\rangle+ \sum_{k=1}^K\frac{1}{N_k}\sum_{j\in\mathcal I_k,j\neq i}\langle G \bar x_j(T),x^*_{j}(T)\rangle\Bigg\}\\
&+\sum_{l=1}^9\varepsilon_l,
\end{aligned}\end{equation}
where
\small\begin{equation*}\left\{\begin{aligned}
&\varepsilon_4=\sum_{k=1}^K\mathbb E\int_0^T\langle (QS+ S^\top Q-S^\top QS)\hat x,x_k^{**}-\delta x_{(k)}\rangle dt,\\
&\varepsilon_5=\sum_{k=1}^K\mathbb E\int_0^T\frac{1}{N_k}\sum_{j\in\mathcal I_k,j\neq i}\langle Q\bar x_j,N_k\delta x_j-x^*_j\rangle dt,\\
&\varepsilon_6=\sum_{k=1}^K\mathbb E\int_0^T\langle (\tilde{Q}\tilde{S}+ \tilde{S}^\top \tilde{Q} -\tilde{S}^\top \tilde{Q}\tilde{S})\hat x(t-\delta),x_k^{**}(t-\delta)-\delta x_{(k)}(t-\delta)\rangle dt,\\
&\varepsilon_7=\sum_{k=1}^K\mathbb E\int_0^T\frac{1}{N_k}\sum_{j\in\mathcal I_k,j\neq i}\langle \tilde{Q}\bar x_j(t-\delta),N_k\delta x_j(t-\delta)-x^*_j(t-\delta)\rangle dt,\\
&\varepsilon_8=\sum_{k=1}^K\mathbb E\langle (G\Gamma+\Gamma ^\top G-\Gamma ^\top G\Gamma)\hat x(T),x_k^{**}(T)-\delta x_{(k)}(T)\rangle,\\
&\varepsilon_9=\sum_{k=1}^K\frac{1}{N_k}\sum_{j\in\mathcal I_k,j\neq i}\langle G\bar x_j(T),N_k\delta x_j(T)-x^*_j(T)\rangle.
\end{aligned}\right.\end{equation*}\normalsize

\textbf{Step 3:}
Finally, we will substitute $x_j^*$ and $x_k^{**}$ by dual method. It is very important to construct an auxiliary control problem for investigating decentralized control in social optimal problem (see e.g. \cite{HCM2012}, \cite{WZ2017}). We may use a duality procedure to break away $\delta \mathcal J_{soc}^{(N)}$ from the dependence on $x_j^*$ and $x_k^{**}$. To this end, we introduce the following adjoint equations $(y_1^j,z_1^j)$ and $(y_2^k,z_2^k)$
of the terms $x_j^*$ and $x_k^{**}$, respectively, which are shown as follows
\begin{equation}\nonumber\left\{\begin{aligned}
&dy_1^j=\alpha_1 dt+z_1^jdW_j(t)+\sum_{l=1,l\neq j}^N z_{1}^{jl}dW_l(t),\\
&dy_2^k=\alpha_2dt+z_2^kdW_k(t),\ t\in[0,T],\\
&y_1^j(T)=G\bar x_j(T),\ y_2^k(T)=-(G\Gamma+\Gamma ^\top G-\Gamma ^\top G\Gamma)\hat x(T),
%\\
%&y_1^j(t)=0,\ z_1^j(t)=0,\ y_2^k(t)=0,\ z_2^k(t)=0,\ t\in(T,T+\delta],
\ j=1,\cdots,N,\ k=1,\cdots,K.
\end{aligned}\right.\end{equation}
Applying It\^{o}'s formula to $\langle y_1^j,x_j^*\rangle$, we have
\begin{equation*}\begin{aligned}
d\langle y_1^j,x_j^*\rangle=\ &\Big[\langle y_1^j,A_kx_j^*+\hat A_kx_j^*(t-\delta)+\tilde{A}\pi_k\delta x_i(t-\delta)+\tilde{A}\pi_k\sum_{l=1}^Kx_l^{**}(t-\delta)+\tilde{B}\pi_k\delta u_i(t-\theta)\rangle\\
&+\langle \alpha_1,x_j^*\rangle\Big]dt+\sum_{j=1}^N(\cdots)dW_j(t).
\end{aligned}\end{equation*}
For $j\in\mathcal I_k$, integrating from $0$ to $T$ and taking expectation, we obtain
\begin{equation}\label{eq13}\begin{aligned}
&\mathbb E\langle G\bar x_j(T), x^*_{j}(T)\rangle\\
=\ &\mathbb E\langle y_1^j(T),x_j^*(T)\rangle-\mathbb E\langle y_1^j(0),x_j^*(0)\rangle\\
=\ &\mathbb E\int_0^T\Big[\langle y_1^j,A_kx_j^*+\hat A_kx_j^*(t-\delta)+\tilde{A}\pi_k\delta x_i(t-\delta)+\tilde{A}\pi_k\sum_{l=1}^Kx_l^{**}(t-\delta)+\tilde{B}\pi_k\delta u_i(t-\theta)\rangle\\
\ &\qquad+\langle \alpha_1,x_j^*\rangle\Big]dt\\
=\ &\mathbb E\int_0^T\Big[\langle\alpha_1+A_k^\top y_1^j,x_j^*\rangle+\langle \hat A_k^\top y^j_1,x_j^*(t-\delta)\rangle+\langle \pi_k\tilde{A}^\top y_1^j,\delta x_i(t-\delta)\rangle\\
&\qquad+\langle \pi_k\tilde{B}^\top y_1^j,\delta u_i(t-\theta)\rangle+\sum_{l=1}^K\langle \pi_k\tilde{A}^\top y_1^j,x_l^{**}(t-\delta)\rangle\Big]dt.
\end{aligned}\end{equation}
Similarly, we have
\begin{equation}\label{eq15}\begin{aligned}
&-\mathbb E\langle (G\Gamma+\Gamma ^\top G-\Gamma ^\top G\Gamma)\hat x(T),x_k^{**}(T)\rangle=\mathbb E\langle y_2^k(T),x_k^{**}(T)\rangle-\mathbb E\langle y_2^k(0),x_k^{**}(0)\rangle\\
=\ &\mathbb E\int_0^T\Big[\langle\alpha_2+A_k^\top y_2^k,x_k^{**}\rangle+\langle\hat A_k^\top y_2^k,x_k^{**}(t-\delta)\rangle+\langle \pi_k\tilde{A}^\top y_2^k,\delta x_i(t-\delta)\rangle\\
&\qquad+\langle \pi_k\tilde{B}^\top y_2^k,\delta u_i(t-\theta)\rangle+\sum_{l=1}^K\langle \pi_k\tilde{A}^\top y_2^k,x_l^{**}(t-\delta)\rangle\Big]dt.
\end{aligned}\end{equation}
Now we should point out the fact that (recall that $x_j^{\ast}=0$, for $t\in[-\delta,0]$)
\begin{equation*}\begin{aligned}
\mathbb E\int_0^T\langle A_k^\top y_1^j,x_j^*(t-\delta)\rangle dt&=\mathbb E\int_0^{T-\delta}\langle \hat A_k^\top(t+\delta) \mathbb E^{\mathcal F_t}[y_1^j(t+\delta)],x_j^*(t)\rangle dt\\
&=\mathbb E\int_0^{T}\langle \hat A_k^\top(t+\delta) \mathbb E^{\mathcal F_t}[y_1^j(t+\delta)]I_{[0,T-\delta]},x_j^*(t)\rangle dt.
\end{aligned}\end{equation*}
Similarly, (recall that $\delta u_i(t)=0$, for $t\in[-\theta,0]$)
\begin{equation*}\begin{aligned}
&\mathbb E\int_0^T\langle \pi_k\tilde{A}^\top y_1^j,\delta x_i(t-\delta)\rangle dt=\mathbb E\int_0^T\langle \pi_k\tilde{A}^\top(t+\delta) \mathbb E^{\mathcal F_t}[y_1^j(t+\delta)]I_{[0,T-\delta]},\delta x_i(t)\rangle dt,\\
&\mathbb E\int_0^T\langle \pi_k\tilde{B}^\top y_1^j,\delta u_i(t-\theta)\rangle dt=\mathbb E\int_0^T\langle \pi_k\tilde{B}^\top(t+\theta) \mathbb E^{\mathcal F_t}[y_1^j(t+\theta)]I_{[0,T-\theta]},\delta u_i(t)\rangle dt,\\
&\mathbb E\int_0^T\langle \pi_k\tilde{A}^\top y_1^j,x_l^{**}(t-\delta)\rangle dt=\mathbb E\int_0^T\langle \pi_k\tilde{A}^\top(t+\delta) \mathbb E^{\mathcal F_t}[y_1^j(t+\delta)]I_{[0,T-\delta]},x_l^{**}(t)\rangle dt,\\
&\mathbb E\int_0^T\sum_{k=1}^K\sum_{l=1}^K\langle \pi_k\tilde{A}^\top y_2^k,x_l^{**}(t-\delta)\rangle dt=\mathbb E\int_0^T\sum_{k=1}^K\sum_{l=1}^K\langle \pi_l\tilde{A}^\top(t+\delta) \mathbb E^{\mathcal F_t}[y_2^l(t+\delta)]I_{[0,T-\delta]},x_k^{**}(t)\rangle dt,\\
&\mathbb E\int_0^T\langle \pi_k\tilde{A}^\top y_2^k,\delta x_i(t-\delta)\rangle dt=\mathbb E\int_0^T\langle \pi_k\tilde{A}^\top(t+\delta) \mathbb E^{\mathcal F_t}[y_2^k(t+\delta)]I_{[0,T-\delta]},\delta x_i(t)\rangle dt,\\
&\mathbb E\int_0^T\langle \pi_k\tilde{B}^\top y_2^k,\delta u_i(t-\theta)\rangle dt=\mathbb E\int_0^T\langle \pi_k\tilde{B}^\top(t+\theta) \mathbb E^{\mathcal F_t}[y_2^k(t+\theta)]I_{[0,T-\theta]},\delta u_i(t)\rangle dt,\\
&\mathbb E\int_0^T\langle \pi_k\hat{A}_k^\top y_2^k,x_k^{**}(t-\delta)\rangle dt=\mathbb E\int_0^T\langle \pi_k\hat{A}_k^\top(t+\delta) \mathbb E^{\mathcal F_t}[y_2^k(t+\delta)]I_{[0,T-\delta]},x_k^{**}(t)\rangle dt.
\end{aligned}\end{equation*}
Letting
\begin{equation*}\left\{\begin{aligned}
\alpha_1=&-A_k^\top y_1^j-\hat A_k^\top(t+\delta) \mathbb E^{\mathcal F_t}[y_1^j(t+\delta)]I_{[0,T-\delta]}-(Q+\tilde{Q}(t+\delta))\bar x_j,\\
\alpha_2=&-A_k^\top y_2^k-\hat A_k^\top(t+\delta) \mathbb E^{\mathcal F_t}[y_2^k(t+\delta)]I_{[0,T-\delta]}+\big(QS+ S^\top Q-S^\top QS+\tilde{Q}(t+\delta)\tilde{S}(t+\delta)\\
&+ \tilde{S}^\top(t+\delta) \tilde{Q}(t+\delta)-\tilde{S}^\top(t+\delta) \tilde{Q}(t+\delta)\tilde{S}(t+\delta)\big)\hat x-\sum_{l=1}^K\pi_l\tilde{A}^\top(t+\delta) \mathbb Ey_1^l(t+\delta)I_{[0,T-\delta]}\\
&-\sum_{l=1}^K\pi_l\tilde{A}^\top (t+\delta) y_2^l(t+\delta)I_{[0,T-\delta]},
\end{aligned}\right.\end{equation*}
and substituting \eqref{eq13}-\eqref{eq15} into \eqref{eq12}, we have
\begin{equation}\label{eq16}\begin{aligned}
\Delta \mathcal J_{soc}^{(N)}
  =&\mathbb E\Bigg\{\int_0^T\Big[\langle Q\bar x_i,\delta x_i\rangle-\langle (QS+ S^\top Q -S^\top QS)\hat x,\delta x_i\rangle+\langle \tilde{Q}\bar x_i(t-\delta),\delta x_i(t-\delta)\rangle\\
&-\langle (\tilde{Q}\tilde{S}+ \tilde{S}^\top \tilde{Q} -\tilde{S}^\top \tilde{Q}\tilde{S})\hat x(t-\delta),\delta x_i(t-\delta)\rangle+\sum_{k=1}^K\langle \pi_k\tilde{A}^\top(t+\delta) \mathbb E\mathbf{y}_k(t+\delta)I_{[0,T-\delta]},\delta x_i\rangle\\
&+\sum_{k=1}^K\langle \pi_k\tilde{A}^\top(t+\delta) y_2^k(t+\delta)I_{[0,T-\delta]},\delta x_i\rangle+\sum_{k=1}^K\langle \pi_k\tilde{B}^\top(t+\theta) \mathbb E\mathbf{y}_k(t+\theta)I_{[0,T-\theta]},\delta u_i\rangle\\
&+\sum_{k=1}^K\langle \pi_k\tilde{B}^\top(t+\theta) y_2^k(t+\theta)I_{[0,T-\theta]},\delta u_i\rangle
+\langle R_{\theta_i}\bar u_i,\delta u_i\rangle+\langle \tilde{R}_{\theta_i}\bar u_i(t-\theta),\delta u_i(t-\theta)\rangle\Big]dt \\
&+\langle G\bar x_i(T),\delta x_i(T)\rangle-\langle (G\Gamma+\Gamma ^\top G-\Gamma ^\top G\Gamma)\hat x(T),\delta x_i(T)\rangle\Bigg\}+\sum_{l=1}^{12}\varepsilon_l,
\end{aligned}\end{equation}
where
\begin{equation}\label{eq17}\left\{\begin{aligned}
&dy_1^j=-\Big[A_k^\top y_1^j+\hat A_k^\top(t+\delta) \mathbb E^{\mathcal F_t}[y_1^j(t+\delta)]I_{[0,T-\delta]}+(Q+\tilde{Q}(t+\delta))\bar x_j\Big] dt\\
&\qquad+z_1^jdW_j(t)+\sum_{l=1,l\neq j}^N z_{1}^{jl}dW_l(t),\\
&dy_2^k=-\Big[A_k^\top y_2^k+\hat A_k^\top(t+\delta) y_2^k(t+\delta)I_{[0,T-\delta]}-\big(QS+ S^\top Q-S^\top QS+\tilde{Q}(t+\delta)\tilde{S}(t+\delta)\\
&\qquad +\tilde{S}^\top(t+\delta) \tilde{Q}(t+\delta)-\tilde{S}^\top(t+\delta) \tilde{Q}(t+\delta)\tilde{S}(t+\delta)\big)\hat x+\sum_{l=1}^K\pi_l\tilde{A}^\top(t+\delta) \mathbb E\mathbf{y}_l(t+\delta)I_{[0,T-\delta]}\\
&\qquad+\sum_{l=1}^K\pi_l\tilde{A}^\top (t+\delta) y_2^l(t+\delta)I_{[0,T-\delta]}\Big]dt,\ t\in[0,T],\\
&y_1^j(T)=G\bar x_j(T),\ y_2^k(T)=-(G\Gamma+\Gamma ^\top G-\Gamma ^\top G\Gamma)\hat x(T),
%y_1^j(t)=0,\ z_1^j(t)=0,\
%y_2^k(t)=0, \ t\in(T,T+\delta],
\ j=1,\cdots,N,\ k=1,\cdots,K,
\end{aligned}\right.\end{equation}
and
\small\begin{equation*}\left\{\begin{aligned}
&\varepsilon_{10}=\sum_{k=1}^K\mathbb E\int_0^T\langle \frac{1}{N_k}\sum_{j\in\mathcal I_k,j\neq i}\pi_k\tilde{A}^\top(t+\delta) \mathbb E^{\mathcal F_t}[y_1^j(t+\delta)]I_{[0,T-\delta]}-\pi_k\tilde{A}^\top(t+\delta) \mathbb E\mathbf{y}_k(t+\delta)I_{[0,T-\delta]},\delta x_i\rangle dt,\\
&\varepsilon_{11}=\sum_{k=1}^K\mathbb E\int_0^T\Big\langle \sum_{l=1}^K\frac{\pi_l}{N_l}\sum_{j\in\mathcal I_l,j\neq i}\tilde{A}^\top(t+\delta) \mathbb E^{\mathcal F_t}[y_1^j(t+\delta)]I_{[0,T-\delta]}-\sum_{l=1}^K\pi_l\tilde{A}^\top(t+\delta)\mathbb E\mathbf{y}_l(t+\delta)I_{[0,T-\delta]},x_k^{**}\Big\rangle dt,\\
&\varepsilon_{12}=\sum_{k=1}^K\mathbb E\int_0^T\langle \frac{1}{N_k}\sum_{j\in\mathcal I_k,j\neq i}\pi_k\tilde{B}^\top(t+\theta) \mathbb E^{\mathcal F_t}[y_1^j(t+\theta)]I_{[0,T-\theta]}-\pi_k\tilde{B}^\top(t+\theta) \mathbb E\mathbf{y}_k(t+\theta)I_{[0,T-\theta]},\delta u_i\rangle dt.
\end{aligned}\right.\end{equation*}\normalsize
Note that the states $y_1^j,\ j\in\mathcal I_k,j\neq i$ are exchangeable. When we consider the expectations, we will use $\mathbf y_k$ denote the process $y_1^j$ defined in \eqref{eq17} of the representative agent of type $k$. Moreover, in Section 5 we still use this kind of notations, i.e., use $\mathbf x_k,\mathbf y_k,\mathbf p_k$ to denote the involved processes of the representative agent of type $k$. $\varepsilon_{1}-\varepsilon_{12}$ are actually $o(1)$ order and the rigorous proof will be shown in Section 5.
Therefore, we introduce the decentralized auxiliary cost functional $J_i$ with perturbation as follows:
\begin{equation}\label{eq18}\begin{aligned}
\Delta \mathcal J_i
  =&\mathbb E\Bigg\{\int_0^T\Big[\langle Q\bar x_i,\delta x_i\rangle-\langle (QS+ S^\top Q -S^\top QS)\hat x,\delta x_i\rangle+\langle \tilde{Q}\bar x_i(t-\delta),\delta x_i(t-\delta)\rangle\\
&\quad-\langle (\tilde{Q}\tilde{S}+ \tilde{S}^\top \tilde{Q} -\tilde{S}^\top \tilde{Q}\tilde{S})\hat x(t-\delta),\delta x_i(t-\delta)\rangle+\sum_{k=1}^K\langle \pi_k\tilde{A}^\top(t+\delta) \hat{y}_k(t+\delta)I_{[0,T-\delta]},\delta x_i\rangle\\
&\quad+\sum_{k=1}^K\langle \pi_k\tilde{A}^\top(t+\delta) y_2^k(t+\delta)I_{[0,T-\delta]},\delta x_i\rangle+\sum_{k=1}^K\langle \pi_k\tilde{B}^\top(t+\theta) \hat{y}_k(t+\theta)I_{[0,T-\theta]},\delta u_i\rangle\\
&\quad+\sum_{k=1}^K\langle \pi_k\tilde{B}^\top(t+\theta) y_2^k(t+\theta)I_{[0,T-\theta]},\delta u_i\rangle
+\langle R_{\theta_i}\bar u_i,\delta u_i\rangle+\langle \tilde{R}_{\theta_i}\bar u_i(t-\theta),\delta u_i(t-\theta)\rangle\Big]dt \\
&\quad+\langle G\bar x_i(T),\delta x_i(T)\rangle-\langle (G\Gamma+\Gamma ^\top G-\Gamma ^\top G\Gamma)\hat x(T),\delta x_i(T)\rangle\Bigg\}.
\end{aligned}\end{equation}
\begin{remark}
It is remarkable that due to $\bar x_j$ is an $\mathcal F_t$-adapted stochastic process, thus $ z_1^j, z_{1}^{jl},1\leq l\leq N,l\neq j$ cannot vanish. By contrast, the drift term of $y_2^k$ is deterministic (By Theorem \ref{the1} \emph{(}see below\emph{)}, we know $\hat x$ is deterministic) and the terminal value $y_2^k(T)$ depends on $\hat x(T)$, thus we derive that $z_2^k$ should be zero, which implies $y_2^k$ is deterministic indeed. Thus $\mathbb E[y_2^k(t+\tau)]=y_2^k(t+\tau)$ and $\mathbb E^{\mathcal F_t}[y_2^k(t+\tau)]=y_2^k(t+\tau)$, $t\in[0,T-\tau]$, for $\tau=\delta\vee\theta$. Therefore, system \eqref{eq11} is a system of two SDEs, while the adjoint system is composed by a BSDE and an ODE.
\end{remark}
\begin{remark}
In above analysis, we introduce $N+K$ adjoint processes to break away $\delta \mathcal J_{soc}^{(N)}$ from the dependence on $x_j^*$ and $x_k^{**}$ (By the existence and uniqueness of solutions of \eqref{eq11}, $x_j^*=x_k^{**}$ indeed). This feature is brought by the existence of $x^{(N)}$ in the drift item of state equations, that is $\tilde{A}(\cdot)$. By contrast, if $\tilde{A}(\cdot)\equiv 0$, then $x_j^*(\cdot)\equiv0$ and $x_k^{**}(\cdot)\equiv0$. There's no additional adjoint process is needed to derive the auxiliary problem.
\end{remark}

\subsection{Decentralized strategy}\label{decentralized strategy}

Motivated by \eqref{eq18}, introduce the following auxiliary LQ control problem with delay:\\

\textbf{Problem 2.} Minimize $J_i(u_i)$ over $u_i\in\mathcal U_i^d$ subject to
\begin{equation}\left\{\label{eq19}\begin{aligned}
dx_i(t)=\ &\Big[A_{\theta_i}(t)x_i(t)+\hat A_{\theta_i}(t)x_i(t-\delta)+\tilde A(t)\hat{x}(t-\delta)+B(t)u_i(t)+\hat B(t)u_i(t-\theta)\\
&+\tilde B(t)\hat{u}(t-\theta)\Big]dt+\Big[D(t)u_i(t)+\hat D(t)u_i(t-\theta)\Big]dW_i(t),\ t\in[0,T],\\
x_i(0)=\ &\xi_i,\ x_i(t)=x^0(t),\ t\in[-\delta,0),\ u_i(t)=u^0(t),\ t\in[-\theta,0),
\end{aligned}\right.\end{equation}
where
\begin{equation}\label{eq21}\begin{aligned}
J_i(u_i)
 =&\frac{1}{2}\Bigg\{\mathbb E\int_0^T\Big[\langle Q x_i, x_i\rangle+\langle \tilde{Q} x_i(t-\delta), x_i(t-\delta)\rangle-2\langle \Theta_1, x_i\rangle-2\langle \Theta_2, x_i(t-\delta)\rangle+\langle R_{\theta_i} u_i, u_i\rangle\\
 &\ +\langle \tilde{R}_{\theta_i} u_i(t-\theta), u_i(t-\theta)\rangle+2\langle \Theta_3, u_i\rangle\Big] dt+\langle G x_i(T),x_i(T)\rangle-2\langle\Theta_4,x_i(T)\rangle\Bigg\},
\end{aligned}\end{equation}
with
\begin{equation}\nonumber\begin{aligned}
&\Theta_1=(QS+ S^\top Q -S^\top QS)\hat x-\sum_{k=1}^K \pi_k\tilde{A}^\top(t+\delta) \hat{y}_k(t+\delta)I_{[0,T-\delta]}-\sum_{k=1}^K \pi_k\tilde{A}^\top(t+\delta) y_2^k(t+\delta)I_{[0,T-\delta]},\\
&\Theta_2=(\tilde{Q}\tilde{S}+ \tilde{S}^\top \tilde{Q} -\tilde{S}^\top \tilde{Q}\tilde{S})\hat x(t-\delta),\\
&\Theta_3=\sum_{k=1}^K\pi_k\tilde{B}^\top(t+\theta) \hat{y}_k(t+\theta)I_{[0,T-\theta]}+\sum_{k=1}^K \pi_k\tilde{B}^\top(t+\theta) y_2^k(t+\theta)I_{[0,T-\theta]},\\
& \Theta_4=(G\Gamma+\Gamma ^\top G-\Gamma ^\top G\Gamma)\hat x(T),\\
\end{aligned}\end{equation}
and $\hat x,\ \hat y_k,\ y_2^k$ will be determined by the consistency condition in the following section.

Similar to \cite{CW10}, \cite{Yu12}, etc, we will apply stochastic maximum principle to study \textbf{Problem 2}. First introduce the following adjoint equation
\begin{equation*}\left\{\begin{aligned}
&dp_i(t)=-\Big[A_{\theta_i}^\top p_i+\hat A_{\theta_i}^\top(t+\delta)\mathbb E^{\mathcal F_t}[p_i(t+\delta)]I_{[0,T-\delta]}+(Q+\tilde{Q}(t+\delta))\bar x_i-\Theta_1-\Theta_2(t+\delta)\Big]dt\\
&\qquad\qquad+q_idW_i(t),\\
&p_i(T)=G \bar x_i(T)-\Theta_4%,\ p_i(t)=0, \ q_i(t)=0,\ q_{ij}(t)=0,\ t\in(T,T+\delta\vee\theta]
.
\end{aligned}\right.\end{equation*}
The global stochastic maximum principle implies that
\begin{equation}\label{eq optimal control}\begin{aligned}
\bar u_i(t)=&-\big(R_{\theta_i}(t)+\tilde{R}_{\theta_i}(t+\theta)\big)^{-1}\big(B^\top(t) p_i(t)+\hat B^\top(t+\theta)\mathbb E^{\mathcal F_t}[p_i(t+\theta)]I_{[0,T-\theta]}+D^\top(t) q_i(t)\\
&\quad+\hat{D}^\top(t+\theta)\mathbb E^{\mathcal F_t}[q_i(t+\theta)]I_{[0,T-\theta]}+\Theta_3\big),\ t\in[0,T],\ \mathbb P-a.s..
\end{aligned}\end{equation}
 The related Hamiltonian system becomes
\small\begin{equation}\label{Hamiltonian system}\left\{\begin{aligned}
&dx_i(t)=\Big[A_{\theta_i}x_i+\hat A_{\theta_i}x_i(t-\delta)+\tilde A\hat{x}(t-\delta)-B\big(R_{\theta_i}+\tilde{R}_{\theta_i}(t+\theta)\big)^{-1}\big(B^\top p_i+\hat B^\top (t+\theta)\mathbb E^{\mathcal F_t}[p_i(t+\theta)]I_{[0,T-\theta]}\\
&\qquad+D^\top q_i+\hat{D}^\top(t+\theta)\mathbb E^{\mathcal F_t}[q_i(t+\theta)]I_{[0,T-\theta]}+\Theta_3\big)-\hat B\big(R_{\theta_i}(t-\theta)+\tilde{R}_{\theta_i}\big)^{-1}\big(B^\top(t-\theta) p_i(t-\theta)\\
&\qquad+\hat B^\top \mathbb E^{\mathcal F_{t-\theta}}[p_i]+D^\top(t-\theta) q_i(t-\theta)+\hat{D}^\top\mathbb E^{\mathcal F_{t-\theta}}[q_i]+\Theta_3(t-\theta)\big)I_{[\theta,T]}+\hat Bu^0I_{[0,\theta]}+\tilde B\hat{u}(t-\theta)I_{[\theta,T]}\\
&\qquad+\tilde{B}u^0I_{[0,\theta]}\Big]dt-\Big[D\big(R_{\theta_i}+\tilde{R}_{\theta_i}(t+\theta)\big)^{-1}\big(B^\top p_i+\hat B^\top(t+\theta)\mathbb E^{\mathcal F_t}[p_i(t+\theta)]I_{[0,T-\theta]}+D^\top q_i\\
&\qquad+\hat{D}^\top(t+\theta)\mathbb E^{\mathcal F_t}[q_i(t+\theta)]I_{[0,T-\theta]}+\Theta_3\big)+\hat D\big(R_{\theta_i}(t-\theta)+\tilde{R}_{\theta_i}(t)\big)^{-1}\big(B^\top(t-\theta) p_i(t-\theta)\\
&\qquad+\hat B^\top \mathbb E^{\mathcal F_{t-\theta}}[p_i]+D^\top(t-\theta) q_i(t-\theta)+\hat{D}^\top\mathbb E^{\mathcal F_{t-\theta}}[q_i]+\Theta_3(t-\theta)\big)I_{[\theta,T]}-\hat Du^0I_{[0,\theta]}\Big]dW_i(t),\ t\in[0,T],\\
&dp_i(t)=-\Big[A_{\theta_i}^\top p_i+\hat A_{\theta_i}^\top(t+\delta)\mathbb E^{\mathcal F_t}[p_i(t+\delta)]I_{[0,T-\delta]}+(Q+\tilde{Q}(t+\delta)) x_i-\Theta_1-\Theta_2(t+\delta)\Big]dt\\
&\qquad+q_idW_i(t),\ t\in[0,T],\\
&x_i(0)=\xi_i,\ x_i(t)=x^0(t),\ t\in[-\delta,0),\ u_i(t)=u^0(t),\ t\in[-\theta,0),\ p_i(T)=G \bar x_i(T)-\Theta_4,%\ p_i(t)=0,\ q_i(t)=0,\ q_{ij}(t)=0,\ t\in(T,T+\delta\vee\theta],
\end{aligned}\right.\end{equation}
\normalsize
where $\Theta_2(t+\delta)=(\tilde{Q}(t+\delta)\tilde{S}(t+\delta)+ \tilde{S}^\top(t+\delta) \tilde{Q}(t+\delta) -\tilde{S}^\top(t+\delta) \tilde{Q}(t+\delta)\tilde{S}(t+\delta))\hat x(t)$, $\Theta_3(t-\theta)=\sum\limits_{k=1}^K\pi_k\tilde{B}^\top \hat{y}_kI_{[\theta,T]}+\sum\limits_{k=1}^K \pi_k\tilde{B}^\top y_2^kI_{[\theta,T]}$.

\section{Consistency condition}\label{CC}

%\subsection{Mean-field FBSDE}
In this section, we focus on the solution of \textbf{Problem 2} by constructing the CC system and solving it. For the sake of presentation, in the following we let $n=m$. There is no essential difference if $n\neq m$.
First let us define some notations:
\begin{equation}\nonumber\begin{aligned}
&\mathbb{R}_k(\cdot)=R_k(\cdot)+\tilde{R}_k(\cdot+\theta),\ \mathbb{Q}(\cdot)=Q(\cdot)+\tilde{Q}(\cdot+\delta),\ \mathbb{S}(\cdot)=Q(\cdot)S(\cdot)+ S^\top(\cdot) Q(\cdot)\\
&\qquad\quad-S^\top(\cdot) Q(\cdot)S(\cdot)+\tilde{Q}(\cdot+\delta)\tilde{S}(\cdot+\delta)+ \tilde{S}^\top(\cdot+\delta) \tilde{Q}(\cdot+\delta)-\tilde{S}^\top(\cdot+\delta) \tilde{Q}(\cdot+\delta)\tilde{S}(\cdot+\delta),\\
&\mathbb{B}_1^k(\cdot)=B(\cdot)\mathbb{R}_k^{-1}(\cdot)B^\top(\cdot),\ \mathbb{B}_2^k(\cdot)=\tilde{B}(\cdot)\mathbb{R}_k^{-1}(\cdot-\theta)B^\top(\cdot-\theta),\ \mathbb{B}_3^k(\cdot)=D(\cdot)\mathbb{R}_k^{-1}(\cdot)B^\top(\cdot),\\
&\mathbb{B}_4^k(\cdot)=\hat{D}(\cdot)\mathbb{R}_k^{-1}(\cdot-\theta)B^\top(\cdot-\theta),\ \mathbb{B}_5^k(\cdot)=B(\cdot)\mathbb{R}_k^{-1}(\cdot)\hat B^\top(\cdot+\theta),\ \mathbb{B}_6^k(\cdot)=\hat B(\cdot)\mathbb{R}_k^{-1}(\cdot-\theta)B^\top(\cdot-\theta),\\
&\mathbb{B}_7^k(\cdot)=\hat B(\cdot)\mathbb{R}_k^{-1}(\cdot-\theta)\hat{B}^\top(\cdot),\ \mathbb{B}_8^k(\cdot)=\tilde B(\cdot)\mathbb{R}_k^{-1}(\cdot-\theta)\hat{B}^\top(\cdot),\ \mathbb{B}_9^k(\cdot)=D(\cdot)\mathbb{R}_k^{-1}(\cdot)\hat B^\top(\cdot+\theta),\\
&\mathbb{B}_{10}^k(\cdot)=\hat D(\cdot)\mathbb{R}_k^{-1}(\cdot-\theta)\hat B^\top(\cdot),\ \mathbb{B}_{11}^k(\cdot)= B(\cdot)\mathbb{R}_k^{-1}(\cdot)\tilde{B}^\top(\cdot+\theta),\ \mathbb{B}_{12}^k(\cdot)=\hat B(\cdot)\mathbb{R}_k^{-1}(\cdot-\theta)\tilde{B}^\top(\cdot),\\
&\mathbb{B}_{13}^k(\cdot)=D(\cdot)\mathbb{R}_k^{-1}(\cdot)\tilde B^\top(\cdot+\theta),\ \mathbb{B}_{14}^k(\cdot)=\hat D(\cdot)\mathbb{R}_k^{-1}(\cdot-\theta)\tilde B^\top(\cdot),\mathbb{D}_1^k(\cdot)=B(\cdot)\mathbb{R}_k^{-1}(\cdot)D^\top(\cdot),\\
&\mathbb{D}_2^k(\cdot)=B(\cdot)\mathbb{R}_k^{-1}(\cdot)\hat{D}^\top(\cdot+\theta),\ \mathbb{D}_3^k(\cdot)=\tilde{B}(\cdot)\mathbb{R}_k^{-1}(\cdot-\theta)D^\top(\cdot-\theta),\ \mathbb{D}_4^k(\cdot)=\tilde{B}(\cdot)\mathbb{R}_k^{-1}(\cdot-\theta)\hat{D}^\top(\cdot),\\
\end{aligned}\end{equation}
\begin{equation}\label{notation}\begin{aligned}
&\mathbb{D}_5^k(\cdot)=D(\cdot)\mathbb{R}_k^{-1}(\cdot)D^\top(\cdot),\ \mathbb{D}_6^k(\cdot)=D(\cdot)\mathbb{R}_k^{-1}(\cdot)\hat{D}^\top(\cdot+\theta),\ \mathbb{D}_7^k(\cdot)=\hat{D}(\cdot)\mathbb{R}_k^{-1}(\cdot-\theta)D^\top(\cdot-\theta),\\
&\mathbb{D}_8^k(\cdot)=\hat{D}(\cdot)\mathbb{R}_k^{-1}(\cdot-\theta)\hat{D}^\top(\cdot),\ \mathbb{D}_9^k(\cdot)=\hat B(\cdot)\mathbb{R}_k^{-1}(\cdot-\theta)D^\top(\cdot-\theta),\ \mathbb{D}_{10}^k(\cdot)=\hat B(\cdot)\mathbb{R}_k^{-1}(\cdot-\theta)\hat D^\top(\cdot),\\ &\mathbb{G}=G\Gamma+\Gamma^\top G-\Gamma^\top G\Gamma.
\end{aligned}\end{equation}

From \eqref{eq17}, \eqref{eq optimal control}, \eqref{Hamiltonian system}, recalling that we use $\hat{x}$ to approximate $\bar x^{N}$ and use $\hat{y}_k$ to approximate $\mathbb{E}y_k$, where $y_k=y_1^j$ for $j\in \mathcal{I}_k$, we have the following theorem

\begin{theorem} \label{the1}
Let \emph{(}A1\emph{)}-\emph{(}A4\emph{)} hold. The parameters in \textbf{Problem 2} can be determined by
$$(\hat x,\ \hat u,\ \hat y_k,\ y_2^k)=\Big(\sum_{l=1}^K\pi_l\mathbb E\alpha_l,\ \sum_{l=1}^K\pi_l\mathbb Ev_l,\ \mathbb E\check y_k,\ \zeta_k\Big),$$ where
\begin{equation*}\begin{aligned}
v_k(t)=-&\big(R_k(t)+\tilde{R}_k(t+\theta)\big)^{-1}\big(B^\top(t) \beta_k(t)+\hat B^\top(t+\theta)\mathbb E^{\mathcal F_t}[\beta_k(t+\theta)]I_{[0,T-\theta]}\\
&+D^\top(t) \gamma_k(t)+\hat{D}^\top(t+\theta)\mathbb E^{\mathcal F_t}[\gamma_k(t+\theta)]I_{[0,T-\theta]}+\sum_{k=1}^K\pi_k\tilde{B}^\top(t+\theta) \mathbb E\check y_k(t+\theta)I_{[0,T-\theta]}\\
&+\sum_{k=1}^K \pi_k\tilde{B}^\top(t+\theta) \zeta_k(t+\theta)I_{[0,T-\theta]}\big),
\end{aligned}\end{equation*}
 and
$(\alpha_k,\ \beta_k,\ \gamma_k,\ \check y_k,\ \check z_k,\ \zeta_k)$ is the solution of the following MF-AFBSDDEs, which is so-called \emph{(}CC\emph{)} system: for $k=1,\cdots,K$,
\begin{equation}\label{CC}\left\{\begin{aligned}
&d\alpha_k(t)=\Big[A_k\alpha_k+\hat A_k\alpha_k(t-\delta)+\tilde{A}\sum_{l=1}^K\pi_l\mathbb E\alpha_l(t-\delta)-\big(\mathbb{B}_1^k \beta_k+\mathbb{B}_5^k \mathbb E^{\mathcal F_t}[\beta_k(t+\theta)]I_{[0,T-\theta]}+\mathbb{D}_1^k \gamma_k\\
&\qquad\qquad+\mathbb{D}_2^k \mathbb E^{\mathcal F_t}[\gamma_k(t+\theta)]I_{[0,T-\theta]}+\mathbb{B}_{11}^k\sum_{l=1}^K\pi_l(\mathbb E\check{y}_l(t+\theta)+\zeta_l(t+\theta))I_{[0,T-\theta]}\big)-\big(\mathbb{B}_6^k \beta_k(t-\theta)\\
&\qquad\qquad+\mathbb{B}_7^k \mathbb E^{\mathcal F_{t-\theta}}[\beta_k]+\mathbb{D}_9^k \gamma_k(t-\theta)+\mathbb{D}_{10}^k \mathbb E^{\mathcal F_{t-\theta}}[\gamma_k]+\mathbb{B}_{12}^k\sum_{l=1}^K\pi_l(\mathbb E\check{y}_l+\zeta_l)\big)I_{[\theta,T]}\\
&\qquad\qquad+\hat Bu^0I_{[0,\theta]}-\sum_{l=1}^K\pi_l \big(\mathbb{B}_2^l \mathbb{E}\beta_l(t-\theta)+\mathbb{B}_8^l \mathbb{E}\beta_l+\mathbb{D}_3^l \mathbb{E}\gamma_l(t-\theta)+\mathbb{D}_4^l \mathbb E\gamma_l\big)I_{[\theta,T]}+\tilde Bu^0I_{[0,\theta]}\Big]dt\\
&\qquad\qquad-\Big[\mathbb{B}_3^k \beta_k+\mathbb{B}_9^k  \mathbb E^{\mathcal F_t}[\beta_k(t+\theta)]I_{[0,T-\theta]}+\mathbb{D}_5^k \gamma_k+\mathbb{D}_6^k \mathbb E^{\mathcal F_t}[\gamma_k(t+\theta)]I_{[0,T-\theta]}\\
&\qquad\qquad+\mathbb{B}_{13}^k\sum_{l=1}^K\pi_l(\mathbb E\check{y}_l(t+\theta)+\zeta_l(t+\theta))I_{[0,T-\theta]}+\big(\mathbb{B}_4^k \beta_k(t-\theta)+\mathbb{B}_{10}^k \mathbb E^{\mathcal F_{t-\theta}}[\beta_k]+\mathbb{D}_7^k \gamma_k(t-\theta)\\
&\qquad\qquad+\mathbb{D}_8^k \mathbb E^{\mathcal F_{t-\theta}}[\gamma_k]+\mathbb{B}_{14}^k\sum_{l=1}^K\pi_l(\mathbb E\check{y}_l+\zeta_l)\big)I_{[\theta,T]}-\hat Du^0I_{[0,\theta]}\Big]dW_k(t),\\
&d\beta_k(t)=-\Big[A_k^\top\beta_k+\hat A_k^\top(t+\delta)\mathbb E^{\mathcal F_t}[\beta_k(t+\delta)]I_{[0,T-\delta]}+\mathbb{Q}\alpha_k-\mathbb{S}\sum_{l=1}^K\pi_l\mathbb E\alpha_l\\
&\qquad\qquad+\sum_{l=1}^K\pi_l\tilde{A}^\top(t+\delta)\mathbb E\check{y}_l(t+\delta)I_{[0,T-\delta]}+\sum_{l=1}^K\pi_l\tilde{A}^\top(t+\delta)\zeta_l(t+\delta)I_{[0,T-\delta]}\Big]dt+\gamma_kdW_k(t),\\
&d\check y_k(t)=-\Big[A_k^\top \check y_k+\hat{A}^\top_k(t+\delta)\mathbb E^{\mathcal F_t}[\check{y}_k(t+\delta)]I_{[0,T-\delta]}+\mathbb{Q}\alpha_k\Big] dt+\check z_kdW_k(t),\\
&d\zeta_k(t)=-\Big[A_k^\top \zeta_k+\hat A_k^\top(t+\delta) \zeta_k(t+\delta)I_{[0,T-\delta]}-\mathbb{S}\sum_{l=1}^K\pi_l\mathbb E\alpha_l+\sum_{l=1}^K\pi_l\tilde{A}^\top(t+\delta) \mathbb E\check y_l(t+\delta)I_{[0,T-\delta]}\\
&\qquad\qquad+\sum_{l=1}^K\pi_l\tilde{A}^\top(t+\delta) \zeta_l(t+\delta)I_{[0,T-\delta]}\Big]dt,\\
&\alpha_k(0)=\xi_k,\  \beta_k(T)=G \alpha_k(T)-\mathbb{G}\sum_{l=1}^K\pi_l\mathbb E\alpha_l(T),\ \check y_k(T)=G \alpha_k(T),\  \zeta_k(T)=-\mathbb{G}\sum_{l=1}^K\pi_l\mathbb E\alpha_l(T),\\
&\alpha_k(t)=x^0(t),\ t\in[-\delta,0)
%,\ \beta_k(t)=0,\ \gamma_k(t)=0,\ t\in(T,T+\delta\vee\theta],\\
%&\check y_k(t)=0,\ \check z_k(t)=0,\ \zeta_k(t)=0,\ t\in(T,T+\delta]
.
\end{aligned}\right.\end{equation}
\end{theorem}

\begin{remark}
(i) It is remarkable that if \eqref{CC} is solved, by the estimates of AFBSDDE, we can easily obtain $\hat x=\sum_{l=1}^K\pi_l\mathbb E\alpha_l$, so does $\hat u=\sum_{l=1}^K\pi_l\mathbb Ev_l$, $\hat y_k=\mathbb E\check y_k$, $y_2^k=\zeta_k$. Actually, the CC system \eqref{CC} is a coupled MF-AFBSDDE composed by not only one SDDE,  two ABSDEs and an ODE, but also the MF terms and conditional expectations (caused by the delay feature). If taking the expectation to \eqref{CC}, we can also obtain $\hat x=\sum_{l=1}^K\pi_l\mathbb E\alpha_l$. However, in the consideration of the generalization, we focus the coupled MF-AFBSDDE here.

(ii) In \eqref{CC}, for $k=1,\cdots,K$, the subscript $k$ \emph{(}e.g.,$\ \alpha_k,\beta_k,\cdots$\emph{)} stands for a representative agent in the $k$-type.
\end{remark}
In the following, we give a proposition to obtain the well-posedness of CC system \eqref{CC}. Before solving \eqref{CC}, let us make some transformations and introduce some notations. Define $X=(\alpha_1^\top,\cdots,\alpha_K^\top)^\top$, $Y=( \beta_1^\top,\cdots,\beta_K^\top,\check y_1^\top,\cdots,\check y_K^\top,\zeta_1^\top,\cdots,\zeta_K^\top)^\top$. Denote $W^{(k)}$ by a representative Brownian motion in the $k$-type, the MF-AFBSDDE system \eqref{CC} then takes the following form
\small\begin{equation}\label{eq23}\left\{\begin{aligned}
&dX=\Big[\mathbf{A}X+\hat{\mathbf{A}}X(t-\delta)+\widetilde{\mathbf{A}}_1^\pi\mathbb EX(t-\delta)-\mathbf B_1Y-\mathbf B_1^\pi \mathbb EY(t-\theta)I_{[\theta,T]}-\mathbf B_2\mathbb E^{\mathcal F_t}[Y(t+\theta)]I_{[0,T-\theta]}\\
&\qquad+\mathbf B_2^\pi \mathbb EY-\mathbf B_3Y(t-\theta)I_{[\theta,T]}-\mathbf B_3^\pi\mathbb EY(t+\theta)I_{[0,T-\theta]}-\mathbf B_4\mathbb E^{\mathcal F_{t-\theta}}[Y]I_{[\theta,T]}-\mathbf D_1(Z)\\
&\qquad-\mathbf D_1^\pi\big(\mathbb EZ(t-\theta)\big)I_{[\theta,T]}-\mathbf D_2\big(\mathbb E^{\mathcal F_t}[Z(t+\theta)]\big)I_{[0,T-\theta]}-\mathbf D_2^\pi(\mathbb EZ)-\mathbf D_3(Z(t-\theta))I_{[\theta,T]}\\
&\qquad-\mathbf D_4(\mathbb E^{\mathcal F_{t-\theta}}[Z])I_{[\theta,T]}+\big(\hat B+\tilde B\big)u^0I_{[0,\theta]}\Big]dt-\Big[\mathbf B_5(Y)+\mathbf B_6 \big(Y(t-\theta)\big)I_{[\theta,T]}\\
&\qquad+\mathbf B_7\big(\mathbb E^{\mathcal F_t}[Y(t+\theta)]\big)I_{[0,T-\theta]}+\mathbf B_8\big(\mathbb E^{\mathcal F_{t-\theta}}[Y]\big)I_{[\theta,T]}+\mathbf B_4^\pi\mathbb EY(t+\theta)I_{[0,T-\theta]}+\mathbf B_5^\pi\mathbb EYI_{[\theta,T]}\\
&\qquad+\mathbf D_5Z+\mathbf D_6Z(t-\theta)I_{[\theta,T]}+\mathbf D_7\mathbb E^{\mathcal F_t}[Z(t+\theta)]I_{[0,T-\theta]}+\mathbf D_8\mathbb E^{\mathcal F_{t-\theta}}[Z]I_{[\theta,T]}-\hat Du^0I_{[0,\theta]}\Big]d\mathbf{W}(t),\\
&dY=-\Big[\mathcal A Y+\hat{\mathcal A}(t+\delta) \mathbb E^{\mathcal F_t}[Y(t+\delta)]I_{[0,T-\delta]}+\mathbf{Q}X-\mathbf{S}^\pi\mathbb EX+\widetilde{\mathbf{A}}_2^\pi\mathbb EY(t+\delta)I_{[0,T-\delta]}\Big]dt\\
&\qquad+Zd\mathbf{W}(t),\ t\in[0,T],\\
&X(0)=\Xi,\  Y(T)=\mathbf{G} X(T)-\mathbf{G}^\pi \mathbb EX(T),\  X(t)=X^0(t),\ t\in[-\delta,0)%,\ Y(t)=0,\ Z(t)=0,\ t\in(T,T+\delta\vee\theta]
.
\end{aligned}\right.\end{equation}\normalsize
where
\tiny\begin{equation*}\begin{aligned}
&\mathbf{A}=\left(
               \begin{smallmatrix}
                 A_1 &  &  \\
                   & \ddots &    \\
                   &   & A_K \\
                \end{smallmatrix}
             \right),
\hat{\mathbf{A}}=\left(
               \begin{smallmatrix}
                 \hat{A}_1 &  &  \\
                   & \ddots &    \\
                   &   & \hat{A}_K \\
                \end{smallmatrix}
             \right),
\widetilde{\mathbf{A}}_1^\pi=\left(
                         \begin{smallmatrix}
                           \tilde{A}\pi_1 & \cdots & \tilde{A}\pi_K   \\
                           \vdots &   & \vdots   \\
                           \tilde{A}\pi_1 & \cdots & \tilde{A}\pi_K \\
                         \end{smallmatrix}
                       \right),
\mathbf B_1=\left(
               \begin{smallmatrix}
                 \mathbb{B}_1^1 & & &0 & & &0& & \\
                     & \ddots &   &   &\ddots & & &\ddots &  \\
                    &   & \mathbb{B}_1^K &   &   &0& & &0   \\
                 \end{smallmatrix}
             \right),
\mathbf B_2=\left(
               \begin{smallmatrix}
                 \mathbb{B}_5^1 & & &0 & & &0& & \\
                     & \ddots &   &   &\ddots & & &\ddots &  \\
                    &   & \mathbb{B}_5^K &   &   &0& & &0   \\
                 \end{smallmatrix}
             \right),\\
&\mathbf B_1^\pi=\left(
                         \begin{smallmatrix}
                           \mathbb{B}_2^1\pi_1 & \cdots & \mathbb{B}_2^K\pi_K & 0 & \cdots & 0  & 0 & \cdots & 0 \\
                           \vdots &   & \vdots &\vdots &   & \vdots  &\vdots &   & \vdots   \\
                           \mathbb{B}_2^1\pi_1 & \cdots & \mathbb{B}_2^K\pi_K & 0 & \cdots & 0  & 0 & \cdots & 0 \\
                         \end{smallmatrix}
                       \right),
\mathbf B_2^\pi=\left(
                         \begin{smallmatrix}
                           -\mathbb{B}_8^1\pi_1 & \cdots & -\mathbb{B}_8^K\pi_K & \mathbb{B}_{12}^1\pi_1 & \cdots & \mathbb{B}_{12}^1\pi_K  & \mathbb{B}_{12}^1\pi_1 & \cdots & \mathbb{B}_{12}^1\pi_K \\
                           \vdots &   & \vdots &\vdots &   & \vdots  &\vdots &   & \vdots   \\
                           -\mathbb{B}_8^1\pi_1 & \cdots & -\mathbb{B}_8^K\pi_K & \mathbb{B}_{12}^K\pi_1 & \cdots & \mathbb{B}_{12}^K\pi_K  & \mathbb{B}_{12}^K\pi_1 & \cdots & \mathbb{B}_{12}^K\pi_K \\
                         \end{smallmatrix}
                       \right),
\mathbf B_3=\left(
               \begin{smallmatrix}
                 \mathbb{B}_6^1 & & &0 & & &0& & \\
                     & \ddots &   &   &\ddots & & &\ddots &  \\
                    &   & \mathbb{B}_6^K &   &   &0& & &0   \\
                 \end{smallmatrix}
             \right),\\
&\mathbf B_3^\pi=\left(
                         \begin{smallmatrix}
                           0 & \cdots & 0 & \mathbb{B}_{11}^1\pi_1 & \cdots & \mathbb{B}_{11}^1\pi_K  & \mathbb{B}_{11}^1\pi_1 & \cdots & \mathbb{B}_{11}^1\pi_K \\
                           \vdots &   & \vdots &\vdots &   & \vdots  &\vdots &   & \vdots   \\
                           0 & \cdots & 0 & \mathbb{B}_{11}^K\pi_1 & \cdots & \mathbb{B}_{11}^K\pi_K  & \mathbb{B}_{11}^K\pi_1 & \cdots & \mathbb{B}_{11}^K\pi_K \\
                         \end{smallmatrix}
                       \right),
\mathbf B_4^\pi=\left(
                         \begin{smallmatrix}
                           0 & \cdots & 0 & \mathbb{B}_{13}^1\pi_1 & \cdots & \mathbb{B}_{13}^1\pi_K  & \mathbb{B}_{13}^1\pi_1 & \cdots & \mathbb{B}_{13}^1\pi_K \\
                           \vdots &   & \vdots &\vdots &   & \vdots  &\vdots &   & \vdots   \\
                           0 & \cdots & 0 & \mathbb{B}_{13}^K\pi_1 & \cdots & \mathbb{B}_{13}^K\pi_K  & \mathbb{B}_{13}^K\pi_1 & \cdots & \mathbb{B}_{13}^K\pi_K \\
                         \end{smallmatrix}
                       \right),
\\
&\mathbf B_4=\left(
               \begin{smallmatrix}
                 \mathbb{B}_7^1 & & &0 & & &0& & \\
                     & \ddots &   &   &\ddots & & &\ddots &  \\
                    &   & \mathbb{B}_7^K &   &   &0& & &0   \\
                 \end{smallmatrix}
             \right),
             \mathbf B_5(Y)=\left(
               \begin{smallmatrix}
                 \mathbb{B}_3^1\beta_1 &  &  \\
                   & \ddots &    \\
                   &   & \mathbb{B}_3^K\beta_K \\
                \end{smallmatrix}
             \right),
\mathbf B_6 \big(Y(t-\theta)\big)=\left(
               \begin{smallmatrix}
                 \mathbb{B}_4^1\beta_1(t-\theta) &  &  \\
                   & \ddots &    \\
                   &   & \mathbb{B}_4^K\beta_K(t-\theta) \\
                \end{smallmatrix}
             \right),
\hat{\mathcal{A}}=\left(
               \begin{smallmatrix}
                 \hat{\mathbf{A}}^\top &  &  \\
                   & \hat{\mathbf{A}}^\top &    \\
                   &   & \hat{\mathbf{A}}^\top \\
                \end{smallmatrix}
             \right),\\
&\mathbf B_7\big(\mathbb E^{\mathcal F_t}[Y(t+\theta)]\big)=\left(
               \begin{smallmatrix}
                 \mathbb{B}_9^1\mathbb E^{\mathcal F_t}[\beta_1(t+\theta)] &  &  \\
                   & \ddots &    \\
                   &   & \mathbb{B}_9^K\mathbb E^{\mathcal F_t}[\beta_K(t+\theta)] \\
                \end{smallmatrix}
             \right),
\mathbf B_8\big(\mathbb E^{\mathcal F_{t-\theta}}[Y]\big)= \left(
               \begin{smallmatrix}
                 \mathbb{B}_{10}^1\mathbb E^{\mathcal F_{t-\theta}}[\beta_1] &  &  \\
                   & \ddots &    \\
                   &   & \mathbb{B}_{10}^K\mathbb E^{\mathcal F_{t-\theta}}[\beta_K] \\
                \end{smallmatrix}
             \right),
\mathcal{A}=\left(
               \begin{smallmatrix}
                 \mathbf{A}^\top &  &  \\
                   & \mathbf{A}^\top &    \\
                   &   & \mathbf{A}^\top \\
                \end{smallmatrix}
             \right),\\
&\mathbf D_1(Z)=\left(
               \begin{smallmatrix}
                 \mathbb{D}_1^1\gamma_1  \\
                    \vdots  \\
                    \mathbb{D}_1^K\gamma_K   \\
                 \end{smallmatrix}
             \right),
\mathbf D_1^\pi\big(\mathbb EZ(t-\theta)\big)=\left(
               \begin{smallmatrix}
                 \sum_{l=1}^K\pi_l\mathbb D_3^l\mathbb E\gamma_l(t-\theta)  \\
                    \vdots  \\
                    \sum_{l=1}^K\pi_l\mathbb D_3^l\mathbb E\gamma_l(t-\theta)    \\
                 \end{smallmatrix}
             \right),
\mathbf D_2\big(\mathbb E^{\mathcal F_t}[Z(t+\theta)]\big)=\left(
               \begin{smallmatrix}
                 \mathbb D_2^1\mathbb E^{\mathcal F_t}[\gamma_1(t+\theta)]  \\
                    \vdots  \\
                 \mathbb D_2^K\mathbb E^{\mathcal F_t}[\gamma_K(t+\theta)]    \\
                 \end{smallmatrix}
             \right),
\mathbf D_2^\pi(\mathbb EZ)=\left(
               \begin{smallmatrix}
                 \sum_{l=1}^K\pi_l\mathbb D_4^l\mathbb E\gamma_l \\
                    \vdots  \\
                 \sum_{l=1}^K\pi_l\mathbb D_4^l\mathbb E\gamma_l    \\
                 \end{smallmatrix}
             \right),\\
&\mathbf D_5=\left(
               \begin{smallmatrix}
                 \mathbb{D}_5^1 & & &0 & & &0& & \\
                     & \ddots &   &   &\ddots & & &\ddots &  \\
                    &   & \mathbb{D}_5^K &   &   &0& & &0   \\
                 \end{smallmatrix}
             \right),
\mathbf D_6=\left(
               \begin{smallmatrix}
                 \mathbb{D}_7^1 & & &0 & & &0& & \\
                     & \ddots &   &   &\ddots & & &\ddots &  \\
                    &   & \mathbb{D}_7^K &   &   &0& & &0   \\
                 \end{smallmatrix}
             \right),
\mathbf D_7=\left(
               \begin{smallmatrix}
                 \mathbb{D}_6^1 & & &0 & & &0& & \\
                     & \ddots &   &   &\ddots & & &\ddots &  \\
                    &   & \mathbb{D}_6^K &   &   &0& & &0   \\
                 \end{smallmatrix}
             \right),
\mathbf D_8=\left(
               \begin{smallmatrix}
                 \mathbb{D}_8^1 & & &0 & & &0& & \\
                     & \ddots &   &   &\ddots & & &\ddots &  \\
                    &   & \mathbb{D}_8^K &   &   &0& & &0   \\
                 \end{smallmatrix}
             \right),\\
&\mathbf D_3\big(Z(t-\theta)\big)=\left(
               \begin{smallmatrix}
                 \mathbb D_9^1\gamma_1(t-\theta)  \\
                    \vdots  \\
                 \mathbb D_9^K\gamma_K(t-\theta)    \\
                 \end{smallmatrix}
             \right),
\mathbf D_4\big(\mathbb E^{\mathcal F_{t-\theta}}[Z]\big)=\left(
               \begin{smallmatrix}
                 \mathbb D_{10}^1\mathbb E^{\mathcal F_{t-\theta}}[\gamma_1]  \\
                    \vdots  \\
                 \mathbb D_{10}^K\mathbb E^{\mathcal F_{t-\theta}}[\gamma_K]    \\
                 \end{smallmatrix}
             \right),
\mathbf D_1=\left(
               \begin{smallmatrix}
                 \mathbb{D}_1^1 & & \\
                     & \ddots &     \\
                    &   & \mathbb{D}_1^K   \\
                 \end{smallmatrix}
             \right),
\mathbf D_1^\pi=\left(
               \begin{smallmatrix}
                 \sum_{l=1}^K\pi_l\mathbb D_3^l \\
                    \vdots  \\
                    \sum_{l=1}^K\pi_l\mathbb D_3^l    \\
                 \end{smallmatrix}
             \right),
\mathbf D_2=\left(
               \begin{smallmatrix}
                 \mathbb{D}_1^2 & & \\
                     & \ddots &     \\
                    &   & \mathbb{D}_2^K   \\
                 \end{smallmatrix}
             \right),\\
&\mathbf D_2^\pi=\left(
               \begin{smallmatrix}
                 \sum_{l=1}^K\pi_l\mathbb D_4^l \\
                    \vdots  \\
                  \sum_{l=1}^K\pi_l\mathbb D_4^l    \\
                 \end{smallmatrix}
             \right),
\mathbf D_3=\left(
               \begin{smallmatrix}
                 \mathbb D_9^1  \\
                    \vdots  \\
                 \mathbb D_9^K   \\
                 \end{smallmatrix}
             \right),
\mathbf D_4=\left(
               \begin{smallmatrix}
                 \mathbb D_{10}^1  \\
                    \vdots  \\
                 \mathbb D_{10}^K  \\
                 \end{smallmatrix}
             \right),
\mathbf B_5=\left(
               \begin{smallmatrix}
                 \mathbb{B}_3^1 &  &  \\
                   & \ddots &    \\
                   &   & \mathbb{B}_3^K \\
                \end{smallmatrix}
             \right),
\mathbf B_6 =\left(
               \begin{smallmatrix}
                 \mathbb{B}_4^1 &  &  \\
                   & \ddots &    \\
                   &   & \mathbb{B}_4^K \\
                \end{smallmatrix}
             \right),
\mathbf B_7=\left(
               \begin{smallmatrix}
                 \mathbb{B}_9^1 &  &  \\
                   & \ddots &    \\
                   &   & \mathbb{B}_9^K \\
                \end{smallmatrix}
             \right),
\mathbf B_8= \left(
               \begin{smallmatrix}
                 \mathbb{B}_{10}^1 &  &  \\
                   & \ddots &    \\
                   &   & \mathbb{B}_{10}^K \\
                \end{smallmatrix}
             \right),\\
\end{aligned}\end{equation*}
\begin{equation*}\begin{aligned}
&Z=\left(
               \begin{smallmatrix}
                 \gamma_1 & & \\
                    &\ddots&  \\
                    & & \gamma_K   \\
                   \check z_1 & & \\
                    &\ddots&  \\
                    & & \check z_K  \\
                     0 & & \\
                    &\ddots&  \\
                    & & 0   \\
                 \end{smallmatrix}
             \right),
\widetilde{\mathbf{A}}_2^\pi=\left(
               \begin{smallmatrix}
                 0 &  & &\tilde{A}^\top(t+\delta)\pi_1& \cdots& \tilde{A}^\top(t+\delta)\pi_K & \tilde{A}^\top(t+\delta)\pi_1& \cdots& \tilde{A}^\top(t+\delta)\pi_K\\
                   & \ddots &   & \vdots  &   & \vdots  & \vdots  & &\vdots  \\
                   &   & 0 & \tilde{A}^\top(t+\delta)\pi_1&\cdots & \tilde{A}^\top(t+\delta)\pi_K &  \tilde{A}^\top(t+\delta)\pi_1&\cdots & \tilde{A}^\top(t+\delta)\pi_K \\
                  0 &   &   & 0 &   &   & 0  &   &   \\
                   & \ddots &   &   & \ddots &   &   & \ddots  &  \\
                   &   & 0  &   &   & 0 &   &   & 0 \\
                   0 &  & &\tilde{A}^\top(t+\delta)\pi_1& \cdots& \tilde{A}^\top(t+\delta)\pi_K & \tilde{A}^\top(t+\delta)\pi_1& \cdots& \tilde{A}^\top(t+\delta)\pi_K\\
                   & \ddots &   & \vdots  &   & \vdots  & \vdots  & &\vdots  \\
                   &   & 0 & \tilde{A}^\top(t+\delta)\pi_1&\cdots & \tilde{A}^\top(t+\delta)\pi_K &  \tilde{A}^\top(t+\delta)\pi_1&\cdots & \tilde{A}^\top(t+\delta)\pi_K \\
               \end{smallmatrix}
             \right),
\mathbf{Q}= \left(
               \begin{smallmatrix}
                 \mathbb{Q} & & \\
                    &\ddots&  \\
                    & & \mathbb{Q}   \\
                   \mathbb{Q} & & \\
                    &\ddots&  \\
                    & & \mathbb{Q}   \\
                     0 & & \\
                    &\ddots&  \\
                    & & 0   \\
                 \end{smallmatrix}
             \right),
\mathbf{S}^\pi= \left(
               \begin{smallmatrix}
                 \mathbb{S}\pi_1 &\cdots & \mathbb{S}\pi_K\\
                   \vdots & & \vdots \\
              \mathbb{S}\pi_1 &\cdots & \mathbb{S}\pi_K\\
                   0 &\cdots &0 \\
                   \vdots & & \vdots \\
                   0 & \cdots & 0   \\
                  \mathbb{S}\pi_1 &\cdots & \mathbb{S}\pi_K\\
                   \vdots & & \vdots \\
              \mathbb{S}\pi_1 &\cdots & \mathbb{S}\pi_K\\
                 \end{smallmatrix}
             \right),\\
&\mathbf B_5^\pi=\left(
                         \begin{smallmatrix}
                           0 & \cdots & 0 & \mathbb{B}_{14}^1\pi_1 & \cdots & \mathbb{B}_{14}^1\pi_K  & \mathbb{B}_{14}^1\pi_1 & \cdots & \mathbb{B}_{14}^1\pi_K \\
                           \vdots &   & \vdots &\vdots &   & \vdots  &\vdots &   & \vdots   \\
                           0 & \cdots & 0 & \mathbb{B}_{14}^K\pi_1 & \cdots & \mathbb{B}_{14}^K\pi_K  & \mathbb{B}_{14}^K\pi_1 & \cdots & \mathbb{B}_{14}^K\pi_K \\
                         \end{smallmatrix}
                       \right),
\mathbf{G}=\left(
               \begin{smallmatrix}
                 G & & \\
                    &\ddots&  \\
                    & & G  \\
                   G & & \\
                    &\ddots&  \\
                    & & G   \\
                     0 & & \\
                    &\ddots&  \\
                    & & 0   \\
                 \end{smallmatrix}
             \right),
\mathbf{G}^\pi= \left(
               \begin{smallmatrix}
                 \mathbb{G}\pi_1 &\cdots & \mathbb{G}\pi_K\\
                   \vdots & & \vdots \\
              \mathbb{G}\pi_1 &\cdots & \mathbb{G}\pi_K\\
                   0 &\cdots &0 \\
                   \vdots & & \vdots \\
                   0 & \cdots & 0   \\
                  \mathbb{G}\pi_1 &\cdots & \mathbb{G}\pi_K\\
                   \vdots & & \vdots \\
              \mathbb{G}\pi_1 &\cdots & \mathbb{G}\pi_K\\
                 \end{smallmatrix}
             \right),
\mathbf{W}=\left(
               \begin{smallmatrix}
                 W^{(1)}  \\
                    \vdots  \\
                 W^{(K)}   \\
                 \end{smallmatrix}
             \right),
\Xi=\left(
               \begin{smallmatrix}
                 \xi^{(1)}  \\
                    \vdots  \\
                 \xi^{(K)}   \\
                 \end{smallmatrix}
             \right),
X^0=\left(
               \begin{smallmatrix}
                 x^0  \\
                    \vdots  \\
                 x^0  \\
                 \end{smallmatrix}
             \right).
\end{aligned}\end{equation*}\normalsize

In the following, we will use the discounting method of \cite{PT1999} to study the global solvability of MF-AFBSDDE \eqref{eq23}. To start, we first give some results for general nonlinear forward-backward system
\small
\begin{equation}\label{eq24}\left\{\begin{aligned}
&dX(t)=b\Big(t,X(t),X(t-\delta),\mathbb EX(t-\delta),Y(t),Y(t-\theta)I_{[\theta,T]},\mathbb EY(t),\mathbb EY(t-\theta)I_{[\theta,T]},\mathbb E^{\mathcal F_t}[Y(t+\theta)]I_{[0,T-\theta]},\\
&\qquad\qquad\mathbb EY(t+\theta)I_{[0,T-\theta]},\mathbb E^{\mathcal F_{t-\theta}}[Y(t)]I_{[\theta,T]},Z(t),Z(t-\theta)I_{[\theta,T]},\mathbb EZ(t),\mathbb EZ(t-\theta)I_{[\theta,T]},\mathbb E^{\mathcal F_t}[Z(t+\theta)]I_{[0,T-\theta]},\\
&\qquad\qquad\mathbb E^{\mathcal F_{t-\theta}}[Z(t)]I_{[\theta,T]}\Big)dt+\sigma\Big(t,Y(t),Y(t-\theta)I_{[\theta,T]},\mathbb E^{\mathcal F_t}[Y(t+\theta)]I_{[0,T-\theta]},\mathbb E^{\mathcal F_{t-\theta}}[Y(t)]I_{[\theta,T]},\\
&\qquad\qquad\mathbb EY(t+\theta)I_{[0,T-\theta]},\mathbb EYI_{[\theta,T]}, Z(t), Z(t-\theta)I_{[\theta,T]},\mathbb E^{\mathcal F_t}[Z(t+\theta)]I_{[0,T-\theta]},\mathbb E^{\mathcal F_{t-\theta}}[Z(t)]I_{[\theta,T]}\Big)d\mathbf{W}(t),\\
&dY(t)=-f\Big(t,X(t),\mathbb EX(t),Y(t),\mathbb EY(t+\delta)I_{[0,T-\delta]},\mathbb E^{\mathcal F_t}[Y(t+\delta)]I_{[0,T-\delta]}\Big)dt+Z(t)d\mathbf{W}(t),\ t\in[0,T],\\
&X(0)=\Xi,\  Y(T)=\Phi(X(T),\mathbb EX(T)),\ X(t)=X^0(t),\ t\in[-\delta,0),%\ Y(t)=0,\ Z(t)=0,\ t\in[-\theta,0)\cup(T,T+\delta\vee\theta],
\end{aligned}\right.\end{equation}\normalsize
where the coefficients satisfy the following conditions
\begin{description}
  \item[(H1)] There exist $\rho_1,\rho_2\in \mathbb R$ and positive constants $k_0,k'_0,k_i,\ i=1,\cdots,31$ such that for $t$ and  all coefficients a.s.
      \begin{enumerate}
        \item $\langle b(t,x_1,x^\delta,\bar x^\delta,y,y^\theta,\bar y,\bar y^\theta,\bar y^{+\theta},\bar y^{\theta+},\bar y^{-\theta},z,z^\theta,\bar z,\bar z^\theta,\bar z^{+\theta},\bar z^{-\theta})$\\
               $-b(t,x_2,x^\delta,\bar x^\delta,y,y^\theta,\bar y,\bar y^\theta,\bar y^{+\theta},\bar y^{\theta+},\bar y^{-\theta},z,z^\theta,\bar z,\bar z^\theta,\bar z^{+\theta},\bar z^{-\theta}),x_1-x_2\rangle\leq \rho_1|x_1-x_2|^2,$
        \item $|b(t,x,x^\delta_1,\bar x^\delta_1,y_1,y_1^\theta,\bar y_1,\bar y^\theta_1,\bar y^{+\theta}_1,\bar y^{\theta+}_1,\bar y^{-\theta}_1,z_1,z_1^\theta,\bar z_1,\bar z^\theta_1,\bar z^{+\theta}_1,\bar z^{-\theta}_1)-b(t,x,x^\delta_2,\bar x^\delta_2,y_2,y_2^\theta,\bar y_2,$\\
            $\bar y^\theta_2,\bar y^{+\theta}_2,\bar y^{\theta+}_2,\bar y^{-\theta}_2,z_2,z_2^\theta,\bar z_2,\bar z^\theta_2,\bar z^{+\theta}_2,\bar z^{-\theta}_2)|\leq k_1|x^\delta_1-x^\delta_2|+k_2|\bar x^\delta_1-\bar x^\delta_2|+k_3|y_1-y_2|+k_4|y^\theta_1-y^\theta_2|+k_5|\bar y_1-\bar y_2|+k_6|\bar y^\theta_1-\bar y^\theta_2|+k_7|\bar y^{+\theta}_1-\bar y^{+\theta}_2|+k_8|\bar y^{\theta+}_1-\bar y^{\theta+}_2|+k_9|\bar y^{-\theta}_1-\bar y^{-\theta}_2|            +k_{10}|z_1-z_2|+k_{11}|z^\theta_1-z^\theta_2|+k_{12}|\bar z_1-\bar z_2|+k_{13}|\bar z^\theta_1-\bar z^\theta_2|+k_{14}|\bar z^{+\theta}_1-\bar z^{+\theta}_2|+k_{15}|\bar z^{-\theta}_1-\bar z^{-\theta}_2| ,$
        \item $|b(t,x,x^\delta,\bar x^\delta,y,y^\theta,\bar y,\bar y^\theta,\bar y^{+\theta},\bar y^{\theta+},\bar y^{-\theta},z,z^\theta,\bar z,\bar z^\theta,\bar z^{+\theta},\bar z^{-\theta})|$\\
            $\leq |b(t,0,x^\delta,\bar x^\delta,y,y^\theta,\bar y,\bar y^\theta,\bar y^{+\theta},\bar y^{\theta+},\bar y^{-\theta},z,z^\theta,\bar z,\bar z^\theta,\bar z^{+\theta},\bar z^{-\theta})|+k_0(1+|x|),$
        \item $\langle f(t,x,\bar x,y_1,\bar y^{\delta+},\bar y^{+\delta})-f(t,x,\bar x,y_2,\bar y^{\delta+},\bar y^{+\delta}),y_1-y_2\rangle\leq \rho_2|y_1-y_2|^2,$
        \item $|f(t,x_1,\bar x_1,y,\bar y^{\delta+}_1,\bar y^{+\delta}_1)-f(t,x_2,\bar x_2,y,\bar y^{\delta+}_2,\bar y^{+\delta}_2)|\leq k_{16}|x_1-x_2|+k_{17}|\bar x_1-\bar x_2|+k_{18}|\bar y^{\delta+}_1-\bar y^{\delta+}_2|+k_{19}|\bar y^{+\delta}_1-\bar y^{+\delta}_2|,$
        \item $|f(t,x,\bar x,y,\bar y^{\delta+},\bar y^{+\delta})|\leq|f(t,x,\bar x,0,\bar y^{\delta+},\bar y^{+\delta})|+k'_0(1+|y|),$
        \item $|\sigma(t,y_1, y^\theta_1,\bar y^{+\theta}_1,\bar y^{-\theta}_1,\bar y^{\theta+}_1,\bar y_1,z_1, z^\theta_1,\bar z^{+\theta}_1,\bar z^{-\theta}_1)-\sigma(t,y_2, y^\theta_2,\bar y^{+\theta}_2,\bar y^{-\theta}_2,\bar y^{\theta+}_2,\bar y_2,z_2, z^\theta_2,\bar z^{+\theta}_2,\bar z^{-\theta}_2)|^2$ $\leq
        k_{20}^2|y_1-y_2|^2+k_{21}^2| y^\theta_1- y^\theta_2|^2+k_{22}^2|\bar y^{+\theta}_1-\bar y^{+\theta}_2|^2+k_{23}^2|\bar y^{-\theta}_1-\bar y^{-\theta}_2|^2+k_{24}^2|\bar y^{\theta+}_1-\bar y^{\theta+}_2|^2+k_{25}^2|\bar y_1-\bar y_2|^2+k_{26}^2|z_1-z_2|^2+k_{27}^2|z^\theta_1-z^\theta_2|^2+k_{28}^2|\bar z^{+\theta}_1-\bar z^{+\theta}_2|^2+k_{29}^2|\bar z^{-\theta}_1-\bar z^{-\theta}_2|^2,$
        \item $|\Phi(x_1,\bar x_1)-\Phi(x_2,\bar x_2)|^2\leq k_{30}^2|x_1-x_2|^2+k_{31}^2|\bar x_1-\bar x_2|^2.$
      \end{enumerate}
  \item[(H2)] $\mathbb E\int_0^T\Big(|b(s,\textbf{0})|^2+|\sigma(s,\textbf{0})|^2+|f(s,\textbf{0})|^2\Big)ds+\mathbb E|\Phi(0,0)|^2<+\infty.$
  \end{description}

Let $\mathcal H$ be a Hilbert space. Recall that $L^2_{\mathcal F}(0,T;\mathcal H)$ denotes the space of $\mathcal H-$valued $\{\mathcal F_s\}-$ progressively measurable processes $\{v(s),s\in[0,T]\}$ such that $\|v\|^2:=\mathbb E\int_0^T|v(s)|^2ds<\infty$. Then for $\rho \in \mathbb R$, we define an equivalent norm on $L^2_{\mathcal F}(0,T;\mathcal H)$:
$$\|v\|_{\rho}:=\Big(\mathbb E\int_0^Te^{-\rho s}|v(s)|^2ds\Big)^{\frac{1}{2}}.$$

Define
\begin{equation}\nonumber\begin{aligned}
L_{\rho,k,\delta}:=2(\rho_1+\rho_2)+k_1+k_2+k_{18}+k_{19}+(k_1+k_2)e^{-(2\rho_1+k_1+k_2)\delta}+(k_{18}+k_{19})e^{-(2\rho_2+k_{18}+k_{19})\delta}.
\end{aligned}\end{equation}
Then we have the following theorem.
\begin{theorem}\label{theo3.2}
Suppose that assumptions \textbf{(H1)} and \textbf{(H2)} hold. Then there exists a $\epsilon_0>0$, which depends on $k_i,\ \rho_1,\ \rho_2,\ T$, for $i=1,2,16,17,18,19,30,31$ such that when $k_j\in[0,\epsilon_0),$ for $j=3,\cdots,15,20,\cdots,29$,  there exists a unique adapted solution $(X,Y,Z)\in L^2_{\mathcal F}(-\delta,T;\mathbb R^n)\times L^2_{\mathcal F}(0,T;\mathbb R^m)\times L^2_{\mathcal F}(0,T;\mathbb R^{m\times d})$ to MF-AFBSDDE \eqref{eq24}. Further, if $L_{\rho,k,\delta}< 0$, there exists a $\epsilon_1>0$, which depends on $k_i,\ \rho_1,\ \rho_2$, for $i=1,2,16,17,18,19,30,31$ and is independent of $T$, such that when $k_j\in[0,\epsilon_1),$ for $j=3,\cdots,15,20,\cdots,29$, there exists a unique adapted solution $(X,Y,Z)$ to MF-AFBSDDE \eqref{eq24}.
\end{theorem}
Before proving Theorem \ref{theo3.2}, we should do some preparing work.
For any given $(Y(\cdot),Z(\cdot))\in L^2_{\mathcal F}(0,T;\mathbb R^m)\times L^2_{\mathcal F}(0,T;\mathbb R^{m\times d})$, the forward equation in the MF-AFBSDDE \eqref{eq24} admits a unique solution $X(\cdot)\in L^2_{\mathcal F}(-\delta,T;\mathbb R^n)$. Thus we introduce a map $\mathcal{M}_1:L^2_{\mathcal F}(0,T;\mathbb R^m)\times L^2_{\mathcal F}(0,T;\mathbb R^{m\times d})\rightarrow L^2_{\mathcal F}(-\delta,T;\mathbb R^n)$, through
\begin{equation}\label{eq25}\left\{\begin{aligned}
&X(t)=\Xi+\int_0^tb\Big(s,X(s),X(s-\delta),\mathbb EX(s-\delta),Y(s),Y(s-\theta)I_{[\theta,T]},\mathbb EY(s),\mathbb EY(s-\theta)I_{[\theta,T]},\\
&\qquad\quad\mathbb E^{\mathcal F_s}[Y(s+\theta)]I_{[0,T-\theta]},\mathbb EY(s+\theta)I_{[0,T-\theta]},\mathbb E^{\mathcal F_{s-\theta}}[Y(s)]I_{[\theta,T]},Z(s),Z(s-\theta)I_{[\theta,T]},\mathbb EZ(s),\\
&\qquad\quad\mathbb EZ(s-\theta)I_{[\theta,T]},\mathbb E^{\mathcal F_s}[Z(s+\theta)]I_{[0,T-\theta]},\mathbb E^{\mathcal F_{s-\theta}}[Z(s)]I_{[\theta,T]}\Big)ds+\int_0^t\sigma\Big(s,Y(s),Y(s-\theta)I_{[\theta,T]},\\
&\qquad\quad\mathbb E^{\mathcal F_s}[Y(s+\theta)]I_{[0,T-\theta]},\mathbb E^{\mathcal F_{s-\theta}}[Y(s)]I_{[\theta,T]},\mathbb EY(s+\theta)I_{[0,T-\theta]},\mathbb EY(s)I_{[\theta,T]},Z(s),Z(s-\theta)I_{[\theta,T]},\\
&\qquad\quad\mathbb E^{\mathcal F_s}[Z(s+\theta)]I_{[0,T-\theta]},\mathbb E^{\mathcal F_{s-\theta}}[Z(s)]\Big)d\mathbf{W}(s),\ t\in[0,T],\\
&X(t)=x(t),\ t\in[-\delta,0).
\end{aligned}\right.\end{equation}
Indeed the well-posedness of \eqref{eq25} can be established by applying the contraction mapping method under the assumptions \textbf{(H1)} and \textbf{(H2)}, although the term $\mathbb EX(s-\delta)$ is involved. We omit the proof here. By BDG inequality and the delay theory, it follows that $\mathbb E\sup_{t\in[-\delta,T]}|X(t)|^2<\infty$.
\begin{lemma}\label {l1}
Let $X_i$ be the solution of \eqref{eq25} corresponding to $(Y_i(\cdot),Z_i(\cdot))\in L^2_{\mathcal F}(0,T;\mathbb R^m)\times L^2_{\mathcal F}(0,T;\mathbb R^{m\times d})$, $i=1,2$, respectively. Then for all $\rho\in \mathbb R$ and some constants $l_j>0$, $j=1,\cdots,13$, we have
\begin{equation}\label{eq26}\begin{aligned}
&e^{-\rho t}\mathbb E|\widehat{X}(t)|^2+\bar \rho_1 \mathbb E\int_0^te^{-\rho s}|\widehat{X}(s)|^2ds\leq \Big(k_3l_1+k_5l_3+k_9l_7+k_{20}^2+k_{23}^2+k_{25}^2+(k_4l_2+k_6l_4\\
&\quad+k_{21}^2)  e^{-\rho\theta}\Big)\mathbb E\int_0^te^{-\rho s}|\widehat{Y}(s)|^2ds+(k_7l_5+k_8l_6+k_{22}^2+k_{24}^2)e^{\rho\theta}\mathbb E\int_0^{t+\theta}e^{-\rho s}|\widehat{Y}(s)|^2I_{[0,T-\theta]}ds\\
&\quad+\Big(k_{10}l_8+k_{12}l_{10}+k_{15}l_{13}+k_{26}^2+k_{29}^2+(k_{11}l_9+k_{13}l_{11}+k_{27}^2) e^{-\rho\theta}\Big)\mathbb E\int_0^te^{-\rho s}|\widehat{Z}(s)|^2ds\\
&\quad+(k_{14}l_{12}+k_{28}^2)e^{\rho\theta}\mathbb E\int_0^{t+\theta}e^{-\rho s}|\widehat{Z}(s)|^2I_{[0,T-\theta]}ds,
\end{aligned}\end{equation}
where $\bar\rho_1=\rho -2\rho_1-(k_1+k_2)(1+e^{-\rho\delta})-\sum_{j=1}^{13} k_{j+2}l_j^{-1}$, and $\widehat\psi=\psi_1-\psi_2,\ \psi=X,Y,Z$. We also have
\begin{equation}\label{eq27}\begin{aligned}
&e^{-\rho t}\mathbb E|\widehat{X}(t)|^2\\
\leq&
\Big(k_3l_1+k_5l_3+k_9l_7+k_{20}^2+k_{23}^2+k_{25}^2+(k_4l_2+k_6l_4+k_{21}^2) e^{-(\rho-\tilde\rho_1)\theta}\Big)\mathbb E\int_0^te^{-\tilde\rho_1(t-s) -\rho s}|\widehat{Y}(s)|^2ds\\
&+(k_7l_5+k_8l_6+k_{22}^2+k_{24}^2)e^{(\rho-\tilde\rho_1)\theta}\mathbb E\int_0^{t+\theta}e^{-\tilde\rho_1(t-s) -\rho s}|\widehat{Y}(s)|^2I_{[0,T-\theta]}ds\\
&+\Big(k_{10}l_8+k_{12}l_{10}+k_{15}l_{13}+k_{26}^2+k_{29}^2+(k_{11}l_9+k_{13}l_{11}+k_{27}^2) e^{-(\rho-\tilde\rho_1)\theta}\Big)\mathbb E\int_0^te^{-\tilde\rho_1(t-s) -\rho s}|\widehat{Z}(s)|^2ds\\
&+(k_{14}l_{12}+k_{28}^2)e^{(\rho-\tilde\rho_1)\theta}\mathbb E\int_0^{t+\theta}e^{-\tilde\rho_1(t-s) -\rho s}|\widehat{Z}(s)|^2I_{[0,T-\theta]}ds,
\end{aligned}\end{equation}
where $\tilde\rho_1$ solves the equation
$\tilde\rho_1=\rho -2\rho_1-(k_1+k_2)(1+e^{-(\rho-\tilde\rho_1)\delta})-\sum_{j=1}^{13} k_{j+2}l_j^{-1}.
$
Moreover,
\begin{equation}\label{eq28}\begin{aligned}
\|\widehat X\|_{\rho}^2\leq& \frac{e^{\tilde\rho_1\theta}-e^{-\tilde\rho_1T}}{\tilde\rho_1}\Big[\Big(k_3l_1+k_5l_3+k_9l_7+k_{20}^2+k_{23}^2+k_{25}^2+(k_4l_2+k_6l_4+k_{21}^2) e^{-(\rho-\tilde\rho_1)\theta}\\
&+(k_7l_5+k_8l_6+k_{22}^2+k_{24}^2)e^{(\rho-\tilde\rho_1)\theta}\Big)\|\widehat Y\|_{\rho}^2+\Big(k_{10}l_8+k_{12}l_{10}+k_{15}l_{13}+k_{26}^2+k_{29}^2\\
&+(k_{11}l_9+k_{13}l_{11}+k_{27}^2) e^{-(\rho-\tilde\rho_1)\theta}+(k_{14}l_{12}+k_{28}^2)e^{(\rho-\tilde\rho_1)\theta}\Big)\|\widehat Z\|_{\rho}^2\Big],
\end{aligned}\end{equation}
and
\begin{equation}\label{eq29}\begin{aligned}
&e^{-\rho T}\mathbb E|\widehat{X}(T)|^2\\
\leq& \big(1\vee e^{-\tilde\rho_1T}\big)\Big[\Big(k_3l_1+k_5l_3+k_9l_7+k_{20}^2+k_{23}^2+k_{25}^2+(k_4l_2+k_6l_4+k_{21}^2)e^{-(\rho-\tilde\rho_1)\theta}+(k_7l_5\\
&+k_8l_6+k_{22}^2+k_{24}^2)e^{(\rho-\tilde\rho_1)\theta}\Big)\|\widehat Y\|_{\rho}^2+\Big(k_{10}l_8+k_{12}l_{10}+k_{15}l_{13}+k_{26}^2+k_{29}^2\\
&+(k_{11}l_9+k_{13}l_{11}+k_{27}^2) e^{-(\rho-\tilde\rho_1)\theta}+(k_{14}l_{12}+k_{28}^2)e^{(\rho-\tilde\rho_1)\theta}\Big)\|\widehat Z\|_{\rho}^2\Big].
\end{aligned}\end{equation}
Specifically, if $\tilde\rho_1>0$,
\begin{equation}\label{eq30}\begin{aligned}
&e^{-\rho T}\mathbb E|\widehat{X}(T)|^2\\
\leq& \Big(k_3l_1+k_5l_3+k_9l_7+k_{20}^2+k_{23}^2+k_{25}^2+(k_4l_2+k_6l_4+k_{21}^2) e^{-(\rho-\tilde\rho_1)\theta}\\
&\ \ +(k_7l_5+k_8l_6+k_{22}^2+k_{24}^2)e^{(\rho-\tilde\rho_1)\theta}\Big)\|\widehat Y\|_{\rho}^2+\Big(k_{10}l_8+k_{12}l_{10}+k_{15}l_{13}+k_{26}^2+k_{29}^2\\
&\ \ +(k_{11}l_9+k_{13}l_{11}+k_{27}^2) e^{-(\rho-\tilde\rho_1)\theta}+(k_{14}l_{12}+k_{28}^2)e^{(\rho-\tilde\rho_1)\theta}\Big)\|\widehat Z\|_{\rho}^2.
\end{aligned}\end{equation}
\end{lemma}

Please refer Appendix A for the proof.

Similarly, for a given $X(\cdot)\in L^2_{\mathcal F}(-\delta,T;\mathbb R^n)$, the backward equation in the MF-AFBSDDE \eqref{eq24} admits a unique solution
$(Y(\cdot),Z(\cdot))\in L^2_{\mathcal F}(0,T;\mathbb R^m)\times L^2_{\mathcal F}(0,T;\mathbb R^{m\times d})$. Thus we introduce a map $\mathcal{M}_2:L^2_{\mathcal F}(-\delta,T;\mathbb R^n)\rightarrow L^2_{\mathcal F}(0,T;\mathbb R^m)\times L^2_{\mathcal F}(0,T;\mathbb R^{m\times d})$, through
\begin{equation}\label{eq33}\begin{aligned}
&dY(t)=\Phi(X(T),\mathbb EX(T))+\int_t^Tf\big(t,X(t),\mathbb EX(t),Y(t),\mathbb EY(t+\delta)I_{[0,T-\delta]},\\
&\qquad\qquad\qquad\mathbb E^{\mathcal F_t}[Y(t+\delta)]I_{[0,T-\delta]}\big)ds-\int_t^TZ(s)d\mathbf{W}(s)
%,\\
%&Y(t)=0,\ Z(t)=0,\ t\in[-\theta,0)\cup(T,T+\delta\vee\theta]
.
\end{aligned}\end{equation}
The well-posedness of \eqref{eq33} under the assumptions \textbf{(H1)} and \textbf{(H2)} can be established by \cite{BLP2009,DP1997,PY09}. The proof is also omitted. Moreover, we have $\mathbb E\sup_{t\in[0,T]}|Y(t)|^2<\infty$.
\begin{lemma}\label {l2}
Let $(Y_i(\cdot),Z_i(\cdot))$ be the solution of \eqref{eq33} corresponding to $X(\cdot)\in L^2_{\mathcal F}(-\delta,T;\mathbb R^n)$, $i=1,2$, respectively. Then for all $\rho\in \mathbb R$ and some constants $l_j>0$, $j=14,15$, we have
\begin{equation}\label{eq34}\begin{aligned}
&e^{-\rho t}\mathbb E|\widehat{Y}(t)|^2+\bar \rho_2 \mathbb E\int_t^Te^{-\rho s}|\widehat{Y}(s)|^2ds+\mathbb E\int_t^Te^{-\rho s}|\widehat{Z}(s)|^2ds\\
\leq &(k_{30}^2+k_{31}^2)e^{-\rho T}\mathbb E|\widehat{X}(T)|^2+(k_{16}l_{14}+k_{17}l_{15})\mathbb E\int_t^Te^{-\rho s}|\widehat{X}(s)|^2ds,
\end{aligned}\end{equation}
where $\bar\rho_2=-\rho -2\rho_2-k_{16}l_{14}^{-1}-k_{17}l_{15}^{-1}-(k_{18}+k_{19})(1+e^{\rho\delta})$, and $\widehat\psi=\psi_1-\psi_2,\ \psi=X,Y,Z$. We also have that
\begin{equation}\label{eq35}\begin{aligned}
&e^{-\rho t}\mathbb E|\widehat{Y}(t)|^2+\mathbb E\int_t^Te^{-\tilde\rho_2(s-t) -\rho s}|\widehat{Z}(s)|^2ds\\
\leq&(k_{30}^2+k_{31}^2)e^{-\tilde\rho_2(T-t)-\rho T}\mathbb E|\widehat{X}(T)|^2+(k_{16}l_{14}+k_{17}l_{15})\mathbb E\int_t^Te^{-\tilde\rho_2(s-t) -\rho s}|\widehat{X}(s)|^2ds,
\end{aligned}\end{equation}
where  $\tilde\rho_2$ solves the equation $\tilde\rho_2=-\rho -2\rho_2-k_{16}l_{14}^{-1}-k_{17}l_{15}^{-1}-(k_{18}+k_{19})(1+e^{(\rho+\tilde\rho_2)\delta})$. Moreover,
\begin{equation}\label{eq366}\begin{aligned}
\|\widehat Y\|_{\rho}^2\leq& \frac{1-e^{-\tilde\rho_2T}}{\tilde\rho_2}\big[(k_{30}^2+k_{31}^2)e^{-\rho T}\mathbb E|\widehat{X}(T)|^2+(k_{16}l_{14}+k_{17}l_{15})\|\widehat X\|_{\rho}^2\big],
\end{aligned}\end{equation}
and
\begin{equation}\label{eq377}\begin{aligned}
\|\widehat Z\|_{\rho}^2\leq& \frac{(k_{30}^2+k_{31}^2)e^{-\tilde\rho_2 T}e^{-\rho T}\mathbb E|\widehat{X}(T)|^2+(k_{16}l_{14}+k_{17}l_{15})(1\vee e^{-\tilde\rho_2T})\|\widehat X\|_{\rho}^2}{1\wedge e^{-\tilde\rho_2T}}.
\end{aligned}\end{equation}
Specifically, if $\tilde\rho_2>0$,
\begin{equation}\label{eq388}\begin{aligned}
\|\widehat Z\|_{\rho}^2\leq& (k_{30}^2+k_{31}^2)e^{-\rho T}\mathbb E|\widehat{X}(T)|^2+(k_{16}l_{14}+k_{17}l_{15})\|\widehat X\|_{\rho}^2.
\end{aligned}\end{equation}
\end{lemma}

The proof is given in Appendix B. Based on Lemmas \ref{l1} and \ref{l2}, we present the proof of Theorem \ref{theo3.2}.

\emph{Proof of Theorem \ref{theo3.2}}. Define $\mathcal M:=\mathcal M_2 \circ\mathcal M_1$. Since $\mathcal M_1$ is defined by \eqref{eq25} and $\mathcal M_2$ is defined by \eqref{eq33}. Therefore $\mathcal M$ maps $L^2_{\mathcal F}(0,T;\mathbb R^m)\times L^2_{\mathcal F}(0,T;\mathbb R^{m\times d})$ into itself. To prove the theorem, we only need to show that $\mathcal M$ is a contraction mapping for some equivalent norm $\|\cdot\|_{\rho}$. For $(Y_i,Z_i)\in L^2_{\mathcal F}(0,T;\mathbb R^m)\times L^2_{\mathcal F}(0,T;\mathbb R^{m\times d})$, let $X_i:=\mathcal M_1(Y_i,Z_i)$ and $(\overline{Y}_i,\overline{Z}_i):=\mathcal M((Y_i,Z_i))$; from \eqref{eq28}-\eqref{eq29}, \eqref{eq366}-\eqref{eq377}, we have
\begin{equation}\nonumber\begin{aligned}
\|\widehat Y\|_{\rho}^2+\|\widehat Z\|_{\rho}^2\leq& \Big[\frac{1-e^{-\tilde\rho_2T}}{\tilde\rho_2}+\frac{1\vee e^{-\tilde\rho_2T}}{1\wedge e^{-\tilde\rho_2T}}\Big]\\
&\times\big[(k_{30}^2+k_{31}^2)e^{-\rho T}\mathbb E|\widehat{X}(T)|^2+(k_{16}l_{14}+k_{17}l_{15})\|\widehat X\|_{\rho}^2\big]\\
\leq & \Big[\frac{1-e^{-\tilde\rho_2T}}{\tilde\rho_2}+\frac{1\vee e^{-\tilde\rho_2T}}{1\wedge e^{-\tilde\rho_2T}}\Big]\\
&\times\Big[(k_{30}^2+k_{31}^2)(1\vee e^{-\tilde\rho_1T})+(k_{16}l_{14}+k_{17}l_{15})\frac{e^{\tilde\rho_1\theta}-e^{-\tilde\rho_1T}}{\tilde\rho_1}\Big]\\
&\times\Big[\Big(k_3l_1+k_5l_3+k_9l_7+k_{20}^2+k_{23}^2+k_{25}^2+(k_4l_2+k_6l_4+k_{21}^2) e^{-(\rho-\tilde\rho_1)\theta}\\
&+(k_7l_5+k_8l_6+k_{22}^2+k_{24}^2)e^{(\rho-\tilde\rho_1)\theta}\Big)\|\widehat Y\|_{\rho}^2+\Big(k_{10}l_8+k_{12}l_{10}+k_{15}l_{13}+k_{26}^2+k_{29}^2\\
&+(k_{11}l_9+k_{13}l_{11}+k_{27}^2) e^{-(\rho-\tilde\rho_1)\theta}+(k_{14}l_{12}+k_{28}^2)e^{(\rho-\tilde\rho_1)\theta}\Big)\|\widehat Z\|_{\rho}^2\Big].
\end{aligned}\end{equation}
Recall that $\tilde\rho_1$ satisfies $\tilde\rho_1-\rho+(k_1+k_2)e^{(\tilde\rho_1-\rho)\delta}=-2\rho_1-(k_1+k_2)-\sum_{j=1}^{13} k_{j+2}l_j^{-1}$, and $\tilde\rho_2$ satisfies  $\tilde\rho_2+\rho+(k_{18}+k_{19})e^{(\tilde\rho_2+\rho)\delta}= -2\rho_2-k_{16}l_{14}^{-1}-k_{17}l_{15}^{-1}-(k_{18}+k_{19})$. By choosing suitable $\rho$ (e.g. choosing $l_j=k_{j+2}$, for $j=1,...,13$ and $\rho$ big enough such that $\tilde\rho_1\ge 1$ and $\tilde\rho_2\leq 0$), the first assertion of Theorem \ref{theo3.2} is obtained.

If $L_{\rho,k,\delta}< 0$, we can choose $\rho\in \mathbb R$ and sufficiently large $l_j\ (j=1,\cdots,15)$ such that $$\tilde\rho_1>0,\quad \tilde\rho_2>0.$$ In fact, recall that $\tilde\rho_1$ satisfies $\tilde\rho_1-\rho+(k_1+k_2)e^{(\tilde\rho_1-\rho)\delta}=-2\rho_1-(k_1+k_2)-\sum_{j=1}^{13} k_{j+2}l_j^{-1}$, and $\tilde\rho_2$ satisfies  $\tilde\rho_2+\rho+(k_{18}+k_{19})e^{(\tilde\rho_2+\rho)\delta}= -2\rho_2-k_{16}l_{14}^{-1}-k_{17}l_{15}^{-1}-(k_{18}+k_{19})$.
Denote $\lambda_1:=\tilde\rho_1-\rho$ and $\lambda_2:=\tilde\rho_2+\rho$, then
$$\lambda_1+(k_1+k_2)e^{\lambda_1\delta}=-2\rho_1-k_1-k_2-\sum_{j=1}^{13} k_{j+2}l_j^{-1}$$ and
$$\lambda_2+(k_{18}+k_{19})e^{\lambda_2\delta}=-2\rho_2-k_{16}l_{14}^{-1}-k_{17}l_{15}^{-1}-k_{18}-k_{19},$$
respectively (concerning the existence of $\lambda_1$, $\lambda_2$, see the proof of Lemma \ref{l1} and  Lemma \ref{l2}). Then it is obvious that $\lambda_1<-2\rho_1-k_1-k_2-\sum_{j=1}^{13} k_{j+2}l_j^{-1}$ and $\lambda_2<-2\rho_2-k_{16}l_{14}^{-1}-k_{17}l_{15}^{-1}-k_{18}-k_{19}$. Thus
\[
\begin{aligned}
-\lambda_1-\lambda_2=&2(\rho_1+\rho_2)+k_1+k_2+k_{18}+k_{19}+(k_1+k_2)e^{\lambda_1\delta}+(k_{18}+k_{19})e^{\lambda_2\delta}\\
\leq&2(\rho_1+\rho_2)+k_1+k_2+k_{18}+k_{19}+(k_1+k_2)e^{(-2\rho_1-k_1-k_2-\sum_{j=1}^{13} k_{j+2}l_j^{-1})\delta}\\
&+(k_{18}+k_{19})e^{(-2\rho_2-k_{16}l_{14}^{-1}-k_{17}l_{15}^{-1}-k_{18}-k_{19})\delta},
\end{aligned}
\]
and by recalling that $L_{\rho,k,\delta}< 0$, one can choose sufficiently large $l_j\ (j=1,\cdots,15)$ such that $-\lambda_1<\lambda_2$. Therefore, by recalling $\tilde\rho_1=\rho+\lambda_1$ and  $\tilde\rho_2=-\rho+\lambda_2$, we can choose suitable $\rho$ (e.g. $\rho:=-\lambda_1+\frac{1}{2}(\lambda_1+\lambda_2)$) such that $\tilde\rho_1>0,\  \tilde\rho_2>0.$

Then from \eqref{eq28}, \eqref{eq30}, \eqref{eq366}, \eqref{eq388}, we have
\begin{equation}\nonumber\begin{aligned}
\|\widehat Y\|_{\rho}^2+\|\widehat Z\|_{\rho}^2\leq& \frac{1+\tilde\rho_2}{\tilde\rho_2}\big[(k_{30}^2+k_{31}^2)e^{-\rho T}\mathbb E|\widehat{X}(T)|^2+(k_{16}l_{14}+k_{17}l_{15})\|\widehat X\|_{\rho}^2\big]\\
\leq & \frac{1+\tilde\rho_2}{\tilde\rho_2}\Big[k_{30}^2+k_{31}^2+(k_{16}l_{14}+k_{17}l_{15})\frac{e^{\tilde\rho_1\theta}}{\tilde\rho_1}\Big]\\
&\times\Big[\Big(k_3l_1+k_5l_3+k_9l_7+k_{20}^2+k_{23}^2+k_{25}^2+(k_4l_2+k_6l_4+k_{21}^2) e^{-(\rho-\tilde\rho_1)\theta}\\
&+(k_7l_5+k_8l_6+k_{22}^2+k_{24}^2)e^{(\rho-\tilde\rho_1)\theta}\Big)\|\widehat Y\|_{\rho}^2+\Big(k_{10}l_8+k_{12}l_{10}+k_{15}l_{13}+k_{26}^2+k_{29}^2\\
&+(k_{11}l_9+k_{13}l_{11}+k_{27}^2) e^{-(\rho-\tilde\rho_1)\theta}+(k_{14}l_{12}+k_{28}^2)e^{(\rho-\tilde\rho_1)\theta}\Big)\|\widehat Z\|_{\rho}^2\Big].
\end{aligned}\end{equation}
This completes the second assertion of Theorem \ref{theo3.2}.

Let $\rho_1^\ast$ and $\rho_2^\ast$ be the largest eigenvalue of $\frac{1}{2}(\mathbf{A}+\mathbf{A}^\top)$ and $\frac{1}{2}(\mathcal{A}+\mathcal{A}^\top)$. Comparing \eqref{eq24} with \eqref{eq23}, we can check that the coefficients of \textbf{(H1)} can be chosen as follows:
\begin{equation}\nonumber\begin{aligned}
&\rho_1=\rho_1^\ast,\ \rho_2=\rho_2^\ast,\ k_0=\|\mathbf{A}\|,\ k'_0=\|\mathcal{A}\|,\ k_1=\|\hat{\mathbf{A}}\|,\ k_2=\|\widetilde{\mathbf{A}}_1^{\pi}\|,\ k_3=\|\mathbf{B}_1\|,\ k_4=\|\mathbf{B}_3\|,\\ &k_5=\|\mathbf{B}_2^\pi\|,\ k_6=\|\mathbf{B}_1^\pi\|,\ k_7=\|\mathbf{B}_2\|,\ k_8=\|\mathbf{B}_3^\pi\|,\ k_9=\|\mathbf{B}_4\|,\ k_{10}=\|\mathbf{D}_1\|,\ k_{11}=\|\mathbf{D}_3\|,\\
& k_{12}=\|\mathbf{D}_2^\pi\|,\ k_{13}=\|\mathbf{D}_1^\pi\|,\ k_{14}=\|\mathbf{D}_2\|,\ k_{15}=\|\mathbf{D}_4\|,\ k_{16}=\|\mathbf{Q}\|,\ k_{17}=\|\mathbf{S}^\pi\|,\ k_{18}=\|\widetilde{\mathbf{A}}_2^{\pi}\|, \\
&k_{19}=\|\hat{\mathcal{A}}\|,\ k_{20}^2=\|\mathbf{B}_5\|^2,\ k_{21}^2=\|\mathbf{B}_6\|^2,\ k_{22}^2=\|\mathbf{B}_7\|^2,\ k_{23}^2=\|\mathbf{B}_8\|^2,\ k_{24}^2=\|\mathbf{B}_4^\pi\|^2,\ k_{25}^2=\|\mathbf{B}_5^\pi\|^2, \\
&k_{26}^2=\|\mathbf{D}_5\|^2,\ k_{27}^2=\|\mathbf{D}_6\|^2,\ k_{28}^2=\|\mathbf{D}_7\|^2,\ k_{29}^2=\|\mathbf{D}_8\|^2,\ k_{30}^2=\|\mathbf{G}\|^2,\ k_{31}^2=\|\mathbf{G}^\pi\|^2,\\
&L_{\rho,k,\delta}=2(\rho_1^\ast+\rho_2^\ast)+\|\hat{\mathbf{A}}\|+\|\widetilde{\mathbf{A}}_1^{\pi}\|+\|\widetilde{\mathbf{A}}_2^{\pi}\|+\|\hat{\mathcal{A}}\|
+(\|\hat{\mathbf{A}}\|+\|\widetilde{\mathbf{A}}_1^{\pi}\|)e^{-(2\rho_1^\ast+\|\hat{\mathbf{A}}\|+\|\widetilde{\mathbf{A}}_1^{\pi}\|)\delta}\\
&+(\|\widetilde{\mathbf{A}}_2^{\pi}\|+\|\hat{\mathcal{A}}\|)e^{-(2\rho_1^\ast+\|\widetilde{\mathbf{A}}_2^{\pi}\|+\|\hat{\mathcal{A}}\|)\delta}
.
\end{aligned}\end{equation}
Thus by applying Theorem \ref{theo3.2}, we obtain the following global well-posedness of \eqref{eq23}.
\begin{theorem}
Suppose that $L_{\rho,k,\delta}< 0$, then there exists a $\epsilon_1>0$, which depends on $\rho_1^\ast,\ \rho_2^\ast,\ \|\hat{\mathbf{A}}\|$, $\ \|\widetilde{\mathbf{A}}_1^{\pi}\|,\ \|\mathbf{Q}\|,\ \|\mathbf{S}^\pi\|,\ \|\widetilde{\mathbf{A}}_2^{\pi}\|,\ \|\hat{\mathcal{A}}\|,\ \|\mathbf{G}\|,\ \|\mathbf{G}^\pi\|$, and is independent of $T$, such that when $\|\mathbf{B}_j\|,\ \|\mathbf{D}_j\|$ $\ (j=1,\cdots,8),\ \|\mathbf{B}_l^\pi\|\ (l=1,\cdots,5),\ \|\mathbf{D}_1^\pi\|,\ \|\mathbf{D}_2^\pi\|\in[0,\epsilon_1)$, there exists a unique adapted solution $(\alpha,\ \beta,\ \gamma,\ \check y,\ \check z,\ \zeta)$ to (CC) system \eqref{CC}.
\end{theorem}

\section{Asymptotic $\varepsilon$-optimality}\label{asymptotic optimality}

We start this section with some estimates, which play an important role in asymptotic optimality.

\subsection{Agent $\mathcal A_i,\ 1\leq i\leq N$ perturbation}\label{Agent's perturbation}

Let $\widetilde u=(\widetilde u_1,\cdots,\widetilde u_N)$ be the decentralized strategy given by
\begin{equation*}\begin{aligned}\widetilde u_i(t)=&-\big(R_{\theta_i}(t)+\tilde{R}_{\theta_i}(t+\theta)\big)^{-1}\big(B^\top(t) p_i(t)+\hat B^\top(t+\theta)\mathbb E^{\mathcal F_t}[p_i(t+\theta)]I_{[0,T-\theta]}+D^\top(t) q_i(t)\\
&\ \ +\hat{D}^\top(t+\theta)\mathbb E^{\mathcal F_t}[q_i(t+\theta)]I_{[0,T-\theta]}+\Theta_3\big),\ 1\leq i\leq N,
\end{aligned}\end{equation*}
where\footnotesize
\begin{equation*}\left\{\begin{aligned}
&dx_i(t)=\Big[A_{\theta_i}x_i+\hat A_{\theta_i}x_i(t-\delta)+\tilde A\sum_{l=1}^K\pi_l\mathbb E\alpha_l(t-\delta)-B\big(R_{\theta_i}+\tilde{R}_{\theta_i}(t+\theta)\big)^{-1}\big(B^\top p_i+\hat B^\top (t+\theta)\mathbb E^{\mathcal F_t}[p_i(t+\theta)]I_{[0,T-\theta]}\\
&\quad+D^\top q_i+\hat{D}^\top(t+\theta)\mathbb E^{\mathcal F_t}[q_i(t+\theta)]I_{[0,T-\theta]}+\sum_{l=1}^K\pi_l\tilde{B}^\top(t+\theta) \hat{y}_l(t+\theta)I_{[0,T-\theta]}+\sum_{l=1}^K \pi_l\tilde{B}^\top(t+\theta) y_2^l(t+\theta)I_{[0,T-\theta]}\big)\\
&\quad-\hat B\big(R_{\theta_i}(t-\theta)+\tilde{R}_{\theta_i}\big)^{-1}\big(B^\top(t-\theta) p_i(t-\theta)+\hat B^\top \mathbb E^{\mathcal F_{t-\theta}}[p_i]+D^\top(t-\theta) q_i(t-\theta)+\hat{D}^\top\mathbb E^{\mathcal F_{t-\theta}}[q_i]+\sum_{l=1}^K\pi_l\tilde{B}^\top \hat{y}_l\\
&\quad+\sum_{l=1}^K \pi_l\tilde{B}^\top y_2^l\big)I_{[\theta,T]}+\hat Bu^0I_{[0,\theta]}+\tilde B\sum_{l=1}^K\pi_l\mathbb Ev_l(t-\theta)I_{[\theta,T]}+\tilde B u^0(t)I_{[0,\theta]}\Big]dt\\
&\quad-\Big[D\big(R_{\theta_i}+\tilde{R}_{\theta_i}(t+\theta)\big)^{-1}\big(B^\top p_i+\hat B^\top(t+\theta)\mathbb E^{\mathcal F_t}[p_i(t+\theta)]I_{[0,T-\theta]}+D^\top q_i+\hat{D}^\top(t+\theta)\mathbb E^{\mathcal F_t}[q_i(t+\theta)]I_{[0,T-\theta]}\\
&\quad+\sum_{l=1}^K\pi_l\tilde{B}^\top(t+\theta) \hat{y}_l(t+\theta)I_{[0,T-\theta]}+\sum_{l=1}^K \pi_l\tilde{B}^\top(t+\theta) y_2^l(t+\theta)I_{[0,T-\theta]}\big)+\hat D(t)\big(R_{\theta_i}(t-\theta)+\tilde{R}_{\theta_i}(t)\big)^{-1}\\
&\quad\cdot\big(B^\top(t-\theta) p_i(t-\theta)+\hat B^\top \mathbb E^{\mathcal F_{t-\theta}}[p_i]+D^\top(t-\theta) q_i(t-\theta)+\hat{D}^\top\mathbb E^{\mathcal F_{t-\theta}}[q_i]+\sum_{l=1}^K\pi_l\tilde{B}^\top \hat{y}_l+\sum_{l=1}^K \pi_l\tilde{B}^\top y_2^l\big)I_{[\theta,T]}\\
&\quad-\hat Du^0I_{[0,\theta]}\Big]dW_i(t),\ t\in[0,T],\\
&dp_i(t)=-\Big[A_{\theta_i}^\top p_i+\hat A_{\theta_i}^\top(t+\delta)\mathbb E^{\mathcal F_t}[p_i(t+\delta)]I_{[0,T-\delta]}+(Q+\tilde{Q}(t+\delta)) x_i-\mathbb S\sum_{l=1}^K\pi_l\mathbb E\alpha_l\\
&\quad+\sum_{l=1}^K \pi_l\tilde{A}^\top(t+\delta) \hat{y}_l(t+\delta)I_{[0,T-\delta]}+\sum_{l=1}^K \pi_l\tilde{A}^\top(t+\delta) y_2^l(t+\delta)I_{[0,T-\delta]}\Big]dt+q_idW_i(t),\ t\in[0,T],\\
&x_i(0)=\xi_i,\ x_i(t)=x^0(t),\ t\in[-\delta,0),\ \ p_i(T)=G x_i(T)-\mathbb G\sum_{l=1}^K\pi_l\mathbb E\alpha_l(T),%\ p_i(t)=0,\ q_i(t)=0,\ q_{ij}(t)=0,\ t\in(T,T+\delta\vee\theta],
\end{aligned}\right.\end{equation*} \normalsize
with $\alpha_l,\ u_l,\ y_2^l,\ \hat y_l$ given by Theorem \ref{the1}.

Correspondingly, the realized state $(\widetilde x_1,\cdots,\widetilde x_N)$ under the decentralized strategy satisfies
\footnotesize
\begin{equation}\label{eq36}\left\{\begin{aligned}
&d\widetilde x_i(t)=\Big[A_{\theta_i}\widetilde x_i+\hat A_{\theta_i}\widetilde x_i(t-\delta)+\tilde A\widetilde x^{(N)}(t-\delta)-B\big(R_{\theta_i}+\tilde{R}_{\theta_i}(t+\theta)\big)^{-1}\big(B^\top p_i+\hat B^\top (t+\theta)\mathbb E^{\mathcal F_t}[p_i(t+\theta)]I_{[0,T-\theta]}\\
&\quad+D^\top q_i+\hat{D}^\top(t+\theta)\mathbb E^{\mathcal F_t}[q_i(t+\theta)]I_{[0,T-\theta]}+\Theta_3\big)-\hat B\big(R_{\theta_i}(t-\theta)+\tilde{R}_{\theta_i}\big)^{-1}\big(B^\top(t-\theta) p_i(t-\theta)+\hat B^\top \mathbb E^{\mathcal F_{t-\theta}}[p_i]\\
&\quad+D^\top(t-\theta) q_i(t-\theta)+\hat{D}^\top\mathbb E^{\mathcal F_{t-\theta}}[q_i]+\Theta_3(t-\theta)\big)I_{[\theta,T]}+\hat B u^0I_{[0,\theta]}+\tilde B\widetilde u^{(N)}(t-\theta)\Big]dt\\
&\quad-\Big[D\big(R_{\theta_i}+\tilde{R}_{\theta_i}(t+\theta)\big)^{-1}\big(B^\top p_i+\hat B^\top(t+\theta)\mathbb E^{\mathcal F_t}[p_i(t+\theta)]I_{[0,T-\theta]}+D^\top q_i+\hat{D}^\top(t+\theta)\mathbb E^{\mathcal F_t}[q_i(t+\theta)]I_{[0,T-\theta]}+\Theta_3\big)\\
&\quad+\hat D(t)\big(R_{\theta_i}(t-\theta)+\tilde{R}_{\theta_i}(t)\big)^{-1}\big(B^\top(t-\theta) p_i(t-\theta)+\hat B^\top \mathbb E^{\mathcal F_{t-\theta}}[p_i]+D^\top(t-\theta) q_i(t-\theta)+\hat{D}^\top\mathbb E^{\mathcal F_{t-\theta}}[q_i]\\
&\quad+\Theta_3(t-\theta)\big)I_{[\theta,T]}-\hat Du^0I_{[0,\theta]}\Big]dW_i(t),\ t\in[0,T],\\
&\widetilde x_i(0)=\xi_i,\ \widetilde x_i(t)=x^0(t),\ t\in[-\delta,0),\ \widetilde u_i(t)=u^0(t),\ t\in[-\theta,0),
\end{aligned}\right.\end{equation}
\normalsize
and $\widetilde x^{(N)}(\cdot)=\frac{1}{N}\sum_{i=1}^N\widetilde x_i(\cdot)$, $\widetilde u^{(N)}(\cdot)=\frac{1}{N}\sum_{i=1}^N\widetilde u_i(\cdot)$.
For $1\leq j\leq N$, denote the perturbation
$$\delta u_j=u_j-\widetilde u_j,\quad \delta x_j=\breve x_j-\widetilde x_j,\quad \Delta\mathcal J_j=\mathcal J_j(u_i,\widetilde u_{-i})-\mathcal J_j(\widetilde u_i,\widetilde u_{-i}).$$

Similar as the computations in Section \ref{p-b-p optimality}, we have
\begin{equation}\label{eq42}\begin{aligned}
\Delta \mathcal J_{soc}^{(N)}
  =&\mathbb E\Bigg\{\int_0^T\Big[\langle Q\widetilde x_i,\delta x_i\rangle-\langle (QS+ S^\top Q -S^\top QS)\sum_{l=1}^K\pi_l\mathbb E\alpha_l,\delta x_i\rangle+\langle \tilde{Q}\widetilde x_i(t-\delta),\delta x_i(t-\delta)\rangle\\
&-\langle (\tilde{Q}\tilde{S}+ \tilde{S}^\top \tilde{Q} -\tilde{S}^\top \tilde{Q}\tilde{S})\sum_{l=1}^K\pi_l\mathbb E\alpha_l(t-\delta),\delta x_i(t-\delta)\rangle+\sum_{k=1}^K\langle \pi_k\tilde{A}^\top(t+\delta) \mathbb E\mathbf{y}_k(t+\delta)I_{[0,T-\delta]},\delta x_i\rangle\\
&+\sum_{k=1}^K\langle \pi_k\tilde{A}^\top(t+\delta) y_2^k(t+\delta)I_{[0,T-\delta]},\delta x_i\rangle+\sum_{k=1}^K\langle \pi_k\tilde{B}^\top(t+\theta) \mathbb E\mathbf{y}_k(t+\theta)I_{[0,T-\theta]},\delta u_i\rangle\\
&+\sum_{k=1}^K\langle \pi_k\tilde{B}^\top(t+\theta) y_2^k(t+\theta)I_{[0,T-\theta]},\delta u_i\rangle
+\langle R_{\theta_i}\widetilde u_i,\delta u_i\rangle+\langle \tilde{R}_{\theta_i}\widetilde u_i(t-\theta),\delta u_i(t-\theta)\rangle\Big]dt \\
&+\langle G\widetilde x_i(T),\delta x_i(T)\rangle-\langle (G\Gamma+\Gamma ^\top G-\Gamma ^\top G\Gamma)\sum_{l=1}^K\pi_l\mathbb E\alpha_l(T),\delta x_i(T)\rangle\Bigg\}+\sum_{l=1}^{12}\varepsilon_l,
\end{aligned}\end{equation}
where
\small\begin{equation}\nonumber\left\{\begin{aligned}
&\varepsilon_1=\mathbb E\int_0^T\big\langle (QS+ S^\top Q- S^\top Q S)(\sum_{l=1}^K\pi_l\mathbb E\alpha_l-\widetilde x^{(N)}),N\delta x^{(N)}\big\rangle dt,\\
&\varepsilon_2=\mathbb E\int_0^T\big\langle (\tilde{Q}\tilde{S}+ \tilde{S}^\top \tilde{Q} -\tilde{S}^\top \tilde{Q}\tilde{S})(\sum_{l=1}^K\pi_l\mathbb E\alpha_l(t-\delta)-\widetilde x^{(N)}(t-\delta)),N\delta x^{(N)}(t-\delta)\big\rangle dt,\\
&\varepsilon_3=\mathbb E\big\langle (G\Gamma+\Gamma ^\top G-\Gamma ^\top G\Gamma)(\sum_{l=1}^K\pi_l\mathbb E\alpha_l(T)-\widetilde x^{(N)}(T)),N\delta x^{(N)}(T)\big\rangle,\\
&\varepsilon_4=\sum_{k=1}^K\mathbb E\int_0^T\langle (QS+ S^\top Q-S^\top QS)\sum_{l=1}^K\pi_l\mathbb E\alpha_l,x_k^{**}-\delta x_{(k)}\rangle dt,\\
&\varepsilon_5=\sum_{k=1}^K\mathbb E\int_0^T\frac{1}{N_k}\sum_{j\in\mathcal I_k,j\neq i}\langle Q\widetilde x_j,N_k\delta x_j-x^*_j\rangle dt,\\
&\varepsilon_6=\sum_{k=1}^K\mathbb E\int_0^T\langle (\tilde{Q}\tilde{S}+ \tilde{S}^\top \tilde{Q} -\tilde{S}^\top \tilde{Q}\tilde{S})\sum_{l=1}^K\pi_l\mathbb E\alpha_l(t-\delta),x_k^{**}(t-\delta)-\delta x_{(k)}(t-\delta)\rangle dt,\\
&\varepsilon_7=\sum_{k=1}^K\mathbb E\int_0^T\frac{1}{N_k}\sum_{j\in\mathcal I_k,j\neq i}\langle \tilde{Q}\widetilde x_j(t-\delta),N_k\delta x_j(t-\delta)-x^*_j(t-\delta)\rangle dt,\\
&\varepsilon_8=\sum_{k=1}^K\mathbb E\langle (G\Gamma+\Gamma ^\top G-\Gamma ^\top G\Gamma)\sum_{l=1}^K\pi_l\mathbb E\alpha_l(T),x_k^{**}(T)-\delta x_{(k)}(T)\rangle,\\
&\varepsilon_9=\sum_{k=1}^K\frac{1}{N_k}\sum_{j\in\mathcal I_k,j\neq i}\langle G\widetilde x_j(T),N_k\delta x_j(T)-x^*_j(T)\rangle,\\
&\varepsilon_{10}=\sum_{k=1}^K\mathbb E\int_0^T\langle \frac{1}{N_k}\sum_{j\in\mathcal I_k,j\neq i}\pi_k\tilde{A}^\top(t+\delta) \mathbb E^{\mathcal F_t}[y_1^j(t+\delta)]I_{[0,T-\delta]}-\pi_k\tilde{A}^\top(t+\delta) \mathbb E\mathbf{y}_k(t+\delta)I_{[0,T-\delta]},\delta x_i\rangle dt,\\
&\varepsilon_{11}=\sum_{k=1}^K\mathbb E\int_0^T\Big\langle \sum_{l=1}^K\frac{\pi_l}{N_l}\sum_{j\in\mathcal I_l,j\neq i}\tilde{A}^\top(t+\delta) \mathbb E^{\mathcal F_t}[y_1^j(t+\delta)]I_{[0,T-\delta]}-\sum_{l=1}^K\pi_l\tilde{A}^\top(t+\delta)\mathbb E\mathbf{y}_l(t+\delta)I_{[0,T-\delta]},x_k^{**}\Big\rangle dt,\\
&\varepsilon_{12}=\sum_{k=1}^K\mathbb E\int_0^T\langle \frac{1}{N_k}\sum_{j\in\mathcal I_k,j\neq i}\pi_k\tilde{B}^\top(t+\theta) \mathbb E^{\mathcal F_t}[y_1^j(t+\theta)]I_{[0,T-\theta]}-\pi_k\tilde{B}^\top(t+\theta) \mathbb E\mathbf{y}_k(t+\theta)I_{[0,T-\theta]},\delta u_i\rangle dt.
\end{aligned}\right.\end{equation}\normalsize

Now, let us consider the case that agent $\mathcal A_i,\ 1\leq i\leq N$ uses an alternative strategy $u_i$ while other agents $\mathcal A_j,j\neq i$ use the strategy $\widetilde u_{-i}$. It is remarkable that hereafter the notation $\breve{x}_i,\ 1\leq i\leq N$ stands for the state of agent $\mathcal A_i$ when applying an alternative strategy $u_i$ while other agents $\mathcal A_l,l\neq j$ use the strategy $\widetilde u_{-i}$.

The realized state with the $i$-th agent's perturbation is
\begin{equation}\left\{\label{eq40}\begin{aligned}
d\breve x_i(t)=&\Big[A_{\theta_i}\breve x_i+\hat{A}_{\theta_i}\breve x_i(t-\delta)+\tilde A\breve x^{(N)}(t-\delta)+Bu_i+\hat{B}u_i(t-\theta)+\tilde B\breve u^{(N)}(t-\theta)\Big]dt\\
&+\Big[Du_i+\hat Du_i(t-\theta)\Big]dW_i(t),\ t\in[0,T],\\
\breve x_i(0)=&\xi_i,\ \breve x_i(t)=x^0(t),\ t\in[-\delta,0),\ u_i(t)=u^0(t),\ t\in[-\theta,0),
\end{aligned}\right.\end{equation}
and
\begin{equation}\left\{\label{eq41}\begin{aligned}
d\breve x_j(t)=&\Big[A_{\theta_j}\breve x_j+\hat{A}_{\theta_j}\breve x_j(t-\delta)+\tilde A\breve x^{(N)}(t-\delta)+B\widetilde u_j+\hat{B}\widetilde u_j(t-\theta)+\tilde B\breve u^{(N)}(t-\theta)\Big]dt\\
&+\Big[D\widetilde u_j+\hat D\widetilde u_j(t-\theta)\Big]dW_i(t),\ t\in[0,T],\\
\breve x_j(0)=&\xi_j,\ \breve x_j(t)=x^0(t),\ t\in[-\delta,0),\ \widetilde u_j(t)=u^0(t),\ t\in[-\theta,0),\ 1\leq j\leq N,\ j\neq i,
\end{aligned}\right.\end{equation}
where $\breve x^{(N)}=\frac{1}{N}\sum_{i=1}^N\breve x_i$, $\breve u^{(N)}=\frac{1}{N}(u_i+\sum_{j=1,j\neq i}^N\widetilde u_j)$.

In order to obtain the asymptotic $\varepsilon-$optimality, firstly we need some estimations. In the following arguments, $C$ will denote a constant whose value may change from line to line. Similar to the proof of \cite[Lemma 5.1]{HHN2018}, by virtue of estimations of AFBSDDE, we derive
\begin{lemma}\label{5.1}
Let \emph{(}A1\emph{)}-\emph{(}A4\emph{)} hold. Then there exists a constant $C\geq 0$ independent of $N$ such that
\begin{equation}\label{estimation-1}\begin{aligned}
&\sum_{l=1}^K\mathbb E\sup_{0\leq t\leq T}\Big[|\alpha_l(t)|^2+|\beta_l(t)|^2+|\check y_l(t)|^2+|\zeta_l(t)|^2\Big]\\
&+\sup_{1\leq i\leq N}\mathbb E\sup_{0\leq t\leq T}|\widetilde x_i(t)|^2+\sum_{l=1}^K\mathbb E\int_0^T\Big[|\gamma_l(t)|^2+|\check z_l(t)|^2\Big]dt\leq C.
\end{aligned}\end{equation}
\end{lemma}
 Similar to Lemma \ref{5.1}, by the $L^2$ boundness of $u_i,\xi_i$ and $\xi_j,p_j,q_j$ $(1\leq j\leq N,j\neq i)$ we have
 \begin{equation}\label{estimation-2}
\sup_{1\leq i\leq N}\mathbb E\sup_{0\leq t\leq T}\left|\breve x_i(t)\right|^2\leq C.
\end{equation}

\begin{lemma}\label{lemma5.2}
Let \emph{(}A1\emph{)}-\emph{(}A4\emph{)} hold. Then there exists a constant $C\geq 0$ independent of $N$ such that
\begin{equation}\label{estimation-3}
\mathbb E\sup_{0\leq t\leq T}\left|\widetilde x^{(N)}(t)-\sum_{l=1}^K\pi_l\mathbb E\alpha_l(t)\right|^2\leq \frac{C}{N}+C\epsilon_N^2,\quad
\mathbb E\sup_{0\leq t\leq T}\left|\widetilde u^{(N)}(t)-\sum_{l=1}^K\pi_l\mathbb Ev_l(t)\right|^2\leq \frac{C}{N}+C\epsilon_N^2,
\end{equation}
where $\epsilon_N=\sup_{1\leq l\leq K}\left|\pi_l^{(N)}-\pi_l\right|.$
\end{lemma}
Please refer to Appendix C for the proof. By Lemma \ref{lemma5.2}, we easily obtain the following result.
\begin{lemma}\label{lemma5.22}
Let \emph{(}A1\emph{)}-\emph{(}A4\emph{)} hold. Then there exists a constant $C\geq 0$ independent of $N$ such that
\begin{equation}\label{estimation-33}
\sup_{1\leq j\leq N}\mathbb E\sup_{0\leq t\leq T}\big|x_j(t)-\widetilde x_j(t)\big|^2\leq \frac{C}{N}+C\epsilon_N^2.
\end{equation}
\end{lemma}
\begin{lemma}\label{Lemma5}
Let \emph{(}A1\emph{)}-\emph{(}A4\emph{)} hold. Then there exists a constant $C\geq 0$ independent of $N$ such that
\begin{equation}\begin{aligned}\label{estimation-5}
&\sup_{1\leq j\leq N,j\neq i}\mathbb E\sup_{0\leq t\leq T}|\delta x_j(t)|^2
\leq\frac{C}{N^2}\mathbb E\int_0^T|\delta u_i|^2ds.
\end{aligned}\end{equation}
\end{lemma}
The proof is given in Appendix D.
\begin{remark}
Note that in \eqref{estimation-5}, the upper bound depends on $\mathbb E\int_0^T|\delta u_i|^2ds$. However, when studying the asymptotic optimality, we only need to consider the perturbations satisfying \eqref{condition for perturbation} below. Hence, in Section \ref{asymptotic optimality} when applying Lemma \ref{Lemma5}, similar estimation still holds while the upper bound is $\frac{C}{N^2}$ and $C$ is a general constant.
\end{remark}

\begin{lemma}\label{lemma5.4}
Let \emph{(}A1\emph{)}-\emph{(}A4\emph{)} hold. Then there exists a constant $C\geq 0$ independent of $N$ such that
\begin{equation}\label{estimation-6}
\sum_{l=1}^K\mathbb E\sup_{0\leq t\leq T}\left|x_l^{**}(t)-\delta x_{(l)}(t)\right|^2\leq \Big(\frac{C}{N^2}+C\epsilon_N^2\Big)  \mathbb E\int_0^T|\delta u_i|^2ds,
\end{equation}
and for $j\in\mathcal I_k,\ j\neq i$, $1\leq k\leq K$,
\begin{equation}\label{estimation-7}
\mathbb E\sup_{0\leq t\leq T}\left|N_k\delta x_j(t)-x^*_j(t)\right|^2\leq \Big(\frac{C}{N^2}+C\epsilon_N^2\Big)  \mathbb E\int_0^T|\delta u_i|^2ds.
\end{equation}
\end{lemma}
Applying the above estimations, by the standard estimations of AFBSDDE, \eqref{eq17} and \eqref{eq36}, we can obtain the following result.
\begin{lemma}\label{Lemma7}
Let \emph{(}A1\emph{)}-\emph{(}A4\emph{)} hold. Then there exists a constant $C\geq0$ independent of $N$ such that
\begin{equation}\label{estimation-9}
\sum_{k=1}^K\mathbb E\sup_{0\leq t\leq T-\delta}\left|\frac{1}{N_k}\sum_{j\in\mathcal I_k,j\neq i} \mathbb E^{\mathcal F_t}[y_1^j(t+\delta)]-\mathbb E\mathbf{y}_k(t+\delta)\right|^2\leq \frac{C}{N}+C\epsilon_N^2.
\end{equation}
\end{lemma}
Please refer to Appendix E and F for the proofs of Lemmas \ref{lemma5.4}-\ref{Lemma7}.
\subsection{Asymptotic optimality}

In order to prove asymptotic optimality, it suffices to consider the perturbations $u_{i}\in\mathcal U_i^c,\ 1\leq i\leq N$ such that $\mathcal J_{soc}^{(N)}(u_1,\cdots,u_N)\leq\mathcal J_{soc}^{(N)}(\widetilde u_1,$ $\cdots,\widetilde u_N).$ It is easy to check that $$\mathcal J_{soc}^{(N)}(\widetilde u_1,\cdots,\widetilde u_N)\leq CN,$$ where $C\geq0$ is a constant independent of $N$. Therefore, in the following we only consider the perturbations $u_i\in\mathcal U_i^{c}$ satisfying
\begin{equation}\label{condition for perturbation}
\sum_{i=1}^N\mathbb E\int_0^T|u_i|^2dt\leq CN.
\end{equation}
Let $\delta u_i=u_i-\widetilde u_i,\ 1\leq i\leq N$,
 and consider a perturbation  ${u} = \widetilde{u} + (\delta u_1,\cdots,\delta u_N):=\widetilde{ u}+ \delta u$.
% Therefore, recalling Lemma \ref{Lemma5} and Lemma \ref{lemma5.4}, there exists a constant $C\geq0$ independent of $N$ such that
%\begin{equation*}\begin{aligned}
%&\sup_{1\leq j\leq N,j\neq i}\mathbb E\sup_{0\leq t\leq T}|\delta x_j(t)|^2
%\leq\frac{C}{N^2},\quad \sum_{l=1}^K\mathbb E\sup_{0\leq t\leq T}\left|x_l^{**}(t)-\delta x_{(l)}(t)\right|^2\leq\frac{C}{N^2}+C\epsilon_N^2 ,\\
%&\mathbb E\sup_{0\leq t\leq T}\left|N_k\delta x_j(t)-x^*_j(t)\right|^2\leq \frac{C}{N^2}+C\epsilon_N^2.
%\end{aligned}\end{equation*}
%{\color{blue}\begin{remark}Without causing an ambiguity, here we still use the notation $\delta x_j(t), x_l^{**}(t), \delta x_{(l)}(t),x^*_j(t)$, we mention that different to the previous Sections where the perturbation takes only for for one components of $\widetilde{u}$, here the perturbation takes for all component of $\widetilde{u}$.\end{remark}}

\begin{theorem}\label{asymptotic optimal}
Let \emph{(}A1\emph{)}-\emph{(}A4\emph{)} hold. Then $\widetilde u=(\widetilde u_1,\cdots,\widetilde u_N)$ is a $\Big(\frac{1}{\sqrt{N}}+\epsilon_N\Big)$-optimal strategy for the agents.
\end{theorem}
\begin{proof}
We denote $\frac{\partial \mathcal{J}^{(N)}_{soc}(\widetilde{u} )}{\partial \widetilde{u}_i}(\delta u_i)$ as the variation of  $\mathcal{J}^{(N)}_{soc}(\widetilde{u} )$ on the $i$-th direction.  From \eqref{eq42}, we know that, $\frac{\partial \mathcal{J}^{(N)}_{soc}(\widetilde{u} )}{\partial \widetilde{u}_i}(\delta u_i)$  exists and is continuous, for each $i=1,2,\ldots, N$. Thus
\begin{equation*}
  \begin{aligned}
&\mathcal{J}^{(N)}_{soc}(\widetilde{u} + \delta u)-\mathcal{J}^{(N)}_{soc}(\widetilde{u} )
=\sum_{i=1}^N \frac{\partial \mathcal{J}^{(N)}_{soc}(\widetilde{u} )}{\partial \widetilde{u}_i}(\delta u_i)+o(\|\delta u\|).
\end{aligned}
\end{equation*}
Here $o(\|\delta u\|)$ is the higher order infinitesimal of $\|\delta u\|$ and the norm of $\delta u$ is given by
\[\|\delta u\|^2:=\sum_{i=1}^N\mathbb E\int_0^T|u_i|^2dt.
 \]
 From \eqref{condition for perturbation}, we know that $o(\|\delta u\|)=o(\sqrt{N})$ which yields that $\frac{1}{N}o(\|\delta u\|)=o(\frac{1}{\sqrt{N}})$.

 Therefore, in order to  prove asymptotic optimality, we only need to show that $$\frac{\partial \mathcal{J}^{(N)}_{soc}(\widetilde{u} )}{\partial \widetilde{u}_i}(\delta u_i)=O\Big(\frac{1}{\sqrt{N}}+\epsilon_N\Big).$$
According to Section \ref{Agent's perturbation}, we derive
\begin{equation}\begin{aligned}
&\frac{\partial \mathcal{J}^{(N)}_{soc}(\widetilde{u} )}{\partial \widetilde{u}_i}(\delta u_i)=\mathbb E\Bigg\{\int_0^T\Big[\langle Q\widetilde x_i,\delta x_i\rangle-\langle (QS+ S^\top Q -S^\top QS)\sum_{l=1}^K\pi_l\mathbb E\alpha_l,\delta x_i\rangle+\langle \tilde{Q}\widetilde x_i(t-\delta),\delta x_i(t-\delta)\rangle\\
&-\langle (\tilde{Q}\tilde{S}+ \tilde{S}^\top \tilde{Q} -\tilde{S}^\top \tilde{Q}\tilde{S})\sum_{l=1}^K\pi_l\mathbb E\alpha_l(t-\delta),\delta x_i(t-\delta)\rangle+\sum_{k=1}^K\langle \pi_k\tilde{A}^\top(t+\delta) \mathbb E\mathbf{y}_k(t+\delta)I_{[0,T-\delta]},\delta x_i\rangle\\
&+\sum_{k=1}^K\langle \pi_k\tilde{A}^\top(t+\delta) y_2^k(t+\delta)I_{[0,T-\delta]},\delta x_i\rangle+\sum_{k=1}^K\langle \pi_k\tilde{B}^\top(t+\theta) \mathbb E\mathbf{y}_k(t+\theta)I_{[0,T-\theta]},\delta u_i\rangle\\
&+\sum_{k=1}^K\langle \pi_k\tilde{B}^\top(t+\theta) y_2^k(t+\theta)I_{[0,T-\theta]},\delta u_i\rangle
+\langle R_{\theta_i}\widetilde u_i,\delta u_i\rangle+\langle \tilde{R}_{\theta_i}\widetilde u_i(t-\theta),\delta u_i(t-\theta)\rangle\Big]dt \\
&+\langle G\widetilde x_i(T),\delta x_i(T)\rangle-\langle (G\Gamma+\Gamma ^\top G-\Gamma ^\top G\Gamma)\sum_{l=1}^K\pi_l\mathbb E\alpha_l(T),\delta x_i(T)\rangle\Bigg\}+\sum_{l=1}^{12}\varepsilon_l.
\end{aligned}\end{equation}
It follows from the optimality of $\widetilde u$ that
\begin{equation*}\begin{aligned}
&\mathbb E\Bigg\{\int_0^T\Big[\langle Q\widetilde x_i,\delta x_i\rangle-\langle (QS+ S^\top Q -S^\top QS)\sum_{l=1}^K\pi_l\mathbb E\alpha_l,\delta x_i\rangle+\langle \tilde{Q}\widetilde x_i(t-\delta),\delta x_i(t-\delta)\rangle\\
&-\langle (\tilde{Q}\tilde{S}+ \tilde{S}^\top \tilde{Q} -\tilde{S}^\top \tilde{Q}\tilde{S})\sum_{l=1}^K\pi_l\mathbb E\alpha_l(t-\delta),\delta x_i(t-\delta)\rangle+\sum_{k=1}^K\langle \pi_k\tilde{A}^\top(t+\delta) \mathbb E\mathbf{y}_k(t+\delta)I_{[0,T-\delta]},\delta x_i\rangle\\
&+\sum_{k=1}^K\langle \pi_k\tilde{A}^\top(t+\delta) y_2^k(t+\delta)I_{[0,T-\delta]},\delta x_i\rangle+\sum_{k=1}^K\langle \pi_k\tilde{B}^\top(t+\theta) \mathbb E\mathbf{y}_k(t+\theta)I_{[0,T-\theta]},\delta u_i\rangle\\
&+\sum_{k=1}^K\langle \pi_k\tilde{B}^\top(t+\theta) y_2^k(t+\theta)I_{[0,T-\theta]},\delta u_i\rangle
+\langle R_{\theta_i}\widetilde u_i,\delta u_i\rangle+\langle \tilde{R}_{\theta_i}\widetilde u_i(t-\theta),\delta u_i(t-\theta)\rangle\Big]dt \\
&+\langle G\widetilde x_i(T),\delta x_i(T)\rangle-\langle (G\Gamma+\Gamma ^\top G-\Gamma ^\top G\Gamma)\sum_{l=1}^K\pi_l\mathbb E\alpha_l(T),\delta x_i(T)\rangle\Bigg\}=0.
\end{aligned}\end{equation*}
Moreover, by Lemmas \ref{lemma5.2}-\ref{Lemma7} (recall that when calculating $\frac{\partial \mathcal{J}^{(N)}_{soc}(\widetilde{u} )}{\partial \widetilde{u}_i}(\delta u_i)$, it only needs to take  perturbation for one component of $\widetilde{u}$), we have $$\sum_{l=1}^{12}\varepsilon_l=
O\Big(\frac{1}{\sqrt{N}}+\epsilon_N\Big).$$
Therefore, $$\frac{\partial \mathcal{J}^{(N)}_{soc}(\widetilde{u} )}{\partial \widetilde{u}_i}(\delta u_i)=O\Big(\frac{1}{\sqrt{N}}+\epsilon_N\Big).$$
\end{proof}

\section{Conclusion}

In this paper, we mainly study a class of stochastic LQ dynamic optimization problems involving a large number of weakly-coupled heterogeneous agents. These agents formalize a team with cooperation to minimize a social cost functional, while the dynamic is driven by SDDE. We formulate an auxiliary control problem and derive the decentralized social strategy based on the CC system. Applying the discounting decoupling method, we establish the wellpoesdness of this MF-AFBSDDE. Finally, we verify the related asymptotic social optimality. In the future, one possible research direction is to study LQ MF social optima with delay in FBSDE scheme, which may involve more applications in practice and bring more challenges in theory. Another research problem is corresponding case in partial information or observation scheme (e.g., behaviors in incomplete market), which may be more valuable and challenging.

\section{Appendix}

\textbf{Appendix A. Proof of Lemma \ref{l1}.}\\
\\
\begin{proof}
Denote \small $\hat b(s):=b(s,X_1(s),X_1(s-\delta),\mathbb EX_1(s-\delta),Y_1(s),Y_1(s-\theta)I_{[\theta,T]},\mathbb EY_1(s),\mathbb EY_1(s-\theta)I_{[\theta,T]},$
$\mathbb E^{\mathcal F_s}[Y_1(s+\theta)]I_{[0,T-\theta]},\mathbb EY_1(s+\theta)I_{[0,T-\theta]},\mathbb E^{\mathcal F_{s-\theta}}[Y_1(s)]I_{[\theta,T]},Z_1(s),Z_1(s-\theta)I_{[\theta,T]},\mathbb EZ_1(s),\mathbb EZ_1(s-\theta)I_{[\theta,T]},$
$\mathbb E^{\mathcal F_s}[Z_1(s+\theta)]I_{[0,T-\theta]},\mathbb E^{\mathcal F_{s-\theta}}[Z_1(s)]I_{[\theta,T]})-b(s,X_2(s),X_2(s-\delta),\mathbb EX_2(s-\delta),Y_2(s),Y_2(s-\theta)I_{[\theta,T]},\mathbb EY_2(s),$
$\mathbb EY_2(s-\theta)I_{[\theta,T]},\mathbb E^{\mathcal F_s}[Y_2(s+\theta)]I_{[0,T-\theta]},\mathbb EY_2(s+\theta)I_{[0,T-\theta]},\mathbb E^{\mathcal F_{s-\theta}}[Y_2(s)]I_{[\theta,T]},Z_2(s),Z_2(s-\theta)I_{[\theta,T]},\mathbb EZ_2(s),$
$\mathbb EZ_2(s-\theta)I_{[\theta,T]},\mathbb E^{\mathcal F_s}[Z_2(s+\theta)]I_{[0,T-\theta]},\mathbb E^{\mathcal F_{s-\theta}}[Z_2(s)]I_{[\theta,T]})$ and $\hat\sigma(s):=\sigma(s,Y_1(s),Y_1(s-\theta)I_{[\theta,T]},\mathbb E^{\mathcal F_s}[Y_1(s+\theta)]I_{[0,T-\theta]},\mathbb E^{\mathcal F_{s-\theta}}[Y_1(s)]I_{[\theta,T]},\mathbb EY_1(s+\theta)I_{[0,T-\theta]},\mathbb EY_1(s)I_{[\theta,T]},Z_1(s),Z_1(s-\theta)I_{[\theta,T]},\mathbb E^{\mathcal F_s}[Z_1(s+\theta)]I_{[0,T-\theta]},$
$\mathbb E^{\mathcal F_{s-\theta}}[Z_1(s)]I_{[\theta,T]})-\sigma(s,Y_2(s),Y_2(s-\theta)I_{[\theta,T]},\mathbb E^{\mathcal F_s}[Y_2(s+\theta)]I_{[0,T-\theta]},\mathbb E^{\mathcal F_{s-\theta}}[Y_2(s)]I_{[\theta,T]},\mathbb EY_2(s+\theta)I_{[0,T-\theta]},$
$\mathbb EY_2(s)I_{[\theta,T]},Z_2(s),Z_2(s-\theta)I_{[\theta,T]},\mathbb E^{\mathcal F_s}[Z_2(s+\theta)]I_{[0,T-\theta]},\mathbb E^{\mathcal F_{s-\theta}}[Z_2(s)]I_{[\theta,T]})$. \normalsize Applying It\^o's formula to $e^{-\rho s}|\widehat X(s)|^2$, we obtain
\begin{equation}\label{eq31}\begin{aligned}
e^{-\rho t}\mathbb E|\widehat{X}(t)|^2=-\rho\mathbb E\int_0^te^{-\rho s}|\widehat{X}(s)|^2ds+2\mathbb E\int_0^te^{-\rho s}\langle\widehat X(s),\hat b(s) \rangle ds+\mathbb E\int_0^te^{-\rho s}|\hat\sigma(s)|^2ds.
\end{aligned}\end{equation}
It follows that\small
\begin{equation}\nonumber\begin{aligned}
&2\langle\widehat X(s),\hat b(s) \rangle\\
=&2\langle\widehat X(s),b(s,X_1(s),X_1(s-\delta),\mathbb EX_1(s-\delta),Y_1(s),Y_1(s-\theta)I_{[\theta,T]},\mathbb EY_1(s),\mathbb EY_1(s-\theta)I_{[\theta,T]},\mathbb E^{\mathcal F_s}[Y_1(s+\theta)]I_{[0,T-\theta]},\\
&\mathbb EY_1(s+\theta)I_{[0,T-\theta]},\mathbb E^{\mathcal F_{s-\theta}}[Y_1(s)]I_{[\theta,T]},Z_1(s),Z_1(s-\theta)I_{[\theta,T]},\mathbb EZ_1(s),\mathbb EZ_1(s-\theta)I_{[\theta,T]},\mathbb E^{\mathcal F_s}[Z_1(s+\theta)]I_{[0,T-\theta]},\\
&\mathbb E^{\mathcal F_{s-\theta}}[Z_1(s)]I_{[\theta,T]})-b(s,X_2(s),X_1(s-\delta),\mathbb EX_1(s-\delta),Y_1(s),Y_1(s-\theta)I_{[\theta,T]},\mathbb EY_1(s),\mathbb EY_1(s-\theta)I_{[\theta,T]},\\
&\mathbb E^{\mathcal F_s}[Y_1(s+\theta)]I_{[0,T-\theta]},\mathbb EY_1(s+\theta)I_{[0,T-\theta]},\mathbb E^{\mathcal F_{s-\theta}}[Y_1(s)]I_{[\theta,T]},Z_1(s),Z_1(s-\theta)I_{[\theta,T]},\mathbb EZ_1(s),\mathbb EZ_1(s-\theta)I_{[\theta,T]},\\
&\mathbb E^{\mathcal F_s}[Z_1(s+\theta)]I_{[0,T-\theta]},\mathbb E^{\mathcal F_{s-\theta}}[Z_1(s)]I_{[\theta,T]})\rangle+2\langle\widehat X(s),b(s,X_2(s),X_1(s-\delta),\mathbb EX_1(s-\delta),Y_1(s),Y_1(s-\theta)I_{[\theta,T]},\\
&\mathbb EY_1(s),\mathbb EY_1(s-\theta)I_{[\theta,T]},\mathbb E^{\mathcal F_s}[Y_1(s+\theta)]I_{[0,T-\theta]},\mathbb EY_1(s+\theta)I_{[0,T-\theta]},\mathbb E^{\mathcal F_{s-\theta}}[Y_1(s)]I_{[\theta,T]},Z_1(s),Z_1(s-\theta)I_{[\theta,T]},\\
&\mathbb EZ_1(s),\mathbb EZ_1(s-\theta)I_{[\theta,T]},\mathbb E^{\mathcal F_s}[Z_1(s+\theta)]I_{[0,T-\theta]},\mathbb E^{\mathcal F_{s-\theta}}[Z_1(s)]I_{[\theta,T]})-b(s,X_2(s),X_2(s-\delta),\mathbb EX_2(s-\delta),Y_2(s),\\
&Y_2(s-\theta)I_{[\theta,T]},\mathbb EY_2(s),\mathbb EY_2(s-\theta)I_{[\theta,T]},\mathbb E^{\mathcal F_s}[Y_2(s+\theta)]I_{[0,T-\theta]},\mathbb EY_2(s+\theta)I_{[0,T-\theta]},\mathbb E^{\mathcal F_{s-\theta}}[Y_2(s)]I_{[\theta,T]},Z_2(s),\\
&Z_2(s-\theta)I_{[\theta,T]},\mathbb EZ_2(s),\mathbb EZ_2(s-\theta)I_{[\theta,T]},\mathbb E^{\mathcal F_s}[Z_2(s+\theta)]I_{[0,T-\theta]},\mathbb E^{\mathcal F_{s-\theta}}[Z_2(s)]I_{[\theta,T]})\rangle\\
\leq&(2\rho_1+k_1+k_2+\sum_{j=1}^{13} k_{j+2}l_j^{-1})|\widehat X(s)|^2+k_1|\widehat X(s-\delta)|^2+k_2\mathbb E|\widehat X(s-\delta)|^2+k_3l_1|\widehat Y(s)|^2+k_4l_2|\widehat Y(s-\theta)|^2I_{[\theta,T]}\\
&+k_5l_3\mathbb E|\widehat Y(s)|^2+k_6l_4\mathbb E|\widehat Y(s-\theta)|^2I_{[\theta,T]}+k_7l_5\mathbb E^{\mathcal F_s}[|\widehat Y(s+\theta)|^2]I_{[0,T-\theta]}+k_8l_6\mathbb E|\widehat Y(s+\theta)|^2I_{[0,T-\theta]}\\
&+k_9l_7\mathbb E^{\mathcal F_{s-\theta}}[|\widehat Y(s)|^2]I_{[\theta,T]}+k_{10}l_8|\widehat Z(s)|^2+k_{11}l_9|\widehat Z(s-\theta)|^2I_{[\theta,T]}+k_{12}l_{10}\mathbb E|\widehat Z(s)|^2+k_{13}l_{11}\mathbb E|\widehat Z(s-\theta)|^2I_{[\theta,T]}\\
&+k_{14}l_{12}\mathbb E^{\mathcal F_s}[|\widehat Z(s+\theta)|^2]I_{[0,T-\theta]}+k_{15}l_{13}\mathbb E^{\mathcal F_{s-\theta}}[|\widehat Z(s)|^2]I_{[\theta,T]},
\end{aligned}\end{equation} \normalsize
and
\begin{equation}\nonumber\begin{aligned}
|\hat\sigma(s)|^2\leq &k_{20}^2|\widehat Y(s)|^2+k_{21}^2|\widehat Y(s-\theta)|^2I_{[\theta,T]}+k_{22}^2|\mathbb E^{\mathcal F_s}[\widehat Y(s+\theta)]|^2I_{[0,T-\theta]}+k_{23}^2|\mathbb E^{\mathcal F_{s-\theta}}[\widehat Y(s)]|^2I_{[\theta,T]}\\
&+k_{24}^2|\mathbb E\widehat Y(s+\theta)|^2I_{[0,T-\theta]}+k_{25}^2|\mathbb E\widehat Y(s)|^2I_{[\theta,T]}+k_{26}^2|\widehat Z(s)|^2+k_{27}^2|\widehat Z(s-\theta)|^2I_{[\theta,T]}\\
&+k_{28}^2|\mathbb E^{\mathcal F_s}[\widehat Z(s+\theta)]|^2I_{[0,T-\theta]}+k_{29}^2|\mathbb E^{\mathcal F_{s-\theta}}[\widehat Z(s)]|^2I_{[\theta,T]}\\
\leq &k_{20}^2|\widehat Y(s)|^2+k_{21}^2|\widehat Y(s-\theta)|^2I_{[\theta,T]}+k_{22}^2\mathbb E^{\mathcal F_s}[|\widehat Y(s+\theta)|^2]I_{[0,T-\theta]}+k_{23}^2\mathbb E^{\mathcal F_{s-\theta}}[|\widehat Y(s)|^2]I_{[\theta,T]}\\
&+k_{24}^2\mathbb E|\widehat Y(s+\theta)|^2I_{[0,T-\theta]}+k_{25}^2\mathbb E|\widehat Y(s)|^2I_{[\theta,T]}+k_{26}^2|\widehat Z(s)|^2+k_{27}^2|\widehat Z(s-\theta)|^2I_{[\theta,T]}\\
&+k_{28}^2\mathbb E^{\mathcal F_s}[|\widehat Z(s+\theta)|^2]I_{[0,T-\theta]}+k_{29}^2\mathbb E^{\mathcal F_{s-\theta}}[|\widehat Z(s)|^2]I_{[\theta,T]}.
\end{aligned}\end{equation}
Note that
\begin{equation}\nonumber\begin{aligned}
\mathbb E\int_0^te^{-\rho s}\mathbb E|\widehat X(s-\delta)|^2ds&=\mathbb E\int_0^te^{-\rho s}|\widehat X(s-\delta)|^2ds=\mathbb E\int_0^{t-\delta}e^{-\rho (s+\delta)}|\widehat X(s)|^2ds\\
&\leq e^{-\rho\delta}\mathbb E\int_0^{t}e^{-\rho s}|\widehat X(s)|^2ds.
\end{aligned}\end{equation}
Similarly,
\begin{equation}\nonumber\begin{aligned}
&\mathbb E\int_0^te^{-\rho s}\mathbb E|\widehat Y(s-\theta)|^2I_{[\theta,T]}ds=\mathbb E\int_0^te^{-\rho s}|\widehat Y(s-\theta)|^2I_{[\theta,T]}ds\leq  e^{-\rho\theta}\mathbb E\int_0^{t}e^{-\rho s}|\widehat Y(s)|^2ds,\\
&\mathbb E\int_0^te^{-\rho s}\mathbb E|\widehat Z(s-\theta)|^2I_{[\theta,T]}ds=\mathbb E\int_0^te^{-\rho s}|\widehat Z(s-\theta)|^2I_{[\theta,T]}ds\leq  e^{-\rho\theta}\mathbb E\int_0^{t}e^{-\rho s}|\widehat Z(s)|^2ds,
\end{aligned}\end{equation}
and
\begin{equation}\nonumber\begin{aligned}
&\mathbb E\int_0^te^{-\rho s}\mathbb E^{\mathcal F_s}[|\widehat Y(s+\theta)|^2]I_{[0,T-\theta]}ds=\mathbb E\int_0^te^{-\rho s}|\widehat Y(s+\theta)|^2I_{[0,T-\theta]}ds
=\mathbb E\int_{\theta}^{t+\theta}e^{-\rho (s-\theta)}|\widehat Y(s)|^2I_{[\theta,T]}ds\\
&\qquad\leq e^{\rho\theta}\mathbb E\int_0^{t+\theta}e^{-\rho s}|\widehat Y(s)|^2I_{[0,T-\theta]}ds,\\
&\mathbb E\int_0^te^{-\rho s}\mathbb E^{\mathcal F_{s-\theta}}[|\widehat Y(s)|^2]ds=\mathbb E\int_0^{t}e^{-\rho s}|\widehat Y(s)|^2ds,
\quad\mathbb E\int_0^te^{-\rho s}\mathbb E^{\mathcal F_{s-\theta}}[|\widehat Z(s)|^2]ds= \mathbb E\int_0^{t}e^{-\rho s}|\widehat Z(s)|^2ds,\\
&\mathbb E\int_0^te^{-\rho s}\mathbb E^{\mathcal F_s}[|\widehat Z(s+\theta)|^2]I_{[0,T-\theta]}ds\leq e^{\rho\theta}\mathbb E\int_0^{t+\theta}e^{-\rho s}|\widehat Z(s)|^2I_{[\theta,T]}ds\leq  e^{\rho\theta}\mathbb E\int_0^{t+\theta}e^{-\rho s}|\widehat Z(s)|^2I_{[0,T-\theta]}ds.
\end{aligned}\end{equation}
By noticing\eqref{eq31}, we obtain \eqref{eq26}.

Similarly,  for $\tilde\rho_1$ which solves the equation
$\tilde\rho_1=\rho -2\rho_1-(k_1+k_2)(1+e^{-(\rho-\tilde\rho_1)\delta})-\sum_{j=1}^{13} k_{j+2}l_j^{-1}$, we
 apply It\^o's formula to $e^{-\tilde\rho_1(t-s)-\rho s}|\widehat X(s)|^2$ for $s\in[0,t]$, we have
\begin{equation}\label{eq32}\begin{aligned}
e^{-\rho t}\mathbb E|\widehat{X}(t)|^2=&-(\rho-\tilde\rho_1)\mathbb E\int_0^te^{-\tilde\rho_1(t-s)-\rho s}|\widehat{X}(s)|^2ds+2\mathbb E\int_0^te^{-\tilde\rho_1(t-s)-\rho s}\langle\widehat X(s),\hat b(s) \rangle ds\\
&+\mathbb E\int_0^te^{-\tilde\rho_1(t-s)-\rho s}|\hat\sigma(s)|^2ds.
\end{aligned}\end{equation}
From the above estimates and \eqref{eq32}, one can prove \eqref{eq27}. We mention that  there indeed exists a $\tilde\rho_1$ solving the equation
$
\tilde\rho_1=\rho -2\rho_1-(k_1+k_2)(1+e^{-(\rho-\tilde\rho_1)\delta})-\sum_{j=1}^{13} k_{j+2}l_j^{-1}.
$
In fact, from above equation, it follows that
\begin{equation}\label{discounting coefficient equation}
\tilde\rho_1-\rho+(k_1+k_2)e^{(\tilde\rho_1-\rho)\delta}=-2\rho_1-(k_1+k_2)-\sum_{j=1}^{13} k_{j+2}l_j^{-1}.
\end{equation}
By noticing that $k_1$, $k_2$, $\delta$ are positive which yields that $g_1(x):=x+(k_1+k_2)e^{x\delta}$ is an increasing continuous function satisfying $g_1(+\infty)=+\infty$ and $g_1(-\infty)=-\infty$. Therefore, for given $\rho$ and $-2\rho_1-(k_1+k_2)-\sum_{j=1}^{13} k_{j+2}l_j^{-1}$, there exists a unique  $\tilde\rho_1$ such that \eqref{discounting coefficient equation} holds by intermediate value theorem.

Integrating from 0 to $T$ on both sides of \eqref{eq27},
we have
\small
\begin{equation}\nonumber\begin{aligned}
&\|\widehat X\|_{\rho}^2\\
\leq& \Big(k_3l_1+k_5l_3+k_9l_7+k_{20}^2+k_{23}^2+k_{25}^2+(k_4l_2+k_6l_4+k_{21}^2) e^{-(\rho-\tilde\rho_1)\theta}\Big)\mathbb E\int_0^T\int_0^te^{-\tilde\rho_1(t-s) -\rho s}|\widehat{Y}(s)|^2dsdt\\
&+(k_7l_5+k_8l_6+k_{22}^2+k_{24}^2)e^{(\rho-\tilde\rho_1)\theta}\mathbb E\int_0^T\int_0^{t+\theta}e^{-\tilde\rho_1(t-s) -\rho s}|\widehat{Y}(s)|^2I_{[0,T-\theta]}(s)dsdt\\
&+\Big(k_{10}l_8+k_{12}l_{10}+k_{15}l_{13}+k_{26}^2+k_{29}^2+(k_{11}l_9+k_{13}l_{11}+k_{27}^2) e^{-(\rho-\tilde\rho_1)\theta}\Big)\mathbb E\int_0^T\int_0^te^{-\tilde\rho_1(t-s) -\rho s}|\widehat{Z}(s)|^2dsdt\\
&+(k_{14}l_{12}+k_{28}^2)e^{(\rho-\tilde\rho_1)\theta}\mathbb E\int_0^T\int_0^{t+\theta}e^{-\tilde\rho_1(t-s) -\rho s}|\widehat{Z}(s)|^2I_{[0,T-\theta]}(s)dsdt\\
=&\Big(k_3l_1+k_5l_3+k_9l_7+k_{20}^2+k_{23}^2+k_{25}^2+(k_4l_2+k_6l_4+k_{21}^2) e^{-(\rho-\tilde\rho_1)\theta}\Big)\mathbb E\int_0^T\frac{1-e^{-\tilde\rho_1(T-s)}}{\tilde\rho_1}e^{-\rho s}|\widehat{Y}(s)|^2ds\\
&+(k_7l_5+k_8l_6+k_{22}^2+k_{24}^2)e^{(\rho-\tilde\rho_1)\theta}\Big(\mathbb E\int_0^{\theta}\frac{e^{\tilde\rho_1s}-e^{-\tilde\rho_1(T-s)}}{\tilde\rho_1}e^{-\rho s}|\widehat{Y}(s)|^2ds+\mathbb E\int_{\theta}^{T}\frac{e^{\tilde\rho_1\theta}-e^{-\tilde\rho_1(T-s)}}{\tilde\rho_1}e^{-\rho s}|\widehat{Y}(s)|^2ds\Big)\\
&+\Big(k_{10}l_8+k_{12}l_{10}+k_{15}l_{13}+k_{26}^2+k_{29}^2+(k_{11}l_9+k_{13}l_{11}+k_{27}^2) e^{-(\rho-\tilde\rho_1)\theta}\Big)\mathbb E\int_0^T\frac{1-e^{-\tilde\rho_1(T-s)}}{\tilde\rho_1}e^{-\rho s}|\widehat{Z}(s)|^2ds\\
&+(k_{14}l_{12}+k_{28}^2)e^{(\rho-\tilde\rho_1)\theta}\Big(\mathbb E\int_0^{\theta}\frac{e^{\tilde\rho_1s}-e^{-\tilde\rho_1(T-s)}}{\tilde\rho_1}e^{-\rho s}|\widehat{Z}(s)|^2ds+\mathbb E\int_{\theta}^{T}\frac{e^{\tilde\rho_1\theta}-e^{-\tilde\rho_1(T-s)}}{\tilde\rho_1}e^{-\rho s}|\widehat{Z}(s)|^2ds\Big).
\end{aligned}\end{equation}
\normalsize
Noticing that for all $s\in[0,T]$,
$$\frac{1-e^{-\tilde\rho_1(T-s)}}{\tilde\rho_1}\leq \frac{1-e^{-\tilde\rho_1T}}{\tilde\rho_1}\leq \frac{e^{\tilde\rho_1\theta}-e^{-\tilde\rho_1T}}{\tilde\rho_1},\quad \frac{e^{\tilde\rho_1\theta}-e^{-\tilde\rho_1(T-s)}}{\tilde\rho_1}\leq \frac{e^{\tilde\rho_1\theta}-e^{-\tilde\rho_1T}}{\tilde\rho_1},$$
and for all $s\in[0,\theta]$,
$$\frac{e^{\tilde\rho_1s}-e^{-\tilde\rho_1(T-s)}}{\tilde\rho_1}\leq \frac{e^{\tilde\rho_1\theta}-e^{-\tilde\rho_1(T-s)}}{\tilde\rho_1},$$
then \eqref{eq28} follows. Letting $t=T$ in \eqref{eq27} and noticing $e^{-\tilde\rho_1(T-s)}\leq 1\vee e^{-\tilde\rho_1T}$ for all $s\in[0,T]$, we have \eqref{eq29}. \eqref{eq30} is obvious by \eqref{eq29}.
\end{proof}
\vspace{1cm}
\textbf{Appendix B. Proof of Lemma \ref{l2}.}\\
\\
\begin{proof}
Denote $\hat f(s):=f(s,X_1(t),\mathbb EX_1(t),Y_1(t),\mathbb EY_1(t+\delta)I_{[0,T-\delta]},\mathbb E^{\mathcal F_t}[Y_1(t+\delta)]I_{[0,T-\delta]})-f(s,X_2(t),\mathbb EX_2(t),Y_2(t),\mathbb EY_2(t+\delta)I_{[0,T-\delta]},\mathbb E^{\mathcal F_t}[Y_2(t+\delta)]I_{[0,T-\delta]})$. Applying It\^o's formula to $e^{-\rho t}|\widehat Y(t)|^2$, we obtain
\begin{equation}\label{eq39}\begin{aligned}
&e^{-\rho t}\mathbb E|\widehat{Y}(t)|^2-\rho\mathbb E\int_t^Te^{-\rho s}|\widehat{Y}(s)|^2ds+\mathbb E\int_t^Te^{-\rho s}|\widehat{Z}(s)|^2ds\\
=&e^{-\rho T}\mathbb E|\widehat{Y}(T)|^2+2\mathbb E\int_t^Te^{-\rho s}\langle\widehat Y(s),\hat f(s) \rangle ds.
\end{aligned}\end{equation}
It follows that $$|\widehat{Y}(T)|^2=|\Phi(X_1(T),\mathbb EX_1(T))-\Phi(X_2(T),\mathbb EX_2(T))|^2\leq k_{30}^2|\widehat{X}(T)|^2+k_{31}^2|\mathbb E\widehat{X}(T)|^2,$$ and
\begin{equation}\nonumber\begin{aligned}
&2\langle\widehat Y(s),\hat f(s) \rangle\\
=&2\langle\widehat Y(s),f(s,X_1(t),\mathbb EX_1(t),Y_1(t),\mathbb EY_1(t+\delta)I_{[0,T-\delta]},\mathbb E^{\mathcal F_t}[Y_1(t+\delta)]I_{[0,T-\delta]})\\
&-f(s,X_1(t),\mathbb EX_1(t),Y_2(t),\mathbb EY_1(t+\delta)I_{[0,T-\delta]},\mathbb E^{\mathcal F_t}[Y_1(t+\delta)]I_{[0,T-\delta]})\rangle\\
&+2\langle\widehat Y(s),f(s,X_1(t),\mathbb EX_1(t),Y_2(t),\mathbb EY_1(t+\delta)I_{[0,T-\delta]},\mathbb E^{\mathcal F_t}[Y_1(t+\delta)]I_{[0,T-\delta]})\\
&-f(s,X_2(t),\mathbb EX_2(t),Y_2(t),\mathbb EY_2(t+\delta)I_{[0,T-\delta]},\mathbb E^{\mathcal F_t}[Y_2(t+\delta)]I_{[0,T-\delta]})\rangle\\
\leq &2\rho_2|\widehat Y(s)|^2+2|\widehat Y(s)|\big(k_{16}|\widehat X(s)|+k_{17}|\mathbb E\widehat X(s)|+k_{18}|\mathbb E\widehat Y(s+\delta)|I_{[0,T-\delta]}+k_{19}|\mathbb E^{\mathcal F_s}[\widehat Y(s+\delta)]I_{[0,T-\delta]}|\big)\\
\leq&(2\rho_2+k_{16}l_{14}^{-1}+k_{17}l_{15}^{-1}+k_{18}+k_{19})|\widehat Y(s)|^2\\
&+k_{16}l_{14}|\widehat X(s)|^2+k_{17}l_{15}\mathbb E|\widehat X(s)|^2+k_{18}\mathbb E|\widehat Y(s+\delta)|^2I_{[0,T-\delta]}+k_{19}\mathbb E^{\mathcal F_s}[|\widehat Y(s+\delta)|^2]I_{[0,T-\delta]}.
\end{aligned}\end{equation}
Note that
\begin{equation}\nonumber\begin{aligned}
\mathbb E\int_t^Te^{-\rho s}\mathbb E|\widehat Y(s+\delta)|^2I_{[0,T-\delta]}ds=\mathbb E\int_{t+\delta}^{T}e^{-\rho (s-\delta)}|\widehat Y(s)|^2ds\leq e^{\rho\delta}\mathbb E\int_t^{T}e^{-\rho s}|\widehat Y(s)|^2ds.
\end{aligned}\end{equation}
Similarly,
\begin{equation}\nonumber\begin{aligned}
\mathbb E\int_t^Te^{-\rho s}\mathbb E^{\mathcal F_s}[|\widehat Y(s+\delta)|^2]I_{[0,T-\delta]}ds
&\leq e^{\rho\delta}\mathbb E\int_t^{T}e^{-\rho s}|\widehat Y(s)|^2ds.
\end{aligned}\end{equation}
By noticing\eqref{eq39}, we obtain \eqref{eq34}.

Similarly, applying It\^o's formula to $e^{-\tilde\rho_2(s-t)-\rho s}|\widehat Y(s)|^2$ for $s\in[0,t]$, we have
\begin{equation}\label{eq400}\begin{aligned}
&e^{-\rho t}\mathbb E|\widehat{Y}(t)|^2-(\rho+\tilde\rho_2)\mathbb E\int_t^Te^{-\tilde\rho_2(s-t)-\rho s}|\widehat{Y}(s)|^2ds+\mathbb E\int_t^Te^{-\tilde\rho_2(s-t)-\rho s}|\widehat{Z}(s)|^2ds\\
=&e^{-\tilde\rho_2(T-t)-\rho T}\mathbb E|\widehat{Y}(T)|^2+2\mathbb E\int_t^Te^{-\tilde\rho_2(s-t)-\rho s}\langle\widehat Y(s),\hat f(s) \rangle ds.
\end{aligned}\end{equation}
From the above estimates and \eqref{eq400}, one can prove \eqref{eq35}.  We mention that  there indeed exists a $\tilde\rho_2$ solving the equation
$\tilde\rho_2=-\rho -2\rho_2-k_{16}l_{14}^{-1}-k_{17}l_{15}^{-1}-(k_{18}+k_{19})(1+e^{(\rho+\tilde\rho_2)\delta})$
In fact, from above equation, it follows that
\begin{equation}\label{discounting coefficient equation 2}
\tilde\rho_2+\rho+(k_{18}+k_{19})e^{(\tilde\rho_2+\rho)\delta}= -2\rho_2-k_{16}l_{14}^{-1}-k_{17}l_{15}^{-1}-(k_{18}+k_{19}).
\end{equation}
By noticing that $k_{18}$, $k_{19}$, $\delta$ are positive which yields that $g_2(x):=x+(k_{18}+k_{19})e^{x\delta}$ is an increasing continuous function satisfying $g_2(+\infty)=+\infty$ and $g_2(-\infty)=-\infty$. Therefore, for given $\rho$ and $-2\rho_2-k_{16}l_{14}^{-1}-k_{17}l_{15}^{-1}-(k_{18}+k_{19})$, there exists a unique  $\tilde\rho_2$ such that \eqref{discounting coefficient equation 2} holds by intermediate value theorem.

Integrating from 0 to $T$ on both sides of \eqref{eq35} and using
$$\frac{1-e^{-\tilde\rho_2s}}{\tilde\rho_2}\leq \frac{1-e^{-\tilde\rho_2T}}{\tilde\rho_2},$$
for all $s\in[0,T]$, we have \eqref{eq366}. Letting $t=0$ in \eqref{eq35} and noticing $1\wedge e^{-\tilde\rho_2T}\leq e^{-\tilde\rho_2s}\leq 1\vee e^{-\tilde\rho_2T}$ for all $s\in[0,T]$, we obtain \eqref{eq377}. On the other hand, if $\tilde\rho_2>0$, letting $t=0$ in \eqref{eq34}, \eqref{eq388} is derived.
\end{proof}
\vspace{3.8cm}

\textbf{Appendix C. Proof of Lemma \ref{lemma5.2}.}\\
\\
\begin{proof}
For $1\leq k\leq K$, denote the $k$-type agent state average by
$ \widetilde x^{(k)}:=\frac{1}{N_k}\sum_{j\in\mathcal I_k}\widetilde x_j,$ thus
\begin{equation}\nonumber\left\{\begin{aligned}
&d\widetilde{x}^{(k)}(t)=\Bigg[A_k\widetilde{x}^{(k)}+\hat A_k\widetilde x^{(k)}(t-\delta)+\tilde A\widetilde x^{(N)}(t-\delta)+\frac{1}{N_k}\sum_{j\in\mathcal I_k}B\widetilde u_j+ \frac{1}{N_k}\sum_{j\in\mathcal I_k}\hat B\widetilde u_j(t-\theta)\\
&\qquad\qquad\qquad+\tilde B\widetilde u^{(N)}(t-\theta) \Bigg]dt+\frac{1}{N_k}\sum_{j\in\mathcal I_k}\Big[D\widetilde u_j+\hat D \widetilde u_j(t-\theta)\Big]dW_j(t),\ t\in[0,T],\\
&\widetilde {x}^{(k)}(0)=\frac{1}{N_k}\sum_{j\in\mathcal I_k}\xi_j,\ \widetilde {x}^{(k)}(t)=x^0(t),\ t\in[-\delta,0).
\end{aligned}\right.\end{equation}
Noticing that
\begin{equation}\nonumber\left\{\begin{aligned}
d\mathbb E\alpha_k(t)=\ &\Bigg[A_k\mathbb E\alpha_k+\hat A_k\mathbb E\alpha_k(t-\delta)+\tilde{A}\sum_{l=1}^K\pi_l\mathbb E\alpha_l(t-\delta)-\mathbb E\big(\mathbb{B}_1^k \mathbf p_k+\mathbb{B}_5^k \mathbf p_k(t+\theta)I_{[0,T-\theta]}\\
&+\mathbb{D}_1^k \mathbf q_k+\mathbb{D}_2^k \mathbf q_k(t+\theta)I_{[0,T-\theta]}+\Theta_3\big)-\mathbb E\big(\mathbb{B}_6^k \mathbf p_k(t-\theta)+\mathbb{B}_7^k \mathbf p_k+\mathbb{D}_9^k \mathbf q_k(t-\theta)\\
&+\mathbb{D}_{10}^k \mathbf q_k+\Theta_3(t-\theta)\big)I_{[\theta,T]}+\hat Bu^0I_{[0,\theta]}-\sum_{l=1}^K\pi_l \mathbb E\big(\mathbb{B}_2^l \mathbf p_l(t-\theta)+\mathbb{B}_8^l \mathbf p_l+\mathbb{D}_3^l \mathbf q_l(t-\theta)\\
&+\mathbb{D}_4^l \mathbf q_l\big)I_{[\theta,T]}+\tilde Bu^0I_{[0,\theta]} \Bigg]dt, \\
\mathbb E\alpha_k(0)=\ &\mathbb E\xi^{(k)},\ \mathbb E\alpha_k(t)=x^0(t),\ t\in[-\delta,0).
\end{aligned}\right.\end{equation}
Under (A2), for $1\leq k\leq K$, $\{\xi_j,j\in\mathcal I_k\}$ are independent identically distributed. Note that $(p_j(\cdot),q_j(\cdot))\in \mathcal F_t^j$,  thus $\{(p_j,q_j),j\in\mathcal I_k\}$ are independent identically distributed. Using Cauchy-Schwartz inequality, Burkholder-Davis-Gundy inequality and estimations of SDE, then we derive
\begin{equation*}\begin{aligned}
&\mathbb E\sup_{0\leq s\leq t}\Big|\widetilde x^{(k)}-\mathbb E\alpha_k\Big|^2\\
\leq\ &\mathbb E\left|\frac{1}{N_k}\sum_{j\in\mathcal I_k}\left(\xi_j-\mathbb E\xi^{(k)}\right)\right|^2+C\mathbb E\int_0^t\left[\left|\widetilde x^{(k)}-\mathbb E\alpha_k\right|^2+\left|(\widetilde x^{(k)}-\mathbb E\alpha_k)(s-\delta)\right|^2\right.\\
&+\left.\left|(\widetilde x^{(N)}-\sum_{l=1}^K\pi_l\mathbb E\alpha_l)(s-\delta)\right|^2\right]ds\\
&+C\mathbb E\int_0^t\left|\frac{1}{N_k}\sum_{j\in\mathcal I_k}\Big(  p_j-\mathbb E \mathbf p_k+p_j(s+\theta)-\mathbb E\mathbf p_k(s+\theta)+q_j-\mathbb E\mathbf q_k+ q_j(s+\theta)-\mathbb E\mathbf q_k(s+\theta)\Big)\right|^2ds\\
\end{aligned}\end{equation*}\
\begin{equation*}\begin{aligned}
&+C\mathbb E\int_0^t\left|\frac{1}{N_k}\sum_{j\in\mathcal I_k}\Big( p_j(s-\theta)-\mathbb E\mathbf p_k(s-\theta)+  p_j-\mathbb E\mathbf p_k+ q_j(s-\theta)-\mathbb E\mathbf q_k(s-\theta)+  q_j-\mathbb E\mathbf q_k\Big)\right|^2ds\\
&+C\mathbb E\int_0^t\frac{1}{N_k^2}\sum_{j\in\mathcal I_k}\left|D\widetilde u_j+\hat D \widetilde u_j(s-\theta)\right|^2ds\\
\leq\ &C\mathbb E\int_0^t\left|\widetilde x^{(k)}-\mathbb E\alpha_k\right|^2ds+C\mathbb E\int_0^t\left|\widetilde x^{(N)}-\sum_{l=1}^K\pi_l\mathbb E\alpha_l\right|^2ds+\frac{C}{N_k} .
\end{aligned}\end{equation*}
Here, for the delay items, we have used the facts similar to
\begin{equation*}\begin{aligned}
\mathbb E\int_0^t\left|(\widetilde x^{(k)}-\mathbb E\alpha_k)(s-\delta)\right|^2ds&=\mathbb E\int_{-\delta}^{t-\delta}\left|\widetilde x^{(k)}-\mathbb E\alpha_k\right|^2ds=\mathbb E\int_{0}^{t-\delta}\left|\widetilde x^{(k)}-\mathbb E\alpha_k\right|^2ds\\
&\leq \mathbb E\int_{0}^{t}\left|\widetilde x^{(k)}-\mathbb E\alpha_k\right|^2ds.
\end{aligned}\end{equation*}
Gronwall inequality implies that
\begin{equation*}\begin{aligned}
\mathbb E\sup_{0\leq s\leq t}\Big|\widetilde x^{(k)}-\mathbb E\alpha_k\Big|^2
\leq C\mathbb E\int_0^t\left|\widetilde x^{(N)}-\sum_{l=1}^K\pi_l\mathbb E\alpha_l\right|^2ds+\frac{C}{N_k}.
\end{aligned}\end{equation*}
Since
\begin{equation*}\begin{aligned}
\widetilde x^{(N)}-\sum_{l=1}^K\pi_l\mathbb E\alpha_l=&\sum_{l=1}^K\Big(\pi_l^{(N)}\widetilde x^{(l)}-\pi_l\mathbb E\alpha_l\Big)=
\sum_{l=1}^K\pi_l^{(N)}\Big(\widetilde x^{(l)}-\mathbb E\alpha_l\Big)+\sum_{l=1}^K\Big(\pi_l^{(N)}-\pi_l\Big)\mathbb E\alpha_l,
\end{aligned}\end{equation*}
we get
\begin{equation*}\begin{aligned}
&\mathbb E\sup_{0\leq s\leq t}\left|\widetilde x^{(N)}-\sum_{l=1}^K\pi_l\mathbb E\alpha_l\right|^2\\
\leq&\
C\sum_{l=1}^K\mathbb E\sup_{0\leq s\leq t}\Big|\widetilde x^{(l)}-\mathbb E\alpha_l\Big|^2+C\epsilon_N^2\leq C\mathbb E\int_0^t\left|\widetilde x^{(N)}-\sum_{l=1}^K\pi_l\mathbb E\alpha_l\right|^2ds+\frac{C}{N}+C\epsilon_N^2.
\end{aligned}\end{equation*}
Therefore, the first result follows from Gronwall inequality.

Define $ \widetilde u^{(k)}:=\frac{1}{N_k}\sum_{j\in\mathcal I_k}\widetilde u_j,$ similarly we have
\begin{equation*}\begin{aligned}
&\mathbb E\sup_{0\leq s\leq t}\Big|\widetilde u^{(k)}-\mathbb Ev_k\Big|^2\\
\leq&C\mathbb E\sup_{0\leq s\leq t}\left|\frac{1}{N_k}\sum_{j\in\mathcal I_k}\Big(  p_j-\mathbb E \mathbf p_k+ p_j(s+\theta)-\mathbb E\mathbf p_k(s+\theta)+ q_j-\mathbb E\mathbf q_k+ q_j(s+\theta)-\mathbb E\mathbf q_k(s+\theta)\Big)\right|^2\\
\leq&\frac{C}{N_k} .
\end{aligned}\end{equation*}
Since
\begin{equation*}\begin{aligned}
\widetilde u^{(N)}-\sum_{l=1}^K\pi_l\mathbb Ev_l=&\sum_{l=1}^K\Big(\pi_l^{(N)}\widetilde u^{(l)}-\pi_l\mathbb Ev_l\Big)=
\sum_{l=1}^K\pi_l^{(N)}\Big(\widetilde u^{(l)}-\mathbb Ev_l\Big)+\sum_{l=1}^K\Big(\pi_l^{(N)}-\pi_l\Big)\mathbb Ev_l,
\end{aligned}\end{equation*}
we get
\begin{equation*}\begin{aligned}
&\mathbb E\sup_{0\leq s\leq t}\left|\widetilde u^{(N)}-\sum_{l=1}^K\pi_l\mathbb Ev_l\right|^2
\leq
C\sum_{l=1}^K\mathbb E\sup_{0\leq s\leq t}\Big|\widetilde u^{(l)}-\mathbb Ev_l\Big|^2+C\epsilon_N^2\leq \frac{C}{N}+C\epsilon_N^2.
\end{aligned}\end{equation*}
\end{proof}
\vspace{2cm}
\textbf{Appendix D. Proof of Lemma \ref{Lemma5}.}\\
\\
\begin{proof}
According to \eqref{eq6}-\eqref{eq8}, it yields
\begin{equation*}\begin{aligned}
\mathbb E\sup_{0\leq s\leq t}|\delta x_i|^2 & \leq C\Bigg(\mathbb E\int_0^t(|\delta u_i|^2+|\delta u_i(s-\theta)|^2)ds\Bigg)+C\mathbb E\int_0^t(|\delta x_i|^2+|\delta x_i(s-\delta)|^2)ds\\
&\quad +C\mathbb E\int_0^t|\delta x^{(N)}(s-\delta)|^2ds\\
&\leq C\mathbb E\int_0^t|\delta u_i|^2ds+C\mathbb E\int_0^t|\delta x_i|^2ds+C\mathbb E\int_0^t|\delta x^{(N)}|^2ds,
\end{aligned}\end{equation*}
for $j\neq i$,
\begin{equation*}\begin{aligned}
\mathbb E\sup_{0\leq s\leq t}|\delta x_j|^2 &\leq \frac{C}{N^2}\mathbb E\int_0^t|\delta u_i(s-\theta)|^2ds+C\mathbb E\int_0^t(|\delta x_j|^2+|\delta x_j(s-\delta)|^2)ds+C\mathbb E\int_0^t\left|\delta x^{(N)}(s-\delta)\right|^2ds\\
&\leq \frac{C}{N^2}\mathbb E\int_0^t|\delta u_i|^2ds+C\mathbb E\int_0^t|\delta x_j|^2ds+C\mathbb E\int_0^t\left|\delta x^{(N)}\right|^2ds,
\end{aligned}\end{equation*}
and
\begin{equation*}\begin{aligned}
&\mathbb E\sup_{0\leq s\leq t}|\delta x_{(k)}|^2\leq C\mathbb E\int_0^t|\delta u_i|^2ds+C\mathbb E\int_0^t|\delta x_{(k)}|^2ds+CN^2\mathbb E\int_0^t|\delta x^{(N)}|^2ds.
\end{aligned}\end{equation*}
Noticing that $$\delta x^{(N)}=\frac{1}{N}\delta x_i+\frac{1}{N}\sum_{l=1}^K\delta x_{(l)},$$
we arrive at
\begin{equation*}\begin{aligned}
&\mathbb E\sup_{0\leq s\leq t}|\delta x_i|^2\leq C\mathbb E\int_0^T|\delta u_i|^2ds+C\mathbb E\int_0^t|\delta x_i|^2ds+\frac{C}{N^2}\sum_{l=1}^K\mathbb E\int_0^t\left|\delta x_{(l)}\right|^2ds,
\end{aligned}\end{equation*}
and
\begin{equation*}\begin{aligned}
&\mathbb E\sup_{0\leq s\leq t}\left|\delta x_{(k)}\right|^2\leq C\mathbb E\int_0^T|\delta u_i|^2ds+ C\mathbb E\int_0^t|\delta x_{(k)}|^2ds+C\mathbb E\int_0^t|\delta x_i|^2ds+
 C\sum_{l=1}^K\mathbb E\int_0^t\left|\delta x_{(l)}\right|^2ds.
\end{aligned}\end{equation*}
Therefore, it follows from Gronwall inequality that
\begin{equation*}\begin{aligned}
\mathbb E\sup_{0\leq s\leq t}|\delta x_i|^2+\sum_{l=1}^K\mathbb E\sup_{0\leq s\leq t}\left|\delta x_{(l)}\right|^2\leq C\mathbb E\int_0^T|\delta u_i|^2ds.
\end{aligned}\end{equation*}
Thus,
\begin{equation*}
\mathbb E\sup_{0\leq s\leq t}\left|\delta x^{(N)}\right|^2\leq\frac{C}{N^2}\mathbb E\int_0^T|\delta u_i|^2ds.
\end{equation*}
By Gronwall inequality again, we obtain \eqref{estimation-5}.
\end{proof}
\vspace{4.5cm}
\textbf{Appendix E. Proof of Lemma \ref{lemma5.4}.}\\
\\
\begin{proof}
First,
\begin{equation*}\left\{\begin{aligned}
&d(x_k^{**}-\delta x_{(k)})(t)=\Bigg[A_k(x_k^{**}-\delta x_{(k)})+\hat A_k(x_k^{**}-\delta x_{(k)})(t-\delta)+\tilde A\Bigg(\pi_k-\frac{N_k-I_{\mathcal I_k}(i)}{N}\Bigg)\delta x_i(t-\delta)\\
&\qquad\qquad\qquad\qquad+\tilde A\pi_k\sum_{l=1}^K(x_l^{**}-\delta x_{(l)})(t-\delta)+\tilde A\Bigg(\pi_k-\frac{N_k-I_{\mathcal I_k}(i)}{N}\Bigg)\sum_{l=1}^K\delta x_{(l)}(t-\delta)\\
&\qquad\qquad\qquad\qquad+\tilde B\Bigg(\pi_k-\frac{N_k-I_{\mathcal I_k}(i)}{N}\Bigg)\delta u_i(t-\theta)\Bigg]dt,\\
&(x_k^{**}-\delta x_{(k)})(0)=0,\quad (x_k^{**}-\delta x_{(k)})(t)=0,\ t\in[-\delta,0),
\end{aligned}\right.\end{equation*}
and for $j\in\mathcal I_k,\ j\neq i$,
\begin{equation*}\left\{\begin{aligned}
&d(x_j^*-N_k\delta x_j)=\Bigg[A_{k}(x_j^*-N_k\delta x_j)+\hat A_{k}(x_j^*-N_k\delta x_j)(t-\delta)+\tilde A\Big(\pi_k-\pi_k^{(N)}\Big)\delta x_i(t-\delta)\\
&\qquad\qquad\qquad\qquad+\tilde A\Big(\pi_k-\pi_k^{(N)}\Big)\sum_{l=1}^Kx_l^{**}(t-\delta)+\tilde A\pi_k^{(N)}\sum_{l=1}^K(x_l^{**}-\delta x_{(l)})(t-\delta)\\
&\qquad\qquad\qquad\qquad+\tilde B\Big(\pi_k-\pi_k^{(N)}\Big)\delta u_i(t-\theta)\Bigg]dt,\\
&(x_j^*-N_k\delta x_j)(0)=0,\quad (x_j^*-N_k\delta x_j)(t)=0,\ t\in[-\delta,0).
\end{aligned}\right.\end{equation*}
Therefore, it follows from Burkholder-Davis-Gundy inequality that
\begin{equation*}\begin{aligned}
&\mathbb E\sup_{0\leq s\leq t}|x_k^{**}-\delta x_{(k)}|^2\\
 \leq\ &C\mathbb E\int_0^t(|x_k^{**}-\delta x_{(k)}|^2+|(x_k^{**}-\delta x_{(k)})(s-\delta)|^2)ds+C\sum_{l=1}^K\mathbb E\int_0^t|(x_l^{**}-\delta x_{(l)})(s-\delta)|^2ds\\
 &+\Big(\frac{C}{N^2}+C\epsilon_N^2\Big)\mathbb E\int_0^t\left[|\delta x_i(s-\delta)|^2+\sum_{l=1}^K|\delta x_{(l)}(s-\delta)|^2+|\delta u_i(s-\theta)|^2\right]ds\\
\leq \ & C\mathbb E\int_0^t|x_k^{**}-\delta x_{(k)}|^2ds+C\sum_{l=1}^K\mathbb E\int_0^t|x_l^{**}-\delta x_{(l)}|^2ds+\Big(\frac{C}{N^2}+C\epsilon_N^2\Big)  \mathbb E\int_0^T|\delta u_i|^2ds.
\end{aligned}\end{equation*}
Thus,
\begin{equation*}\begin{aligned}
\sum_{l=1}^K\mathbb E\sup_{0\leq s\leq t}|x_l^{**}-\delta x_{(l)}|^2
\leq& C\sum_{l=1}^K\mathbb E\int_0^t|x_l^{**}-\delta x_{(l)}|^2+\Big(\frac{C}{N^2}+C\epsilon_N^2\Big) \mathbb E\int_0^T|\delta u_i|^2ds.
\end{aligned}\end{equation*}
It then follows from Gronwall inequality that
$$\sum_{l=1}^K\mathbb E\sup_{0\leq t\leq T}|x_l^{**}(t)-\delta x_{(l)}(t)|^2\leq\Big(\frac{C}{N^2}+C\epsilon_N^2\Big) \mathbb E\int_0^T|\delta u_i|^2ds.$$
Similarly,
\begin{equation*}\begin{aligned}
\mathbb E\sup_{0\leq s\leq t}|x_j^*-N_k\delta x_j|^2 \leq& C\mathbb E\int_0^t|x_j^*-N_k\delta x_j|^2ds+C\sum_{l=1}^K\mathbb E\int_0^t|x_l^{**}-\delta x_{(l)}|^2ds\\
&\ +C\epsilon_N^2\mathbb E\int_0^t\left[|\delta x_i|^2+\sum_{l=1}^K|x_l^{**}|^2+|\delta u_i|^2\right]ds.
\end{aligned}\end{equation*}
By the first equation of \eqref{eq11}, we derive
\begin{equation*}\begin{aligned}
&\mathbb E\sup_{0\leq s\leq t}|x_k^{**}|^2 \leq\ C\mathbb E\int_0^t|x_k^{**}|^2ds+C\sum_{l=1}^K\mathbb E\int_0^t|x_l^{**}|^2ds+C\mathbb E\int_0^t|\delta x_i|^2ds+C\mathbb E\int_0^t|\delta u_i|^2ds.
\end{aligned}\end{equation*}
By virtue of Gronwall inequality, we have
\begin{equation*}\begin{aligned}
&\sum_{l=1}^K\mathbb E\sup_{0\leq s\leq t}|x_l^{**}|^2 \leq C\mathbb E\int_0^t|\delta x_i|^2ds+C\mathbb E\int_0^t|\delta u_i|^2ds.
\end{aligned}\end{equation*}
Then by noticing \eqref{estimation-6} and Lemma \ref{Lemma5}, we get
\begin{equation*}\begin{aligned}
&\qquad\mathbb E\sup_{0\leq s\leq t}|x_j^*-N_k\delta x_j|^2\\
&\leq C\mathbb E\int_0^t|x_j^*-N_k\delta x_j|^2ds+C\sum_{l=1}^K\mathbb E\int_0^t|x_l^{**}-\delta x_{(l)}|^2ds+C\epsilon_N^2\mathbb E\int_0^t|\delta x_i|^2ds+C\epsilon_N^2\mathbb E\int_0^t|\delta u_i|^2ds\\
&\leq C\mathbb E\int_0^t|x_j^*-N_k\delta x_j|^2ds+\Big(\frac{C}{N^2}+C\epsilon_N^2\Big) \mathbb E\int_0^T|\delta u_i|^2ds,
\end{aligned}\end{equation*}
which implies \eqref{estimation-7}.
\end{proof}
\vspace{1cm}
\textbf{Appendix F. Proof of Lemma \ref{Lemma7}.}\\
\\
\begin{proof}
It follows from \eqref{eq17} that
\small\begin{equation}\nonumber\left\{\begin{aligned}
&d\Bigg(\frac{1}{N_k}\sum_{j\in\mathcal I_k,j\neq i}y_1^j\Bigg)=-\Bigg[ \frac{A_k^\top}{N_k}\sum_{j\in\mathcal I_k,j\neq i}y_1^j+\frac{\hat A_k^\top}{N_k}\sum_{j\in\mathcal I_k,j\neq i}\mathbb E^{\mathcal F_t}[y_1^j(t+\delta)]I_{[0,T-\delta]}+\frac{Q+\tilde Q(t+\delta)}{N_k}\sum_{j\in\mathcal I_k,j\neq i}\widetilde x_j\Bigg] dt\\
&\qquad\qquad\qquad\qquad\qquad+\frac{1}{N_k}\sum_{j\in\mathcal I_k,j\neq i}z_1^jdW_j(t)+\frac{1}{N_k}\sum_{j\in\mathcal I_k,j\neq i}\sum_{l=1,l\neq j}^N z_1^{jl}dW_l(t),\\
&\frac{1}{N_k}\sum_{j\in\mathcal I_k,j\neq i}y_1^j(T)=\frac{G}{N_k}\sum_{j\in\mathcal I_k,j\neq i}\widetilde x_j(T)%,\ \frac{1}{N_k}\sum_{j\in\mathcal I_k,j\neq i}y_1^j(t)=0,\ t\in(T,T+\delta]
.
\end{aligned}\right.\end{equation}\normalsize
By the definition of $\mathbf{y}_k$, we have

\begin{equation}\nonumber\left\{\begin{aligned}
&d\mathbb E\mathbf{y}_k=-\Big[A_k^\top \mathbb E\mathbf{y}_k+\tilde A_k^\top(t+\delta) \mathbb E\mathbf{y}_k(t+\delta)I_{[0,T-\delta]}+(Q+\tilde Q(t+\delta))\mathbb E\widetilde x_k\Big] dt,\\
&\mathbb E\mathbf{y}_k(T)= G \mathbb{E}\widetilde x_k(T),%\ \mathbb E\mathbf{y}_k(t)=0,\ t\in(T,T+\delta],
\ k=1,\cdots,K,
\end{aligned}\right.\end{equation}
where $\widetilde x_k$ denotes the optimal state of $k-$type corresponding to \eqref{eq36} and $\mathbb E\widetilde{x}_k$ satisfies
\begin{equation}\nonumber\left\{\begin{aligned}
d\mathbb E\widetilde{x}_k=\ &\Bigg[A_k\mathbb E\widetilde{x}_k+\hat A_k\mathbb E\widetilde{x}_k(t-\delta)+\tilde{A}\mathbb E\widetilde{x}^{(N)}(t-\delta)-\mathbb E\big(\mathbb{B}_1^k \mathbf p_k+\mathbb{B}_5^k \mathbf p_k(t+\theta)I_{[0,T-\theta]}+\mathbb{D}_1^k \mathbf q_k\\
&\quad+\mathbb{D}_2^k \mathbf q_k(t+\theta)I_{[0,T-\theta]}+\Theta_3\big)-\mathbb E\big(\mathbb{B}_6^k \mathbf p_k(t-\theta)+\mathbb{B}_7^k \mathbf p_k+\mathbb{D}_9^k \mathbf q_k(t-\theta)+\mathbb{D}_{10}^k \mathbf q_k\\
&\quad+\Theta_3(t-\theta)\big)I_{[\theta,T]}+\hat Bu^0I_{[0,\theta]}+\tilde{B}\mathbb E\widetilde{u}^{(N)}(t-\theta) \Bigg]dt, \\
\mathbb E\widetilde{x}_k(0)=\ &\mathbb E\xi^{(k)},\ \mathbb E\widetilde{x}_k(t)=x^0(t),\ t\in[-\delta,0).
\end{aligned}\right.\end{equation}
Recall notation $\widetilde x^{(k)}$ defined in proof of Lemma \ref{lemma5.2}. Noticing
$$\frac{1}{N_k}\sum_{j\in\mathcal I_k,j\neq i}\widetilde x_j=\widetilde x^{(k)}-\frac{I_{\mathcal I_k}(i)}{N_k}\widetilde x_i,$$
we have
\begin{equation*}\begin{aligned}
&\mathbb E\sup_{t\leq s\leq T}\left|\frac{1}{N_k}\sum_{j\in\mathcal I_k,j\neq i}y_1^j-\mathbb E\mathbf{y}_k\right|^2\leq C\mathbb E\left|\frac{1}{N_k}\sum_{j\in\mathcal I_k,j\neq i}\widetilde x_j(T)-\mathbb{E}\widetilde x_k(T)\right|^2\\
&\qquad\qquad+C\mathbb E\int_t^T\left|\frac{1}{N_k}\sum_{j\in\mathcal I_k,j\neq i}y_1^j-\mathbb E\mathbf{y}_k\right|^2ds+C\mathbb E\int_t^T\left|\frac{1}{N_k}\sum_{j\in\mathcal I_k,j\neq i}\mathbb E^{\mathcal F_s}[y_1^j(s+\delta)]I_{[0,T-\delta]}\right.\\
&\qquad\qquad\left.-\mathbb E\mathbf{y}_k(s+\delta)I_{[0,T-\delta]}\right|^2ds+C\mathbb E\int_t^T\left|\frac{1}{N_k}\sum_{j\in\mathcal I_k,j\neq i}\widetilde x_j-\mathbb E\widetilde x_k\right|^2ds\\
\end{aligned}\end{equation*}
\begin{equation*}\begin{aligned}
&\leq C\mathbb E\int_t^T\left|\frac{1}{N_k}\sum_{j\in\mathcal I_k,j\neq i}y_1^j-\mathbb E\mathbf{y}_k\right|^2ds+C\mathbb E\left(|\widetilde x^{(k)}(T)-\mathbb E\alpha_k(T)|^2+|\mathbb E\widetilde x_k(T)-\mathbb E\alpha_k(T)|^2\right.\\
&\qquad\qquad\left.+\frac{1}{N_k^2}|\widetilde x_i(T)|^2\right)+C\mathbb E\int_t^T\left(|\widetilde x^{(k)}-\mathbb E\alpha_k|^2+|\mathbb E\widetilde x_k-\mathbb E\alpha_k|^2+\frac{1}{N_k^2}|\widetilde x_i|^2\right)ds.
\end{aligned}\end{equation*}
Here we have used
\begin{equation*}\begin{aligned}
&\quad\mathbb E\int_t^T\left|\frac{1}{N_k}\sum_{j\in\mathcal I_k,j\neq i}\mathbb E^{\mathcal F_s}[y_1^j(s+\delta)]I_{[0,T-\delta]}-\mathbb E\mathbf{y}_k(s+\delta)I_{[0,T-\delta]}\right|^2ds\\
&\leq \mathbb E\int_t^T\mathbb E^{\mathcal F_s}\left[\left|\frac{1}{N_k}\sum_{j\in\mathcal I_k,j\neq i}y_1^j(s+\delta)-\mathbb E\mathbf{y}_k(s+\delta)\right|^2I_{[0,T-\delta]}\right]ds\\
& = \mathbb E\int_t^T\left|\frac{1}{N_k}\sum_{j\in\mathcal I_k,j\neq i}y_1^j(s+\delta)-\mathbb E\mathbf{y}_k(s+\delta)\right|^2I_{[0,T-\delta]}ds\\
&=\mathbb E\int_{t+\delta}^{T+\delta}\left|\frac{1}{N_k}\sum_{j\in\mathcal I_k,j\neq i}y_1^j-\mathbb E\mathbf{y}_k\right|^2I_{[\delta,T]}ds\\
&\leq \mathbb E\int_{t}^{T}\left|\frac{1}{N_k}\sum_{j\in\mathcal I_k,j\neq i}y_1^j-\mathbb E\mathbf{y}_k\right|^2ds.
\end{aligned}\end{equation*}
Noticing
\begin{equation*}\begin{aligned}
&\sup_{0\leq s\leq t}|\mathbb E\widetilde x_k-\mathbb E\alpha_k|^2 \leq C\mathbb E\int_0^t\left|\widetilde x^{(N)}-\sum_{l=1}^K\pi_l\mathbb E\alpha_l\right|^2ds+C\mathbb E\int_0^t\left|\widetilde u^{(N)}-\sum_{l=1}^K\pi_l\mathbb Ev_l\right|^2ds.
\end{aligned}\end{equation*}
With the help of the proof of Lemma \ref{lemma5.2}, one gets
\begin{equation*}\begin{aligned}
&\mathbb E\sup_{0\leq t\leq T}\Big|\widetilde x^{(k)}-\mathbb E\alpha_k\Big|^2\leq \frac{C}{N}+C\epsilon_N^2.
\end{aligned}\end{equation*}
Then \eqref{estimation-9} is obtained based on above inequalities, $L^2$ boundness of $\mathbf{y}_k$, Lemma \ref{5.1}, Lemma \ref{lemma5.2} and Gronwall inequality.
\end{proof}
%\vspace{0.5cm}
%\textbf{Appendix H. Proof of Theorem \ref{asymptotic optimal}.}\\
%\\

\end{document}